# PROJECTIVE GROUP STRUCTURES
# AS ABSOLUTE GALOIS STRUCTURES
# WITH BLOCK APPROXIMATION*

by

Dan Haran**

*School of Mathematics, Tel Aviv University*

*Ramat Aviv, Tel Aviv 69978, Israel*

*e-mail: haran@post.tau.ac.il*

and

Moshe Jarden**

*School of Mathematics, Tel Aviv University*

*Ramat Aviv, Tel Aviv 69978, Israel*

*e-mail: jardenm@post.tau.ac.il*

and

Florian Pop

*Mathematisches Institut, Universität Bonn*

*Bering Str. 6, 53115 Bonn, Germany*

*e-mail: pop@math.uni-bonn.de*

MR Classification: 12E30

Directory: \Jarden\Diary\HJPa

24 September, 2002

---

\* Research supported by the Minkowski Center for Geometry at Tel Aviv University, established by the Minerva Foundation

\*\* Research partially done at the Max-Planck-Institute for Mathematic in Bonn.

# Table of Contents

**Introduction**

A. BACKGROUND AND MOTIVATION. One of the main features of Field Arithmetic is the interplay between the arithmetic-geometrical properties of a field and the profinite group theoretic properties of its absolute Galois group. Here is the prototype for this kind of results:

BASIC THEOREM:

(a) *If a field $K$ is PAC, then $\mathrm{Gal}(K)$ is projective (Ax, [FrJ, Thm. 10.17]).*

(b) *For every projective group $G$ there exists a field $K$ with $\mathrm{Gal}(K) \cong G$ (Lubotzky-v.d.Dries [FrJ, Cor. 20.16]).*

Here we say that a field $K$ is **PAC** if every absolutely irreducible variety $V$ over $K$ has a $K$-rational point. By an **absolutely irreducible variety** over $K$ we mean a geometrically integral scheme of finite type over $K$. We denote the separable closure of $K$ by $K_s$ and its algebraic closure by $\tilde{K}$. Then we call $\mathrm{Gal}(K) = \mathrm{Gal}(K_s/K)$ the **absolute Galois group** of $K$.

A profinite group $G$ is **projective** if every finite embedding problem

$$(1) \qquad\qquad (\varphi\colon G \to A,\ \alpha\colon B \to A)$$

for $G$ is solvable. Here $A$ and $B$ are finite groups, $\varphi$ is a homomorphism, and $\alpha$ is an epimorphism. A **solution** of (1) is a homomorphism $\gamma\colon G \to B$ with $\alpha \circ \gamma = \varphi$.

Both concepts "projective group" and "PAC field" have relative counterparts which we now describe.

Let $G$ be a profinite group and $\mathcal{G}$ a collection of closed subgroups of $G$. Call $G$ $\mathcal{G}$-**projective** if every finite embedding problem (1) for $G$ has a solution provided for each $\Gamma \in \mathcal{G}$ there is a homomorphism $\gamma_\Gamma\colon \Gamma \to B$ with $\alpha \circ \gamma_\Gamma = \varphi|_\Gamma$.

Let $K$ be a field and $\mathcal{K}$ a collection of separable algebraic extensions of $K$. Call $K$ P$\mathcal{K}$C (**pseudo $\mathcal{K}$-closed**) if every smooth absolutely irreducible variety $V$ over $K$ with a $K'$-rational point for each $K' \in \mathcal{K}$ has a $K$-rational point.

In both cases we have local-global principles. Thus, $G$ is $\mathcal{G}$-projective if the existence of local solutions of embedding problems guaranties the existence of global so-lutions. Analogously, $K$ is P$\mathcal{K}$C if the existence of local points on smooth absolutely irreducible varieties give global points on them.

It is desirable to generalize the Basic Theorem to the relative case:

TARGET:

(a) Let $K$ be a field and $\mathcal{K}$ a collection of separable algebraic extensions of $K$. Put $\mathcal{G} = \{\mathrm{Gal}(K') \mid K' \in \mathcal{K}\}$. Suppose $K$ is P$\mathcal{K}$C. Then $\mathrm{Gal}(K)$ is $\mathcal{G}$-projective.

(b) Let $G$ be a profinite group and $\mathcal{G}$ a collection of closed subgroups of $G$. Suppose $G$ is $\mathcal{G}$-projective and for each $\Gamma \in \mathcal{G}$ there exists a field $F_\Gamma$ with $\mathrm{Gal}(F_\Gamma) \cong \Gamma$. Then there exists a field $K$ and an isomorphism $\varphi\colon G \to \mathrm{Gal}(K)$. Moreover, for each $\Gamma \in \mathcal{G}$ let $K_\Gamma$ be the fixed field of $\varphi(\Gamma)$ in $K_s$. Put $\mathcal{K} = \{K_\Gamma \mid \Gamma \in \mathcal{G}\}$. Then $K$ is P$\mathcal{K}$C.

The Basic Theorem is a special case of the Target in which both $\mathcal{K}$ and $\mathcal{G}$ are empty.

Another special case of the Target occurs when $\mathcal{K}$ is the collection of all real closures of $K$ and $\mathcal{G}$ is the collection of all subgroups of $G$ which are isomorphic to $\mathrm{Gal}(\mathbb{R})$ [HaJ1, p. 450, Thm.]. In this case P$\mathcal{K}$C are referred to as **PRC fields**. However, in order for Part (b) of the Target to hold, we must assume 1 does not lie in the closure of $\mathcal{G}$; that is, $G$ has an open subgroup $U$ which contains no $\Gamma \in \mathcal{G}$.

Similarly, the Target is reached when $\mathcal{K}$ is the collection of all $p$-adic closures of $K$ for some fixed prime number $p$ and $\mathcal{G}$ is the collection of all subgroups of $G$ which are isomorphic to $\mathrm{Gal}(\mathbb{Q}_p)$ [HaJ2, p. 148, Thm.]. Again, we must assume 1 does not belong to the closure of $\mathcal{G}$. Then P$\mathcal{K}$C fields are just **P$p$C fields**.

Another instance where the Target is obtained is when $\mathcal{K} = \{K_1, \ldots, K_n\}$ and each $K_i$ is Henselian with respect to a valuation $v_i$ such that $v_1|_K, \ldots, v_n|_K$ are independent ([Koe, Thm. 2'] or [HaJ3, Theorems A and B]). Here one starts in Part B with a profinite group $G$ which is projective with respect to $n$ closed subgroups $G_1, \ldots, G_n$, each of which is isomorphic to the absolute Galois group of a field. Then one constructs $K$ and $\varphi$ such that the fixed field $K_i$ of $\varphi(G_i)$ is Henselian with respect to a valuation $v_i$, $i = 1, \ldots, n$. Moreover, the restrictions of $v_1, \ldots, v_n$ to $K$ are independent.

In general it is possible to prove Part (a) of the Target under some mild compactness assumption on $\mathcal{K}$ [Pop, Thm. 3.3]. We are therefore allowed to make the same assumption on $\mathcal{G}$ in Part (b) of the Target. Nevertheless, when we try to realize $G$ as an absolute Galois group, we are forced to solve certain infinite embedding problems and not only finite ones. So, we must assume $G$ is "strongly $\mathcal{G}$-projective" rather than only $\mathcal{G}$-projective. This has actually been done in [Pop, Thm. 3.4] (Note however that the adjective "strongly" is mistakenly ommited in the formulation of [Pop, Thm. 3.4]). But replacing "$G$ is $\mathcal{G}$-projective" by "$G$ is strongly $\mathcal{G}$-projective" in Part (b) brings the Target out of balance. In general we allow to add extra conditions to (a) and to (b). The rule is that each assumption we make on $\mathcal{K}$ in (a) should appear as a consequence in (b). Similarly, each assumption we make on $\mathcal{G}$ in (b) should appear as a consequence in (a). The disturbed balance in [Pop] is restored only when "large quotients' exist, as in the case of $p$-adically closed fields [Pop, Section 1, Lemma and Definition]. The general case is left unbalanced in [Pop].

The goal of this work is to achieve a very general balanced Target. Like in the above mentioned three instances, we extend both $\mathcal{K}$ and $\mathcal{G}$ to "structures" over a profinite space $X$ and let each field in $\mathcal{K}$ be a Henselian closure of a valuation of the base field $K$. In order to prove projectivity of the group structure in (a) we must assume a strong form of the weak approximation theorem. We call it the "block approximation condition". One of the main achievments of this work is the realization of the strucutre in (b) as an "absolute Galois strucutre" of a "field-valuation strucucture" satisfying the block approximation condition.

B. The main theorem. For the convenience of the reader we state the main result of this work, define all concepts appearing in it, and describe the most essential ingredients of the proof.

Main Theorem:

(a) Let $\mathbf{K} = (K, X, K_x, v_x)_{x \in X}$ be a proper field-valuation structure. Suppose $\mathbf{K}$ satisfies the block approximation condition. Then $\mathrm{Gal}(\mathbf{K}) = (\mathrm{Gal}(K), X, \mathrm{Gal}(K_x))_{x \in X}$ is a proper projective group structure.

(b) Let $\mathbf{G} = (G, X, G_x)_{x \in X}$ be a proper projective group structure and $\bar{\kappa} \colon \mathbf{G} \to \mathrm{Gal}(\bar{\mathbf{K}})$ be a Galois approximation of $\mathbf{G}$. Then there exists a proper field-valuation structure $\mathbf{K} = (K, X, K_x, v_x)_{x \in X}$ satisfying the block approximation condition and there is an isomorphism $\kappa \colon \mathbf{G} \to \mathrm{Gal}(\mathbf{K})$ which lifts $\bar{\kappa}$.

Here are the definitions of the notions which occur in the Main Theorem.

We call $\mathbf{G} = (G, X, G_x)_{x \in X}$ a **group structure** if $G$ is a profinite group, $X$ is a profinite space, and for each $x \in X$, $G_x$ is a closed subgroup of $G$ satisfying these conditions:

(2a) $G$ acts continuously on $X$ from the right.

(2b) $G_{x^g} = G_x^g$ for all $x \in X$ and $g \in G$.

(2c) Let $\mathrm{Subgr}(G)$ be the space of all closed subgroups of $G$ equipped with the **étale topology**. (A basis of the étale topology consists of all sets $\mathrm{Subgr}(U)$ with $U$ open in $G$.) Then the map $\delta_{\mathbf{G}} \colon X \to \mathrm{Subgr}(G)$ defined by $\delta_{\mathbf{G}}(x) = G_x$ is continuous in the étale topology.

(2d) $\{g \in G \mid x^g = x\} \le G_x$ for each $x \in X$.

We say $\mathbf{G}$ is **proper** if the map $\delta_{\mathbf{G}} \colon X \to \{G_x \mid x \in X\}$ is a homeomorphism in the étale topology.

A group structure $\mathbf{G}$ is **projective** if every finite embedding problem

(4) $$(\varphi \colon \mathbf{G} \to \mathbf{A}, \ \alpha \colon \mathbf{B} \to \mathbf{A})$$

for $\mathbf{G}$ is solvable. Here we call (4) an **embedding problem** if the following holds:

(5a) $\mathbf{A} = (A, I, A_i)_{i \in I}$ and $\mathbf{B} = (B, J, B_j)_{j \in J}$ are **finite group structures**, i.e., $A$, $B$, $I$, and $J$ are finite.

(5b) $\varphi \colon \mathbf{G} \to \mathbf{A}$ is a **morphism**; that is, $\varphi$ is a pair consisting of a homomorphism $\varphi \colon G \to A$ and a continuous map $\varphi \colon X \to I$ such that $\varphi(x^g) = \varphi(x)^{\varphi(g)}$ and $\varphi(G_x) \le A_{\varphi(x)}$ for all $x \in X$ and $g \in G$.

(5c) $\alpha \colon \mathbf{B} \to \mathbf{A}$ is a **cover**; that is, $\alpha$ is a morphism, $\alpha(B) = A$, $\alpha(J) = I$, $\alpha \colon B_j \to A_{\alpha(j)}$ is an isomorphism for each $j \in J$, and for all $j_1, j_2 \in J$ with $\varphi(j_1) = \varphi(j_2)$ there is $b \in \mathrm{Ker}(\alpha)$ with $j_1^b = j_2$.

A **solution** of (4) is a morphism $\gamma\colon \mathbf{G} \to \mathbf{B}$ satisfying $\alpha \circ \gamma = \varphi$.

We call $(K, X, K_x)_{x \in X}$ a **field structure** if $K$ is a field, $X$ is a profinite space, and $K_x$ is a separable algebraic extension, $x \in X$, such that

$$\mathrm{Gal}(\mathbf{K}) = (\mathrm{Gal}(K), X, \mathrm{Gal}(K_x))_{x \in X}$$

is a group structure.

A **Galois approximation** of a group structure $\mathbf{G} = (G, X, G_x)_{x \in X}$ is a morphism $\bar{\kappa}\colon \mathbf{G} \to \mathrm{Gal}(\bar{\mathbf{K}})$ where $\bar{\mathbf{K}} = (\bar{K}, \bar{X}, \bar{K}_{\bar{x}})_{\bar{x} \in \bar{X}}$ is a field structure, $\bar{\kappa}(G) = \mathrm{Gal}(\bar{K})$, $\bar{\kappa}(X) = \bar{X}$, and $\bar{\kappa}\colon G_x \to \mathrm{Gal}(\bar{K}_{\bar{\kappa}(x)})$ is an isomorphism for each $x \in X$.

We call $\mathbf{K} = (K, X, K_x, v_x)_{x \in X}$ a **field-valuation structure** if $(K, X, K_x)_{x \in X}$ is a field structure and $v_x$ is a valuation of $K_x$ satisfying these conditions:

(6a) $v_{x^\sigma} = v_x^\sigma$ for all $x \in X$ and $\sigma \in \mathrm{Gal}(K)$.

(6b) For each finite separable extension $L$ the map $\nu_L\colon X_L \to \mathrm{Val}(L)$ given by $\nu_L(x) = v_x|_L$ is continuous. Here $X_L = \{x \in X \mid L \subseteq K_x\}$ and $\mathrm{Val}(L)$ is the space of all valuation of $L$ including the trivial one. A subbasis for the topology of $\mathrm{Val}(L)$ is the collection of all sets

$$U = \{w \in \mathrm{Val}(L) \mid w(a) > 0\} \quad \text{and} \quad U' = \{w \in \mathrm{Val}(L) \mid w(a) \geq 0\}$$

with $a \in L$.

A **block approximation problem** for $\mathbf{K}$ is a data $(V, X_i, L_i, \mathbf{a}_i, c_i)_{i \in I_0}$ satisfying these conditions:

(7a) $I_0$ is a finite set.

(7b) $L_i$ is a finite separable extension of $K$ contained in $K_x$ for all $x \in X_i$ and $i \in I_0$.

(7c) $X_i$ is an open-closed subset of $X$, $i \in I_0$.

(7d) $\mathrm{Gal}(L_i) = \{\sigma \in \mathrm{Gal}(L) \mid X_i^\sigma = X_i\}$, $i \in I_0$.

(7e) For each $i \in I_0$ let $R_i$ be a subset of $\mathrm{Gal}(K)$ satisfying $\mathrm{Gal}(K) = \bigcup_{\rho \in R_i} \mathrm{Gal}(L_i)\rho$. Then $X = \bigcup_{i \in I_0} \bigcup_{\rho \in R_i} X_i^\rho$.

(7f) $V$ is a smooth absolutely irreducible variety over $K$.

(7g) $\mathbf{a}_i \in V(L_i)$, $i \in I_0$.

(7h) $c_i \in K^\times$, $i \in I_0$.

A **solution** of the block approximation problem is a point $\mathbf{a} \in V(K)$ satisfying $v_x(\mathbf{a} - \mathbf{a}_i) > v(c_i)$ for all $i \in I_0$ and $x \in X_i$. We say that $\mathbf{K}$ satisfies the **block approximation condition** if every block approximation problem for $\mathbf{K}$ has a solution.

Finally, in the notation of (b) of the Main Theorem, we say that $\kappa$ **lifts** $\bar{\kappa}$ if $K$ is a regular extension of $\bar{K}$ and the epimorphism res: $\mathrm{Gal}(K) \to \mathrm{Gal}(\bar{K})$ extends to a morphism $\rho\colon \mathrm{Gal}(\mathbf{K}) \to \mathrm{Gal}(\bar{\mathbf{K}})$ with $\rho \circ \kappa = \bar{\kappa}$.

In the rest of the introduction we explain some of the main points of the proof. This will partially explain why the notions in the Main Theorem are so involved.

In the proof of Part (b) of the Target we have to solve embedding problems of the type $(\varphi\colon G \to \mathrm{Gal}(K), \alpha\colon \mathrm{Gal}(L) \to \mathrm{Gal}(K))$. Since $\mathrm{Gal}(K)$ and $\mathrm{Gal}(L)$ are infinite, it does not follow immediately from the projectivity of $G$ that a solution $\gamma$ exists. However, a result of Gruenberg [FrJ, Lemma 20.8] does give $\gamma$ in the setup of the Basic Theorem. In all other cases of the Target Theorem proved prior to this work it is needed that for each $\Gamma \in \mathcal{G}$, $\gamma(\Gamma)$ belongs to a subset of $\mathrm{Subgr}(G)$ given in advance. Therefore, the profinite groups $G$ have been equipped with certain group structures and homomorphisms have been replaced by morphisms such that solvability of finite embedding problems in the so obtained category implies solvability of arbitrary embedding problems.

Each of these structures consisted of a profinite group $G$ acting on a profinite space and local objects parametrized by $X$. It was further assumed that the action of $G$ on $X$ is regular; that is $x^g = x$ for $x \in X$ and $g \in G$ implies $g = 1$. This gave a closed system of representatives for the $G$-orbits of $X$ [HaJ2, Lemma 2.4]. But in general, closed system of representatives do not exist. Instead we find representatives modulo each open normal subgroup of $G$. More precisely, let $\mathbf{G} = (G, X, G_x)_{x \in X}$ be a group structure as in the Main Theorem and $N$ an open normal subgroup of $G$. Then we find a finite system of triples $(G_i, X_i, R_i)_{i \in I_0}$ which we call a **special partition** of $G$. It satisfies the following conditions:

(8a)  $I_0$ is a finite set, disjoint from $X$.

(8b)  $X_i$ is an open-closed subset of $X$, $i \in I_0$.

(8c)  $G_i$ is an open subgroup of $G$ containing $G_x$ for all $x \in X_i$, $i \in I_0$.

(8d) $G_i = \{\sigma \in \mathrm{Gal}(L) \mid X_i^\sigma = X_i\}$, $i \in I_0$.

(8e) $R_i$ is finite, $G = \bigcup_{\rho \in R_i} G_i \rho$, and $X = \bigcup_{i \in I_0} \bigcup_{\rho \in R_i} X_i^\rho$.

The existence of special partitions goes back to [Pop, Prop. 4.9].

We use special partitions on several occasions:

(9a) to extend each homomorphism $\varphi\colon G \to A$ with a finite group $A$ to a morphism $\varphi\colon \mathbf{G} \to \mathbf{A}$ where $\mathbf{A} = (A, I, A_i)_{i \in I}$ is a finite group structure given in advance (Lemma 3.7);

(9b) in the definition of "unirational arithmetical problem" (Section 6) and "block approximation problem" (Section 12) and in the proof of Part (a) of the Main Theorem (Lemma 14.2); and

(9d) in the proof of Part (b) of the Main Theorem (Lemma 15.1).

A second essential ingredient in the proof of Part (a) of the Main Theorem is the local homeomorphism theorem for étale morphisms of varieties over Henselian fields [GPR, Thm. 9.4]. A special partition, a "locally uniform Hensel's lemma" (Corollary 10.4), and block approximation prepare the use of the local homeomorphism theorem. The idea to use this set up goes back to [HaJ3, Prop. 3.2]. Block approximation can be found in [FHV, Prop. 2.1] in the context of real closed fields.

In a subsequent work we intend to apply the Main Theorem to prove the Target Theorem in a general $p$-adic setting which will make a far reaching generalization of [HaJ1].

## 1. Étale Topology

Let $G$ be a profinite group. Denote the collection of all closed (resp. open, open normal) subgroups of $G$ by $\mathrm{Subgr}(G)$ (resp. $\mathrm{Open}(G)$, $\mathrm{OpenNormal}(G)$). We introduce two topologies on $\mathrm{Subgr}(G)$, the strict topology and the étale topology, and relate them to each other.

A basis of the **strict topology** is the collection of all sets

$$(1) \qquad \nu(H, N) = \{A \in \mathrm{Subgr}(G) \mid AN = HN\},$$

with $H \in \mathrm{Open}(G)$ and $N \in \mathrm{OpenNormal}(G)$. When $G$ is finite, the strict topology is the discrete topology. In general, $\mathrm{Subgr}(G) \cong \varprojlim \mathrm{Subgr}(G/N)$ with $N$ ranging over all open normal subgroups of $G$. Thus, $\mathrm{Subgr}(G)$ is a profinite space under the strict topology. Indeed, each of the sets $\nu(H, N)$ is also closed in the strict topology. We use the adverb "strictly" as a substitute for "in the strict topology". For example, given a subset $\mathcal{G}$ of $\mathrm{Subgr}(G)$, we say $\mathcal{G}$ is **strictly open** (resp. **closed**, **compact**, **Hausdorff**) if it is open (resp. closed, compact, Hausdorff) in the strict topology. Likewise, for a function $f$ from a topological space $X$ into $\mathrm{Subgr}(G)$ we say $f$ is **strictly continuous** if $f$ is continuous when $\mathrm{Subgr}(G)$ is equipped with the strict topology.

A basis of the **étale topology** is the collection of all sets

$$\{\mathrm{Subgr}(U) \mid U \in \mathrm{Open}(G)\}$$

with $U \in \mathrm{Open}(G)$. As above, for a subset $\mathcal{G}$ of $\mathrm{Subgr}(G)$ we say $\mathcal{G}$ is **étale open** (**closed**, **compact**, **Hausdorff**, etc) if $\mathcal{G}$ is open (closed, compact, Hausdorff, etc) in the étale topology. Likewise, for a function $f$ from a topological space $X$ into $\mathrm{Subgr}(G)$ we say $f$ is **étale continuous** if $f$ is continuous when $\mathrm{Subgr}(G)$ is equipped with the étale topology.

Note: We use the adjective **compact** for a topological space $X$ in the sense of Hewitt-Ross [HRo]. Thus, every open covering of $X$ has a finite subcovering (but, in contrast to the terminology of Bourbaki, $X$ need not be Hausdorff).

*Remark 1.1: Categorical properties of the étale topology.*

(a) Subgroups: Let $H$ be a closed subgroup of $G$. Then a subgroup $H_0$ of $H$ is open in $H$ if and only if $H_0 = H \cap G_0$ with $G_0 \in \mathrm{Open}(G)$. Moreover, $\mathrm{Subgr}(H_0) = \mathrm{Subgr}(H) \cap \mathrm{Subgr}(G_0)$. Thus, the étale topology of $\mathrm{Subgr}(H)$ is the one induced from the étale topology of $\mathrm{Subgr}(G)$.

(b) Quotients: Let $N$ be a closed normal subgroup of $G$. Put $\bar{G} = G/N$ and let $\pi \colon G \to \bar{G}$ be the quotient map. Given $\bar{U} \in \mathrm{Open}(\bar{G})$, put $U = \pi^{-1}(\bar{U})$ and observe that $\pi^{-1}(\mathrm{Subgr}(\bar{U})) = \mathrm{Subgr}(U)$. It follows that the étale topology of $\mathrm{Subgr}(\bar{G})$ is the quotient topology of $\mathrm{Subgr}(G)$ via the quotient map $\pi \colon \mathrm{Subgr}(G) \to \mathrm{Subgr}(\bar{G})$. ∎

*Remark 1.2: Étale versus strict.* The strict topology of $\mathrm{Subgr}(G)$ is finer than the étale topology. Indeed, consider an open subgroup $U$ of $G$. Choose an open normal subgroup $N$ of $G$ in $U$. List the subgroups between $N$ and $U$ as $H_1, \ldots, H_n$. Then $\mathrm{Subgr}(U) = \bigcup_{i=1}^{n} \{A \in \mathrm{Subgr}(G) \mid AN = H_i\}$. Hence, $\mathrm{Subgr}(U)$ is strictly open (and closed).

Since $\mathrm{Subgr}(G)$ is strictly profinite, this gives the following chain of implications for a subset $\mathcal{G}$ of $\mathrm{Subgr}(G)$: $\mathcal{G}$ is étale closed $\Longrightarrow$ $\mathcal{G}$ is strictly closed $\Longleftrightarrow$ $\mathcal{G}$ is strictly compact $\Longrightarrow$ $\mathcal{G}$ is étale compact. ∎

The intersection of two étale open basic sets contains the trivial group. So, if $G \neq 1$, the étale topology of $\mathrm{Subgr}(G)$ is not Hausdorff. However, a subset $\mathcal{G}$ of $\mathrm{Subgr}(G)$ can be étale Hausdorff. Indeed, we will be looking for such $\mathcal{G}$ which are even étale profinite.

Denote the strict closure of a subset $\mathcal{G}$ of $\mathrm{Subgr}(G)$ (resp. a point $H \in \mathrm{Subgr}(G)$) by $\mathrm{StrictClosure}(\mathcal{G})$ (resp. $\mathrm{StrictClosure}(H)$).

Lemma 1.3: *Let $\mathcal{G}$ be a subset of $\mathrm{Subgr}(G)$.*

(a) *Let $H, H' \in \mathcal{G}$. Suppose $H \cap H'$ contains no $L$ which belongs to $\mathrm{StrictClosure}(\mathcal{G})$. Then $H$ and $H'$ can be separated by the étale topology of $\mathcal{G}$.*

(b) *Suppose $H \cap H'$ contains no $L \in \mathrm{StrictClosure}(\mathcal{G})$ for all distinct $H, H' \in \mathcal{G}$. Then $\mathcal{G}$ is étale Hausdorff.*

*Proof:* Statement (b) follows from (a). So, we prove (a). Assume $H$ and $H'$ cannot be separated by the étale topology of $\mathcal{G}$. Denote the set of all pairs $(U, U') \in \mathrm{Open}(G) \times$

Open($G$) with $H \leq U$ and $H' \leq U'$ by $\mathcal{U}$. Then, Subgr($U$) $\cap$ Subgr($U'$) $\cap \mathcal{G} \neq \emptyset$ for all $(U, U') \in \mathcal{U}$. Hence, Subgr($U$) $\cap$ Subgr($U'$) $\cap$ StrictClosure($\mathcal{G}$) $\neq \emptyset$ for all $(U, U') \in \mathcal{U}$. Each of the sets Subgr($U$) $\cap$ Subgr($U'$) $\cap$ StrictClosure($\mathcal{G}$) is strictly closed (Remark 1.2). The intersection of finitely many of them is a set of the same type. Hence, the intersection is nonempty. Since StrictClosure($\mathcal{G}$) is strictly compact, there is $L \in \bigcap_{(U,U') \in \mathcal{U}}$ Subgr($U$) $\cap$ Subgr($U'$) $\cap$ StrictClosure($\mathcal{G}$). It satisfies $L \leq H \cap H'$. This contradicts the assumption of the lemma. ∎

COROLLARY 1.4: *Let $\mathcal{G}$ be a subset of* Subgr($G$) *with* $1 \notin$ StrictClosure($\mathcal{G}$).

(a) *Let $H, H' \in \mathcal{G}$. Suppose $H \cap H' = 1$. Then $H$ and $H'$ can be separated by the étale topology of $\mathcal{G}$.*

(b) *Suppose $H \cap H' = 1$ for all distinct $H, H' \in \mathcal{G}$. Then $\mathcal{G}$ is étale Hausdorff .*

Here is a certain converse to Corollary 1.4:

LEMMA 1.5: *Let $G$ be a profinite group and $\mathcal{G}$ a subset of* Subgr($G$). *Suppose $\mathcal{G}$ is étale Hausdorff and contains at least two groups. Then $1 \notin$ StrictClosure($\mathcal{G}$).*

*Proof:* Let $H_1$ and $H_2$ be distinct groups in $\mathcal{G}$. Then there are disjoint étale open subsets $\mathcal{U}_1$ and $\mathcal{U}_2$ of $\mathcal{G}$ such that $H_i \in \mathcal{U}_i$, $i = 1, 2$. For each $i$ there is $U_i \in$ Open($G$) with $H_i \in \mathcal{G} \cap$ Subgr($U_i$) $\subseteq \mathcal{U}_i$. Let $U = U_1 \cap U_2$. Then $U \in$ Open($G$) and

$$\mathcal{G} \cap \text{Subgr}(U) \subseteq \mathcal{G} \cap \text{Subgr}(U_1) \cap \text{Subgr}(U_2) \subseteq \mathcal{U}_1 \cap \mathcal{U}_2 = \emptyset.$$

It follows, $1 \notin$ StrictClosure($\mathcal{G}$). ∎

## 2. Group Structures

The profinite group structures we introduce in this section replace the Artin-Schreier Structures of [HaJ1], the $\Gamma$-structures of [HaJ2], and the étale spaces of [Har]. The category of profinite group structures admits quotients (Example 2.6), fiber products (Construction 2.10), and certain inverse limits (Remark 2.9). These are the necessary tools to prove that solvability of finite embedding problems of a finite group structure **G** implies the solvability of arbitrary embedding problems for **G** (Proposition 4.2).

A **profinite group space** (also called a **profinite transformation group** [HaJ1, Section 1]) is a pair $(X, G)$ consisting of a profinite space $X$, a profinite group $G$, and a continuous action of $G$ on $X$ from the right (which we right exponentially). A **morphism** $\varphi \colon (X, G) \to (Y, H)$ of group spaces is a couple consisting of a continuous map $\varphi \colon X \to Y$ and a continuous group homomorphism $\varphi \colon G \to H$ satisfying $\varphi(x^g) = \varphi(x)^{\varphi(g)}$ for all $x \in X$ and $g \in G$. Composition of morphisms of profinite group spaces and the identity maps are morphisms of profinite group spaces satisfying the associativity law. Thus, the class of profinite groups spaces with their morphisms form a category.

Giving a profinite group space $(X, G)$ and an element $x \in X$, we let $S_x = \{g \in X \mid x^g = x\}$. It is a closed subgroup of $G$ called the **stabilizer** of $x$. If $\varphi \colon (X, G) \to (Y, H)$ is a morphism and $x \in X$, then $\varphi(S_x) \leq S_{\varphi(x)}$.

Every profinite group $G$ acts on $\mathrm{Subgr}(G)$ by conjugation. This action is both strictly continuous and étale continuous. Therefore, $(\mathrm{Subgr}(G), G)$ with $\mathrm{Subgr}(G)$ equipped either with the strict topology or the étale topology is a profinite group space. In Section 6 we will encounter our second basic example of profinite group spaces arising in the context of absolute Galois groups.

A **profinite group structure** is a triple $\mathbf{G} = (G, X, \delta)$ consisting of a profinite group space $(X, G)$ and an étale continuous map $\delta \colon X \to \mathrm{Subgr}(G)$. This object must satisfy the following conditions:

(1a) $G_{x^g} = G_x^g$ for all $x \in X$ and $g \in G$; thus $\delta$ is a morphism of group spaces.

(1b) $S_x \leq G_x$ for each $x \in X$.

Denote $\delta$ also by $\delta_{\mathbf{G}}$. The continuity of $\delta_{\mathbf{G}}$ means that $\{x \in X \mid G_x \leq U\}$ is an open subset of $X$ for each $U \in \text{Open}(G)$.

We write $\mathbf{G}$ also as $(G, X, G_x)_{x \in X}$ and refer to $\mathbf{G}$ as a **group structure** (ommiting "profinite").

A **morphism** of group structures

$$(2) \qquad \varphi \colon (G, X, G_x)_{x \in X} \to (H, Y, H_y)_{y \in Y}$$

is a couple consisting of a continuous group homomorphism $\varphi \colon G \to H$ and a continuous map $\varphi \colon X \to Y$ such that $\varphi(G_x) \leq H_{\varphi(x)}$ and $\varphi(x^g) = \varphi(x)^{\varphi(g)}$ for all $x \in X$ and $g \in G$. Thus, $\varphi(S_x) \leq S_{\varphi(x)}$ for each $x \in X$.

We call $\varphi$ an **epimorphism** if $\varphi(G) = H$, $\varphi(X) = Y$, and for each $y \in Y$ there is $x \in X$ with $\varphi(x) = y$ and $\varphi(G_x) = H_y$.

We call $\varphi$ a **cover** if $\varphi$ is an epimorphism with the following properties:

(3a) $\varphi$ maps each $G_x$ isomorphically onto $G_{\varphi(x)}$.

(3b) $\varphi(x) = \varphi(x')$ implies $x^k = x'$ for some $k \in \text{Ker}(\varphi)$.

If an epimorphism $\varphi \colon \mathbf{G} \to \mathbf{H}$ satisfies (3a) (but not necessarily (3b)), we say $\varphi$ is **rigid**. We call $\mathbf{G}$ **finite**, if both $G$ and $X$ are finite.

*Remark 2.1: Proper group structures.* Let $\mathbf{G} = (G, X, G_x)_{x \in X}$ be a group structure. Write $\mathcal{G} = \{G_x \mid x \in X\}$. We say $\mathbf{G}$ is **proper**, if $\delta_{\mathbf{G}} \colon X \to \mathcal{G}$ is an étale homeomorphism. Then $\mathcal{G}$ is étale profinite. Moreover, $S_x = G_x$ for each $x \in X$. Indeed, if $g \in G_x$, then $G_{x^g} = G_x^g = G_x$, hence $x^g = x$. Thus, $N_G(\Gamma) = \Gamma$ for each $\Gamma \in \mathcal{G}$. If $X = \{x\}$ consists of one element and $\sigma \in G$, then $x^\sigma = x$, hence $\sigma \in S_x = G_x$. Therefore, $G_x = G$. If $X$ contains at least two points, then $1 \notin \text{StrictClosure}(\mathcal{G})$ (Lemma 1.5).

Let $\mathbf{H} = (H, Y, H_y)_{y \in Y}$ be another proper group structure and $\varphi \colon G \to H$ a group homomorphism. Put $\mathcal{H} = \{H_y \mid y \in Y\}$. Suppose $\varphi(\mathcal{G}) \subseteq \mathcal{H}$. Then $\delta_{\mathbf{H}}^{-1} \circ \varphi \circ \delta_{\mathbf{G}}$ is a continuous map from $X$ into $Y$ which is compatible with the action of $G$ and $H$. This gives a unique extension of $\varphi \colon G \to H$ to a morphism $\varphi \colon \mathbf{G} \to \mathbf{H}$ satisfying $\varphi(G_x) = H_{\varphi(x)}$ for each $x \in X$.

Consider now a third proper group structure $\mathbf{A} = (A, I, A_i)_{i \in I}$. Let $\alpha \colon \mathbf{G} \to \mathbf{A}$ and $\beta \colon \mathbf{H} \to \mathbf{A}$ be morphisms. Suppose $\varphi(G_x) = H_{\varphi(x)}$, $\beta(H_y) = A_{\beta(y)}$, and $\alpha(G_x) = A_{\alpha(x)}$ for all $x \in X$ and $y \in Y$. Then $\alpha = \beta \circ \varphi$ as homomorphisms of groups implies $\alpha = \beta \circ \varphi$ as morphisms of group structures.

Finally suppose $\varphi\colon \mathbf{G} \to \mathbf{H}$ is a rigid epimorphism of group structures with $\mathbf{G}$ proper and $\mathrm{Ker}(\varphi) = 1$. Then $\varphi$ is an isomorphism and $\mathbf{H}$ is proper. Indeed, $\varphi\colon G \to H$ is an isomorphism. It remains to prove that $\varphi\colon X \to Y$ is an isomorphism. Since both $X$ and $Y$ are profinite spaces and $\varphi\colon X \to Y$ is continuous and surjective, it suffices to prove that $\varphi\colon X \to Y$ is injective. Consider $x, x' \in X$ with $\varphi(x) = \varphi(x')$. Then $\varphi(G_x) = \varphi(G_{x'})$. Hence, $G_x = G_{x'}$. Since $\mathbf{G}$ is proper, $x = x'$, as desired. ∎

Lemma 2.2: *Suppose (2) is a cover of groups structures. Then $\varphi(S_x) = S_{\varphi(x)}$ for each $x \in X$. In particular, if $H_y = S_y$ for each $y \in Y$, then $G_x = S_x$ for each $x \in X$.*

Proof: Let $x \in X$ and $y = \varphi(x)$. We have already mentioned that $\varphi(S_x) \leq S_y$. Also, $\varphi\colon G_x \to H_y$ is an isomorphism. Hence, in order to prove that $\varphi(S_x) = S_y$, it suffices to consider $g \in \varphi^{-1}(S_y) \cap G_x$ and to prove that $g \in S_x$.

Indeed, $\varphi(x^g) = y^{\varphi(g)} = y = \varphi(x)$. Hence, there is $k \in \mathrm{Ker}(\varphi)$ with $x^{gk} = x$. Thus, $gk \in S_x$. Therefore, $k \in \mathrm{Ker}(\varphi) \cap G_x = 1$. It follows that $x^g = x$, so $g \in S_x$.

Now assume $H_y = S_y$. Then, by the preceding paragraph, $\varphi(G_x) = H_y = \varphi(S_x)$. Since, $\varphi\colon G_x \to H_y$ is an isomorphism, $G_x = S_x$, as claimed. ∎

Lemma 2.3: *Let $(G, X, G_x)_{x \in X}$ be a group structure and $Y$ a closed subset of $X$. Then $\bigcup_{x \in Y} G_x$ is closed in $G$.*

Proof (After [Gil, Lemma 1.4]): Let $g \in G \smallsetminus \bigcup_{x \in Y} G_x$. For each $x \in Y$ there is an open normal subgroup $N_x$ of $G$ with $gN_x \cap G_x = \emptyset$. Thus, $g \notin G_x N_x$. As $G_x N_x \in \mathrm{Open}(G)$, continuity of $\delta_{\mathbf{G}}$ implies $V_x = \{y \in Y \mid G_y \leq G_x N_x\}$ is an open neighborhood of $x$ in $Y$. As $Y$ is compact, the covering $\{V_x \mid x \in Y\}$ has a finite subcovering $\{V_{x_1}, \ldots, V_{x_n}\}$. Then $N = \bigcap_{i=1}^{n} N_{x_i}$ is an open normal subgroup of $G$ and $g \notin G_y N$ for each $y \in Y$. Therefore, $gN \subseteq G \smallsetminus \bigcup_{y \in Y} G_y$. ∎

Proper group structures are our main subject of research. We have introduced the more general concept of group structures in order to be able to extend the basic operations of the category of profinite groups to the category of group structures. This is not always possible in the category of proper group structures. For example, a quotient of a proper group structure need not be proper (Example 2.5).

*Example 2.4: Absolute Galois group structures.* An **absolute Galois group structure** is a group structure $\mathbf{G} = (\mathrm{Gal}(K), X, \mathrm{Gal}(K_x))_{x \in X}$ where each $K_x$ is a separable algebraic extension of $K$. Let $\mathbf{H} = (\mathrm{Gal}(L), Y, L_y)_{y \in Y}$ be another proper absolute Galois group structure. Suppose $K \subseteq L$ and for each $y \in Y$ there is $x \in X$ with $L_y \cap K_s = K_x$. By Remark 2.1, $\mathrm{res}_{L_s/K_s} \colon \mathrm{Gal}(L) \to \mathrm{Gal}(K)$ extends to a unique morphism $\rho \colon \mathbf{H} \to \mathbf{G}$ of group structures satisfying $\mathrm{res}_{L_s/K_s}(\mathrm{Gal}(L_y)) = \mathrm{Gal}(K_x)$ for all $y \in Y$ and $x = \rho(y)$. We denote this morphism by res if the reference to $K$ and $L$ is clear from the context. ∎

*Example 2.5: Quotient maps.* Let $\mathbf{G} = (G, X, G_x)_{x \in X}$ be a group structure and $N$ a closed normal subgroup of $G$. Put $\bar{G} = G/N$ and $\bar{X} = X/N$. Let $\pi \colon G \to \bar{G}$ and $\pi \colon X \to \bar{X}$ the quotient maps: $\pi(g) = \bar{g} = gN$ and $\pi(x) = \bar{x} = \{x^\nu \mid \nu \in N\}$. Then $\bar{X}$ is a profinite space [HaJ1, Claim 1.6]. For each $x \in X$ let $\bar{G}_{\bar{x}} = \pi(G_x) = G_x N/N$.

Consider $\bar{U} \in \mathrm{Open}(\bar{G})$. Put $U = \pi^{-1}(\bar{U})$. Then $\pi^{-1}(\{\bar{x} \in \bar{X} \mid \bar{G}_{\bar{x}} \leq \bar{U}\}) = \{x \in X \mid G_x \leq U\}$. Hence, the map $\delta_{\bar{\mathbf{G}}} \colon \bar{X} \to \mathrm{Subgr}(\bar{G})$ given by $\delta_{\bar{\mathbf{G}}}(\bar{x}) = \bar{G}_{\bar{x}}$ is étale continuous. Also, $\bar{x}^{\bar{\sigma}} = \bar{x}$ implies $\bar{\sigma} \in \bar{G}_{\bar{x}}$. Thus, $\bar{\mathbf{G}} = (\bar{G}, \bar{X}, \bar{G}_{\bar{x}})_{\bar{x} \in \bar{X}}$ is a group structure which we denote by $\mathbf{G}/N$ and $\pi \colon \mathbf{G} \to \bar{\mathbf{G}}$ is an epimorphism. Moreover, $\pi(G_x) = G_{\pi(x)}$ for every $x \in X$. We call $\pi$ the **quotient map**. If $G_x \cap N = 1$ for each $x \in X$, then $\pi$ is a cover.

Let $\mathcal{G} = \{G_x \mid x \in X\}$ and $\bar{\mathcal{G}} = \{\bar{G}_{\bar{x}} \mid \bar{x} \in \bar{X}\}$. Then $\pi$ induces a strictly continuous map of $\mathrm{Subgr}(G)$ onto $\mathrm{Subgr}(\bar{G})$ and $\pi(\mathcal{G}) = \bar{\mathcal{G}}$. Thus, if $1 \notin \mathrm{StrictClosure}(\bar{\mathcal{G}})$, then $1 \notin \mathrm{StrictClosure}(\mathcal{G})$.

Conversely, every cover $\varphi \colon \mathbf{G} \to \mathbf{H}$ of group structures is isomorphic to the quotient map $\mathbf{G} \to \mathbf{G}/\mathrm{Ker}(\varphi)$. Indeed, let $\mathbf{H} = (H, Y, H_y)_{y \in Y}$. Then $\varphi$ induces a bijective continuous map $\bar{\varphi} \colon \bar{X} \to Y$. As both $\bar{X}$ and $Y$ are profinite, $\bar{\varphi}$ is a homeomorphism.

Consider now the case where $N = G$. Suppose $|\bar{X}| > 1$. Then $\bar{G} = 1$, so the forgetful map $\delta_{\bar{\mathbf{G}}}$ is not injective. Thus, $\bar{\mathbf{G}}$ need not be proper even if $G$ is proper. This is one of the reasons why we work in the category of group structures and not in the category of proper group structures, which may look at first glance more attractive. Another reason is the need to use morphisms called "Galois approximations" (Section 14). The target objects of Galois approximations are group structures which need not be proper.

Quotient maps of group structures has the universal property of quotient maps of groups. Thus, if $\pi \colon \mathbf{G} \to \bar{\mathbf{G}}$ and $\varphi \colon \mathbf{G} \to \mathbf{H}$ are quotient maps satisfying $\mathrm{Ker}(\varphi) \leq \mathrm{Ker}(\pi)$, then there is a unique quotient map $\psi \colon \mathbf{H} \to \mathbf{G}/N$ satisfying $\psi \circ \varphi = \pi$. Moreover, $\pi$ is a cover if and only if $\varphi$ and $\psi$ are covers. Also, if $N'$ is a closed normal subgroup of $G$ which is contained in $N$, then there is a natural isomorphism $\mathbf{G}/N \cong (\mathbf{G}/N')/(N/N')$. ∎

*Remark 2.6: Sub-group-structures.* Let $\mathbf{G} = (G, X, G_x)_{x \in X}$ and $\mathbf{H} = (H, Y, H_y)_{y \in Y}$ be group structures. We say $\mathbf{H}$ is a **sub-group-structure** of $\mathbf{G}$ if $H \leq G$, $Y$ is a subspace of $X$, and $H_y = G_y$ for each $y \in Y$. If $\mathbf{G}$ proper, then so is $\mathbf{H}$.

Suppose we start with a group $G$, a profinite space $X$, and for each $x \in X$ a closed subgroup $G_x$ of $G$. Consider a closed subgroup $H$ of $G$ which contains all $G_x$. If $U$ is an open subgroup of $G$, then $V = U \cap H$ is an open subgroup of $H$. Conversely, for each open subgroup $V$ of $H$ there is an open subgroup $U$ of $G$ with $V = U \cap H$. In each case $\{x \in X \mid G_x \leq U\} = \{x \in X \mid G_x \leq V\}$. So, if one of the sets is open, so is the other. Thus, $(G, X, G_x)_{x \in X}$ is a group structure if and only if $(H, X, G_x)_{x \in X}$ is a group structure. ∎

*Remark 2.7: Inverse limit of group structures.* Let $\mathbf{G}_i = (G_i, X_i, G_{i,x})_{x \in X_i}$, $i \in I$, be an inverse system of group structures with connecting homomorphisms $\pi_{ji} \colon \mathbf{G}_j \to \mathbf{G}_i$. Suppose all $\pi_{ji}$ are surjective. Put $G = \varprojlim G_i$, $X = \varprojlim X_i$, and let $\pi_i$ be the projections on the $i$th coordinate of $G$ and $X$. Then $\pi_i$ is surjective. Since the $\pi_{ji}$'s commute with the action of $G_i$ on $X_i$, they define a continous action of $G$ on $X$. Since the $\pi_{ji}$'s commute with the maps $\delta_{\mathbf{G}_i} \colon X_i \to \mathrm{Subgr}(G_i)$, they define an map $\delta \colon X \to \mathrm{Subgr}(G)$ which is étale continuous. Specifically, for each $x = (x_i)_{i \in I}$ in $X$ we have $\delta(x) = G_x = \varprojlim G_{i,x_i}$.

Indeed, let $U$ be an open subgroup of $G$. Then there is $i \in I$ with $\mathrm{Ker}(\pi_i) \leq U$.

Hence,

(4) $$\{x \in X \mid G_x \leq U\} = \pi_i^{-1}\big(\{x_i \in X_i \mid G_{x_i} \leq \pi_i(U)\}\big).$$

Since $\pi_i$ is surjective and $\delta_{\mathbf{G}_i}$ is étale continuous, $\pi_i(U)$ is open in $G_{x_i}$ and the argument of $\pi^{-1}$ in (4) is open. Thus, the left hand side of (4) is open.

Since each $\delta_{\mathbf{G}_i}$ commutes with the action of $\mathbf{G}_i$, the map $\delta$ commutes with the action of $G$. Finally, with $x = (x_i)_{i \in I}$, it follows from $S_{x_i} \leq G_{x_i}$ that $S_x \leq G_x$. Therefore, $\mathbf{G} = (G, X, G_x)_{x \in X}$ is a group structure.

If each $\pi_{ji}$ is rigid, then each $\pi_i$ is rigid. If each $\pi_{ji}$ is a cover, then so is each $\pi_i$. Indeed, Let $x = (x_k)_{k \in I}$ and $y = (y_k)_{k \in I}$ be elements of $X$ satisfying $x_i = y_i$. Then, for each $j \geq i$ the closed subset $K_j = \{\kappa \in \mathrm{Ker}(\pi_{ji}) \mid x_j^\kappa = y_j\}$ of $G_j$ is not empty. If $k \geq j$, then $\pi_{kj}(K_k) \subseteq K_j$. Therefore, there is $\kappa \in G$ with $\pi_j(\kappa) \in K_j$ for all $j \geq i$. This $\kappa$ belongs to $\mathrm{Ker}(\pi_i)$ and $x^\kappa = y$, as claimed. ∎

LEMMA 2.8: *Let* $\mathbf{G} = (G, X, G_x)_{x \in X}$ *be a group structure and* $\mathcal{N}$ *an inductive collection of closed normal subgroups of* $G$. *Then* $\mathbf{G} = \varprojlim \mathbf{G}/N$ *where* $N$ *ranges over* $\mathcal{N}$.

*Proof:* The only point which is perhaps not clear is $X = \varprojlim X/N$. To prove this equality define a map $f \colon X \to \varprojlim X/N$ by $f(x) = (x^N)_{N \in \mathcal{N}}$, where $x^N = \{x^\nu \mid \nu \in N\}$. Then $f$ is continuous. Compactness of $X$ implies $f$ is surjective. Since both $X$ and $\varprojlim X/N$ are profinite spaces, it suffices now to prove that $f$ is injective.

Consider distinct elements $x, y \in X$. Choose disjoint open subsets $U$ and $V$ of $X$ with $x \in U$ and $y \in V$. Since the action of $G$ on $X$ is continuous, $x$ has an open neigborhood $U_0$ and there is $N \in \mathcal{N}$ with $U_0^N \subseteq U$. Then $x^\nu \notin V$, so $x^\nu \neq y$ for all $\nu \in N$. Therefore, $f(x) \neq f(y)$. ∎

Construction 2.9: *Fiber products.* Let $\mathbf{A} = (A, I, A_i)_{i \in I}$, $\mathbf{B} = (B, J, B_j)_{j \in J}$, and $\mathbf{G} = (G, X, G_x)_{x \in X}$ be group structures. Let $\alpha \colon \mathbf{B} \to \mathbf{A}$ and $\varphi \colon \mathbf{G} \to \mathbf{A}$ be morphisms of group structures. Put

(5a) $H = B \times_A G = \{(b, g) \in B \times G \mid \alpha(b) = \varphi(g)\}$,

(5b) $Y = J \times_I X = \{(j, x) \in J \times X \mid \alpha(j) = \varphi(x)\}$, and

(5c) $H_y = B_j \times_A G_x = \{(b, g) \in B_j \times G_x \mid \alpha(b) = \varphi(g)\}$ for $y = (j, x) \in Y$.

Define a continuous action of $H$ on $Y$ by $(j,x)^{(b,g)} = (j^b, x^g)$. We claim: $\mathbf{H} = (H, Y, H_y)_{y \in Y}$ is a group structure.

To verify the claim it suffices to prove that the map $y \mapsto H_y$ is étale continuous. Indeed, let $y = (j, x)$. Consider an open subgroup $W$ of $H$ which contains $H_y = B_j \times_A G_x$. Let $\mathcal{U}$ be the set of all open subgroups of $B$ which contain $B_j$. Let $\mathcal{V}$ be the set of all open subgroups of $G$ which contain $G_x$. The intersection of all $U \times_A V$ with $U \in \mathcal{U}$ and $V \in \mathcal{V}$ is $B_j \times_A G_x$. Since $H \smallsetminus W$ is closed, there are an open subgroup $U$ of $B$ and an open subgroup $V$ of $G$ with $H_y \leq U \times_A V \leq W$. The set $Y_0 = \{(j', x') \in Y \mid B_{j'} \leq U, \ \ G_{x'} \leq G\}$ is an open neighborhood of $y$ in $Y$ and $H_{(j', x')} \leq W$ for each $(j', x') \in Y_0$. So, the above map is continuous.

Finally let $\beta \colon H \to B$, $\beta \colon Y \to J$, $\psi \colon H \to G$, and $\psi \colon Y \to X$ be the projections on the coordinates. Then the following diagram of group structures is commutative:

(6)
$$
\begin{array}{ccc}
\mathbf{H} & \xrightarrow{\ \psi\ } & \mathbf{G} \\
\Big\downarrow{\scriptstyle \beta} & & \Big\downarrow{\scriptstyle \varphi} \\
\mathbf{B} & \xrightarrow{\ \alpha\ } & \mathbf{A}
\end{array}
$$

If both $\mathbf{B}$ and $\mathbf{G}$ are finite, then so is $\mathbf{H}$. ∎

*Definition 2.10: Cartesian squares.* Let (6) be a commutative diagram of group structures. Call (6) a **cartesian square** if this holds: For all group structures $\mathbf{F}$ and morphisms $\beta' \colon \mathbf{F} \to \mathbf{B}$ and $\psi' \colon \mathbf{F} \to \mathbf{G}$ with $\alpha \circ \beta' = \varphi \circ \psi'$ there is a unique morphism $\varepsilon \colon \mathbf{F} \to \mathbf{H}$ satisfying $\beta \circ \varepsilon = \beta'$ and $\psi \circ \varepsilon = \psi'$. ∎

LEMMA 2.11: *Let (6) be a commutative diagram of group structures.*

(a) *Suppose $\mathbf{H} = \mathbf{B} \times_{\mathbf{A}} \mathbf{G}$ and $\beta$, $\psi$ are the coordinate projections. Then (6) is a cartesian square.*

(b) *Suppose (6) is a cartesian square. Put $\mathbf{H}' = \mathbf{B} \times_{\mathbf{A}} \mathbf{G}$. Let $\psi' \colon \mathbf{H}' \to \mathbf{B}$ and $\beta' \colon \mathbf{H}' \to \mathbf{G}$ be the projection maps. Then there is a unique isomorphism $\gamma \colon \mathbf{H}' \to \mathbf{H}$ with $\psi \circ \gamma = \psi'$ and $\beta \circ \gamma = \beta'$.*

*Proof:* Statement (a) follows from the definition of $\mathbf{B} \times_{\mathbf{A}} \mathbf{G}$. Statement (b) follows from (a) and from the uniqueness of $\varepsilon$ in Definition 2.10. ∎

LEMMA 2.12: *Suppose (6) is a cartesian square of group structures. Then:*

(a) $\beta \colon \mathrm{Ker}(\psi) \to \mathrm{Ker}(\alpha)$ *is an isomorphism.*

(b) *For each $y \in Y$, $\psi \colon H_y \to G_{\psi(y)}$ is injective if and only if $\alpha \colon B_{\beta(y)} \to A_{\alpha(\beta(y))}$ is injective.*

(c) *If $\alpha$ is a cover, then $\psi$ is a cover.*

(d) *If $\psi$ is a cover and $\varphi$ is an epimorphism, then $\alpha$ is a cover.*

*Proof:* By Lemma 2.11(b) we may assume that $\mathbf{H}$ is $\mathbf{B} \times_{\mathbf{A}} \mathbf{G}$ and $\beta, \psi$ are the projections.

*Proof of (a) and (b):* By assumption, $\mathrm{Ker}(\psi) = \mathrm{Ker}(\alpha) \times \{1\}$, which gives (a). Similarly, for $y = (j, x) \in Y$, (5c) implies that $\beta$ maps $\mathrm{Ker}(\psi) \cap H_y$ isomorphically onto $\mathrm{Ker}(\alpha) \cap B_j$. This gives (b).

*Proof of (c):* Suppose $\alpha$ is a cover. Then $\alpha(B) = A$. Hence, for each $g \in G$ there is $b \in B$ with $\alpha(b) = \varphi(g)$. Therefore, $(b, g) \in H$ and $\psi(b, g) = g$. Thus, $\psi(H) = G$. Similarly, $\psi(Y) = X$ and $\psi(H_y) = G_{\psi(y)}$ for each $y \in Y$. Since $\alpha \colon B_{\beta(y)} \to A_{\alpha(\beta(y))}$ is an isomorphism, (b) implies $\psi \colon H_y \to G_{\psi(y)}$ is an isomorphism.

Finally, suppose $\psi(j, x) = \psi(j', x')$. Then $x = x'$ and $\alpha(j) = \alpha(j')$. The rigidity of $\alpha$ gives $b \in \mathrm{Ker}(\alpha)$ with $j' = j^b$. Then $(b, 1) \in \mathrm{Ker}(\psi)$ and $(j', x') = (j, x)^{(b,1)}$. This proves $\psi$ is a cover.

*Proof of (d):* By assumption, $\alpha(\beta(H)) = \varphi(\psi(H)) = A$ and $\alpha(\beta(Y)) = \varphi(\psi(Y)) = I$. Hence, $\alpha(B) = A$ and $\alpha(J) = I$.

Now let $j \in J$ and $i = \alpha(j)$. As $\varphi$ is an epimorphism, there is $x \in X$ with $\varphi(x) = i$ and $\varphi(G_x) = A_i$. Put $y = (j, x)$. As $\psi$ is a cover, $\psi \colon H_y \to G_x$ is an isomorphism. Hence, $A_i \geq \alpha(B_j) \geq \alpha(\beta(H_y)) = \varphi(\psi(H_y)) = A_i$, so $\alpha(B_j) = A_i$. We conclude from (b) that $\alpha \colon B_j \to A_i$ is an isomorphism.

Finally, consider $j, j' \in J$ with $\alpha(j) = \alpha(j')$. Choose $x \in X$ with $\alpha(j) = \alpha(j') = \varphi(x)$. Then $\psi(j, x) = \psi(j', x)$. So, there is $(b, 1) \in \mathrm{Ker}(\psi)$ with $(j', x) = (j, x)^{(b,1)}$. Hence, $b \in \mathrm{Ker}(\alpha)$ and $j' = j^b$. This proves $\alpha$ is a cover. ∎

## 3. Completion of a Cover to a Cartesian Square

There are many occasions in this work where a group structure $\mathbf{G} = (G, X, G_x)_{x \in X}$ is given and we need to define a morphism $f \colon (X, G) \to (Y, H)$ of group spaces, where $f \colon G \to H$ is a given homomorphism. If the set of $G$-orbits of $X$ has a closed system of representatives $X'$ (also called a **fundamental domain**), we may first define $f$ on $X'$ and then extend it to $X$ by the rule

$$(0) \qquad\qquad f(x^\sigma) = f(x)^{f(\sigma)}, \qquad \sigma \in G.$$

This could considerably simplify the proof of the Main Theorem. Unfortunately, fundamental domains do not always exist. One may find a counter example of J. L. Kelly on page 473 of [ArK]. Instead we produce a "special partition" of $\mathbf{G}$ giving rise to a subset $X'$ of $X$ that "approximates" a fundamental domain in a way that allows the definition of the desired function $f$.

Lemmas 3.1, 3.2, 3.3, and 3.4 below prepare ingredients of the construction of special partitions in Lemma 3.6. The definition of "special partition" appears in Lemma 3.5. It follows by a specification of the above mentioned set $X'$.

LEMMA 3.1: *Let $G$ be a profinite group acting continuously on a compact Hausdorff space $X$. Then:*

(a) *$S_x$ is a closed subgroup of $G$.*

(b) *The map $x \mapsto S_x$ from $X$ to $\mathrm{Subgr}(G)$ is étale continuous.*

*Proof of (a):* The action $a \colon X \times G \to X$ and the projection $p \colon X \times G \to X$ are continuous, so $\{x\} \times S_x = p^{-1}(x) \cap a^{-1}(x)$ is closed in $X \times G$. Therefore, $S_x$ is closed in $G$.

*Proof of (b):* Let $N$ be an open normal subgroup of $G$ and let $x \in X$. We have to find an open neighborhood $V$ of $x$ with $S_y \leq S_x N$ for all $y \in V$.

CASE A: *$G$ is finite and $N = 1$.* Consider $\sigma \in G \smallsetminus S_x$. Then, $x^\sigma \neq x$. Since $X$ is Hausdorff, it has disjoint open subsets $U_1, U_2$ with $x \in U_1$ and $x^\sigma \in U_2$. Then $V_\sigma = U_1 \cap U_2^{\sigma^{-1}}$ is an open neighborhood of $x$. If $y \in V_\sigma$, then $y \in U_1$ and $y^\sigma \in U_2$, so $y^\sigma \neq y$. Since $G$ is finite, $V = \bigcap_{\sigma \in G \smallsetminus S_x} V_\sigma$ is open. Each $y \in V$ satisfies $S_y \leq S_x$.

CASE B: *The general case.* The quotient space $\bar{X} = X/N$ is Hausdorff [Bre, Thm. 3.1(1)] and $\bar{G} = G/N$ acts continuously on $\bar{X}$. Use a bar for reduction modulo $N$. Case A gives an open neighborhood $\bar{V}$ of $\bar{x}$ in $\bar{X}$ with $S_{\bar{y}} \leq S_{\bar{x}}$ for each $\bar{y} \in \bar{V}$. Then the preimage $V$ of $\bar{V}$ in $X$ is an open neighborhood of $x$ in $X$. For each $y \in V$ we have $S_{\bar{y}} = S_y N/N$. Hence, $S_y \leq S_x N$. ∎

LEMMA 3.2: *Let $Y$ be a profinite space and $A, B$ disjoint closed subsets. Then there are disjoint open-closed subsets $U, V$ with $A \subseteq U$ and $B \subseteq V$.*

*Proof:* As a profinite space, $Y$ is compact and Hausdorff. Hence, it has disjoint open subsets $U', V'$ with $A \subseteq U'$ and $B \subseteq V'$. The set $U'$ is a union of open-closed subsets. Since $A$ is compact, finitely many of them cover $A$. Their union $U$ is an open-closed subset satisfying $A \subseteq U \subseteq U'$. Similarly, $Y$ has an open-closed subset $V$ with $B \subseteq V \subseteq V'$. It satisfies $U \cap V = \emptyset$. ∎

LEMMA 3.3: *Let $(X, G)$ be a profinite group space, $X$ a profinite space, $x \in X$, and $V$ an open neighborhood of $x$. Suppose $x^G \subseteq V$. Then $x$ has an open-closed $G$-invariant neighborhood $W$ with $W \subseteq V$.*

*Proof:* Denote the images of points and subsets of $X$ under the quotient map $\pi\colon X \to X/G$ by a bar. Since $F = X \smallsetminus V$ is closed in $X$, $\bar{F}$ is closed in $X/G$. Moreover, $\bar{x} \notin \bar{F}$. Lemma 3.2 gives an open-closed subset $\bar{W}$ with $\bar{x} \in \bar{W}$ and $\bar{W} \cap \bar{F} = \emptyset$. Put $W = \pi^{-1}(\bar{W})$. Then $W$ is open-closed in $X$, invariant under $G$, and $x \in W \subseteq V$. ∎

LEMMA 3.4: *Let $(X, G)$ be a profinite group space, $X$ a profinite space, $x \in X$, and $H$ an open subgroup of $G$. Suppose $S_x \leq H$. Write $G = \bigcup_{\rho \in R} H\rho$. Then $x^H$ has an $H$-invariant open-closed neighborhood $U$ satisfying $U^G = \bigcup_{\rho \in R} U^\rho$.*

*Proof:* The closed sets $x^{H\rho}$, $\rho \in R$, are disjoint, because $S_x \leq H$. Hence, $X$ has open disjoint sets $V_\rho$ satisfying $x^{H\rho} \subseteq V_\rho$, $\rho \in R$. For each $\rho \in R$ we have $x^H \subseteq V_\rho^{\rho^{-1}}$. By Lemma 3.3, with $H$ replacing $G$, there is an $H$ invariant open-closed set $U_\rho$ with $x^H \subseteq U_\rho \subseteq V_\rho^{\rho^{-1}}$.

Now consider the $H$-invariant open-closed set $U = \bigcap_{\rho \in R} U_\rho$. It satisfies $U^\rho \subseteq V_\rho$ for each $\rho \in R$, so the $U^\rho$ are disjoint. Therefore, $U^G = \bigcup_{\rho \in R} U^\rho$. ∎

*Definition 3.5: Special partition.* Let $\mathbf{G} = (G, X, G_x)_{x \in X}$ be a group structure. A **special partition** of $\mathbf{G}$ is a data $(G_i, X_i, R_i)_{i \in I_0}$ satisfying these conditions:

(2a) $I_0$ is a finite set which is disjoint from $X$.

(2b) $X_i$ is an open-closed subset of $X$, $i \in I_0$.

(2c) $G_i$ is an open subgroup of $G$ containing $G_x$ for all $x \in X_i$ and $i \in I_0$.

(2d) $G_i = \{\sigma \in G \mid X_i^{\sigma} = X_i\}$, $i \in I_0$.

(2e) $R_i$ is a finite subset of $G$ and $G = \bigcup_{\rho \in R_i} G_i \rho$, $i \in I_0$.

(2f) $X = \bigcup_{i \in I_0} \bigcup_{\rho \in R_i} X_i^{\rho}$.

Here is a consequence of (2a)-(2f):

(2g) Suppose $i, j \in I_0$ and $X_i^{\sigma} \cap X_j \neq \emptyset$. Then, $i = j$ and $\sigma \in G_i$.

To prove (2g) write $\sigma = \zeta \rho$ with $\zeta \in G_i$ and $\rho \in R_i$. By (2d), $X_i^{\rho} \cap X_j \neq \emptyset$. Hence, by (2f), $i = j$ and $X_i^{\rho} = X_i$. By (2d), $\rho \in G_i$. Therefore, $\sigma \in G_i$.

By (2d), each $R_i$ can be replaced by every set $R_i'$ satisfying $G = \bigcup_{\rho \in R_i'} G_i \rho$. Thus, we also call $(G_i, X_i)_{i \in I_0}$ a **special partition** of $\mathbf{G}$ is there exist $R_i$, $i \in I_0$ satisfying (2a)-(2f). ∎

Suppose now $\mathbf{G} = (G, X, G_x)_{x \in X}$ is a group structure, $(Y, H)$ is a profinite group space, and $\varphi \colon G \to H$ is a homomorphism which we wish to extend to a morphism $\varphi \colon (X, G) \to (Y, H)$ of profinite group spaces. We construct a special partition $(X_i, G_i, R_i)_{i \in I_0}$ of $\mathbf{G}$ such that $\varphi$ has a natural definition on $X' = \bigcup_{i \in I_0} X_i$ satisfying $\varphi(x) = \varphi(x)^{\varphi(\sigma)}$ for all $x \in X'$ and $\sigma \in G$ with $x^{\sigma} \in X'$. Then $\varphi(x^{\tau}) = \varphi(x)^{\varphi(\tau)}$ for arbitrary $\tau \in G$ will define the desired extension $\varphi$.

It allows extending each epimorphism $\varphi$ of $G$ onto a finite group $A$ to an epimorphisms of $\mathbf{G}$ onto a finite group structure $\mathbf{A} = (A, I, A_i)_{i \in I}$ (Lemma 3.7). Consequently, each cover $\psi \colon \mathbf{H} \to \mathbf{G}$ of group structures with a finite kernel can be completed to a cartesian square as in (6) of Section 2 such that $\alpha \colon \mathbf{B} \to \mathbf{A}$ is a cover of finite group structures (Lemma 3.9). The latter result is a main ingredient in the transition from solving finite embedding problems to solving arbitrary embedding problems of projective group structures (Proposition 4.2).

LEMMA 3.6: *Let* $\mathbf{G} = (G, X, G_x)_{x \in X}$ *be a group structure, $Y$ be a subset of $X$, and $Y_0$ a finite subset of $Y$. Suppose $X = Y^G$ and the elements of $Y_0$ belong to distinct $G$-orbits. For each $y \in Y$ let $G'_y$ be an open subgroup of $G$ containing $G_y$ and $V_y$ an open neighborhood of $y^{G'_y}$ in $X$. Then there exists a finite subset $\{y_i \mid i \in I_0\}$ of $Y$ containing $Y_0$ and a special partition $(G'_{y_i}, X_i)_{i \in I_0}$ of $\mathbf{G}$ such that $y_i \in X_i \subseteq V_{y_i}$ for all $i \in I_0$.*

*Proof:* We may assume $Y$ is a (not necessarily closed) system of representatives of the $G$-orbits of $X$. For each $y \in Y$ use Lemma 3.4 to replace $V_y$ by another set, if necessary, to assume:

(3a) $V_y$ is open-closed, $G'_y$-invariant and $y^{G'_y} \subseteq V_y$.

(3b) Writing $G = \bigcup_{\rho \in R_y} G'_y \rho$, we have $V_y^G = \bigcup_{\rho \in R_y} V_y^\rho$.

The rest of the proof has three parts.

PART A: *Finite covering of $X$.* By assumption, $X = \bigcup_{y \in Y} y^G \subseteq \bigcup_{y \in Y} V_y^G$. Hence, by compactness, there is a finite subset $\{y_i \mid i \in I_0\}$ of $Y$ with $X = \bigcup_{i \in I_0} V_{y_i}^G$. Add the elements of $Y_0$ to $\{y_i \mid i \in I_0\}$, if necessary, to assume that $Y_0 \subseteq \{y_i \mid i \in I_0\}$. By our choice of $Y$, the sets $y_i^G$, $i \in I_0$, are closed and disjoint. Hence, there are disjoint open subsets $W'_i$ with $y_i^G \subseteq W'_i$, $i \in I_0$. For each $i \in I_0$ Lemma 3.3 gives a $G$-invariant open-closed set $W_i$ with $y_i^G \subseteq W_i \subseteq V_{y_i}^G \cap W'_i$.

PART B: *Making $V_{y_i}$ smaller.* By Part A, $y_i \in W_i \smallsetminus \bigcup_{j \neq i} W_j \subseteq V_{y_i}^G \smallsetminus \bigcup_{j \neq i} W_j$ and $\bigcup_{j \neq i} W_j$ is $G$-invariant. Let $V_i = V_{y_i} \smallsetminus \bigcup_{j \neq i} W_j$. Then $V_i$ is a $G'_{y_i}$-invariant open-closed set which, by (3), satisfies

(4) $V_i^G = \bigcup_{\rho \in R_i} V_i^\rho$

where $R_i = R_{y_i}$. Moreover, $y_i \in V_i \smallsetminus \bigcup_{j \neq i} V_j^G$. Indeed $y_i \in W_i$. If $y_i \in V_j^G$ for $j \neq i$, then there is $\sigma \in G$ with $y_i^\sigma \in V_j$, so $y_i^\sigma \notin W_i$. But $W_i$ is $G$-invariant. Hence, $y_i \in W_i$, a contradiction.

We claim that $X = \bigcup_{i \in I_0} V_i^G$. Indeed, let $x \in X$. If there is $i$ with $x \in W_i$, then $x \notin \bigcup_{j \neq i} W_j$. So, $x \in V_i^G$. Else, $x \notin \bigcup_{j \in I_0} W_j$ but there is $i$ with $x \in V_{y_i}^G$ (Part A). So, $x \in V_{y_i}^G \smallsetminus \bigcup_{j \neq i} W_j = V_i^G$.

PART C: *Separating $V_i$.* Let $X_i = V_i \smallsetminus \bigcup_{j=1}^{i-1} V_j^G$, $i \in I_0$. Then $X = \bigcup_{i \in I_0} X_i^G$. Also, $X_i$ is a $G'_{y_i}$-invariant open-closed neighborhood of $y_i$ and $X_i \subseteq V_{y_i}$. By (4), $X_i^G = \bigcup_{\rho \in R_i} X_i^\rho$.

Finally, consider $\sigma \in G$ with $X_i^\sigma = X_i$. Write $\sigma = \zeta \rho$ with $\zeta \in G'_{y_i}$ and $\rho \in R_i$. Then $X_i = X_i^\sigma = X_i^\rho$. By the preceding paragraph, $\rho \in G'_{y_i}$. Therefore, $\sigma \in G'_{y_i}$. ∎

LEMMA 3.7: *Let* $\mathbf{G} = (G, X, G_x)_{x \in X}$ *be a group structure, $A$ a finite group, and $\varphi \colon G \to A$ an epimorphism. Then:*

(a) $\varphi$ *extends to an epimorphism $\varphi$ of $\mathbf{G}$ onto a finite group structure $\mathbf{A} = (A, I, A_i)_{i \in I}$.*

(b) *Let $X_0$ be a finite subset of $X$. Then $\varphi$ may be constructed in (a) with $\varphi(G_x) = A_{\varphi(x)}$ for each $x \in X_0$.*

(c) *Suppose $X = \bigcup_{j \in J} Y_j$ with $J$ finite, each $Y_j$ is open-closed, $G$ permutes the $Y_j$'s, and $Y_j^\nu = Y_j$ for all $j \in J$ and $\nu \in \mathrm{Ker}(\varphi)$. Then $\varphi$ may be constructed in (a) such that $\varphi(Y_j)$, $j \in J$, are disjoint.*

(d) *Let $x_1, \dots, x_m$ be elements of $X$ lying in distinct $G$-orbits. Then $\varphi$ may be constructed in (a) such that $\varphi(x_1), \dots, \varphi(x_n)$ lie in distinct $A$-orbits.*

*Proof of (a):* We may assume $A = G/N$ with $N = \mathrm{Ker}(\varphi)$. Both maps $x \mapsto G_x$ and $x \mapsto S_x$ of $X$ into $\mathrm{Subgr}(G)$ are étale continuous (by definition and by Lemma 3.1). Hence, for each $y \in X$ the set $V_y = \{x \in X \mid S_x \leq S_y N, \ G_x \leq G_y N\}$ is open and contains $y^N = y^{S_y N}$. Lemma 3.6, with $S_y N$ replacing $G'_y$, gives a finite subset $\{y_i \mid i \in I_0\}$ of $X$ and a special partition $(S_{y_i} N, X_i)_{i \in I_0}$ of $\mathbf{G}$ such that

(5) $y_i \in X_i \subseteq V_{y_i}$ for all $i \in I_0$.

Thus, the following holds:

(6a) $X_i$ is open closed, $i \in I_0$.

(6b) $S_{y_i} N = \{\sigma \in G \mid X_i^\sigma = X_i\}$.

(6c) $X = \bigcup_{i \in I_0}^n \bigcup_{\rho \in R_i} X_i^\rho$, where $G = \bigcup_{\rho \in R_i} S_{y_i} N \rho$.

Set $I = \bigcup_{i \in I_0} \{X_i^\sigma \mid \sigma \in G\} = \bigcup_{i \in I_0} \{X_i^\rho \mid \rho \in R_i\}$. Since $R_i$ are finite, $I$ is finite and $G$ acts on $I$ from the right. For $i \in I_0$, $\sigma \in G$, and $\nu \in N$, (6b) implies $X_i^{\sigma \nu} = (X_i^{\sigma \nu \sigma^{-1}})^\sigma = X_i^\sigma$. Hence, the action of $G$ induces an action of $G/N$ on $I$.

Next define a map $\varphi\colon X \to I$ such that $\varphi(x) = X_i^\rho$ for all $i \in I_0$ and $\rho \in R_i$ and each $x \in X_i^\rho$. Since (6c) is a partition of $X$ into open-closed sets, $\varphi$ is surjective and continuous. Let $i \in I_0$, $y \in X_i$, and $\sigma \in G$. Write $\sigma = \tau\rho$ with $\tau \in S_{y_i}N$ and $\rho \in R_i$. Then $y^\sigma \in X_i^\rho$ (by (6c)) and $\varphi(y^\sigma) = X_i^\rho = X_i^\sigma$. Thus,

$$(7) \qquad\qquad \varphi(y^\sigma) = X_i^\sigma \qquad \text{for } y \in X_i,\ \sigma \in G.$$

It follows that $\varphi(x^\sigma) = \varphi(x)^{\varphi(\sigma)}$ for all $x \in X$ and $\sigma \in G$.

For each $X_i^\sigma \in I$ put $A_{X_i^\sigma} = \varphi(G_{y_i}^\sigma) = \varphi(G_{y_i})^{\varphi(\sigma)}$. This is a good difinition: If $X_i^\sigma = X_j^{\sigma'}$, then, by (6), $i = j$ and $\sigma' = \zeta\nu\sigma$ with $\zeta \in S_{y_i} \leq G_{y_i}$ and $\nu \in N$. Then $\varphi(G_{y_i}^{\sigma'}) = \varphi(G_{y_i}^\zeta)^{\varphi(\nu)\varphi(\sigma)} = \varphi(G_{y_i})^{\varphi(\sigma)}$.

We claim that $\varphi(G_x) \leq A_{\varphi(x)}$ for all $x \in X$. Indeed, there are $y \in X_i$ and $\sigma \in G$ such that $x = y^\sigma$. By (5), $y \in V_{y_i}$. Hence, $G_y \leq G_{y_i}N$, so $\varphi(G_y) \leq \varphi(G_{y_i})$. Therefore,

$$\varphi(G_x) = \varphi(G_y^\sigma) = \varphi(G_y)^{\varphi(\sigma)} \leq \varphi(G_{y_i})^{\varphi(\sigma)} = A_{X_i^\sigma} = A_{\varphi(x)}.$$

Finally, by (6b), the stabilizer of $X_i \in I$ in $G/N$ is contained in $S_{y_i}N/N = \varphi(S_{y_i})$. Therefore it is contained in $\varphi(G_{y_i}) = A_{X_i}$. Consequently, $(G/N, I, A_i)_{i \in I}$ is a finite group structure.

*Proof of (b):* Let $Y_0$ be a subset of $X_0$ with $X_0 \leq Y_0^G$ and $y^\sigma \neq y'$ for all distinct $y, y' \in Y_0$ and $\sigma \in G$. By Lemma 3.6 we may assume $\{y_i \mid i \in I_0\}$ contains $Y_0$. Write each $x \in X_0$ as $x = y^\sigma$ with $y \in Y_0$ and $\sigma \in G$. Then $y = y_i$ for some $i \in I_0$ and $\varphi(G_x) = \varphi(G_{y_i}^\sigma) = A_{X_i^\sigma} = A_{\varphi(x)}$.

*Proof of (c):* For each $y \in Y$ we may choose $V_y$ at the beginning of the proof of (a) such that $V_y$ is contained in the unique $Y_j$ which contains $y^N$. Since $G$ permutes the $Y_j$'s, each $Y_j$ can be partitioned as $Y_j = \bigcup_{(i,\rho) \in S_j} X_i^\rho$ with disjoint subsets $S_j$ of $\{(i,\rho) \mid i \in I_0,\ \rho \in R_i\}$. Therefore, $\varphi(Y_j) = \{A_{X_i^\rho} \mid (i,\rho) \in S_j\}$ are disjoint.

*Proof of (d):* Lemma 3.6 allows us to choose $I_0$ at the beginning of the proof of of (a) such that $\{1, \ldots, m\} \subseteq I_0$. By (7), $\varphi(x_i) = A_{X_i}$ belong then to distinct $A$-orbits. ∎

Lemma 3.7 has several consequences.

LEMMA 3.8: *Let $\varphi\colon \mathbf{G} \to \mathbf{A}$ a morphism of group structures. Then there are a morphism $\bar{\varphi}\colon \hat{\mathbf{A}} \to \mathbf{A}$ of group structures and an epimorphism $\hat{\varphi}\colon \mathbf{G} \to \hat{\mathbf{A}}$ satisfying $\varphi = \bar{\varphi} \circ \hat{\varphi}$.*

Proof: Let $\mathbf{G} = (G, X, G_x)_{x \in X}$ and $\mathbf{A} = (A, I, A_i)_{i \in I}$. For each $i \in I$ let $X_i = \varphi^{-1}(i)$. Then $G$ permutes the finite set $\{X_i \mid i \in I\}$. Hence, $G$ has an open normal subgroup $N$ such that $N \le \mathrm{Ker}(\varphi)$ and $X^\nu = X$ for each $\nu \in N$.

Set $\hat{A} = G/N$ and let $\hat{\varphi}\colon G \to \hat{A}$ be the quotient map. Use Lemma 3.7 to extend $\hat{\varphi}\colon G \to \hat{A}$ to an epimorphism $\hat{\varphi}\colon \mathbf{G} \to \hat{\mathbf{A}}$ such that $\hat{\mathbf{A}} = (\hat{A}, J, \hat{A}_j)_{j \in J}$ is finite and $\hat{\varphi}(X_i)$, $i \in I$, are disjoint.

Now define $\hat{\varphi}\colon \hat{A} \to A$ to be the map induced by $\varphi$. Define $\bar{\varphi}\colon J \to I$ by $\bar{\varphi}(j) = i$ for all $j \in \hat{\varphi}(X_i)$ and $i \in I$. Then $\hat{\varphi}\colon \hat{\mathbf{A}} \to \mathbf{A}$ is a morphism of finite group structures and $\varphi = \bar{\varphi} \circ \hat{\varphi}$.    ∎

LEMMA 3.9: *Let $\psi\colon \mathbf{H} \to \mathbf{G}$ be a cover of group structures with a finite kernel. Then there is a cartesian square of group structures*

(8)
$$
\begin{array}{ccc}
\mathbf{H} & \xrightarrow{\ \psi\ } & \mathbf{G} \\
{\scriptstyle \beta}\big\downarrow & & \big\downarrow{\scriptstyle \varphi} \\
\mathbf{B} & \xrightarrow{\ \alpha\ } & \mathbf{A}
\end{array}
$$

*in which $\mathbf{A}$ and $\mathbf{B}$ are finite and $\alpha$ is a cover.*

Proof: Let $\mathbf{G} = (G, X, G_x)_{x \in X}$ and $\mathbf{H} = (H, X, H_y)_{y \in Y}$. By Lemma 2.3, $\bigcup_{y \in Y} H_y$ is a closed subset of $H$. By assumption, $K = \mathrm{Ker}(\psi)$ is a finite group and $\bigcup_{y \in Y} H_y \cap (K \smallsetminus 1) = \emptyset$. Hence, $H$ has an open normal subgroup $N$ with $\left( \bigcup_{y \in Y} H_y \right) N \cap (K \smallsetminus 1) = \emptyset$. Thus, $N \cap K = 1$ and $H_y N \cap KN = N$ for each $y \in Y$.

Let $B = H/N$ and let $\beta\colon H \to B$ be the quotient map. Use Lemma 3.7 to complete $\beta$ to an epimorphism $\beta\colon \mathbf{H} \to \mathbf{B}$ with $\mathbf{B} = (B, J, B_j)_{j \in J}$ a finite group structure.

By Example 2.5, we may assume $\psi$ is the quotient map $\mathbf{H} \to \mathbf{H}/K$. Put $\mathbf{A} = (A, I, A_i)_{i \in I} = \mathbf{B}/\beta(K)$. Then let $\alpha\colon \mathbf{B} \to \mathbf{A}$ be the quotient map and $\varphi\colon \mathbf{G} \to \mathbf{A}$ the epimorphism which $\beta$ induces. This gives the commutative diagram (8). The assumption $H_y N \cap KN = N$ implies $B_j \cap \mathrm{Ker}(\alpha) = 1$ for each $j \in J$. Hence, by Example 2.5, $\alpha$ is a cover.

To prove (8) is cartesian, it suffices to prove that the unique morphism $\varepsilon\colon \mathbf{H} \to \mathbf{B} \times_{\mathbf{A}} \mathbf{G}$ induced by $\beta$ and $\psi$ is an isomorphism. Indeed, the group homomorphism $\varepsilon\colon H \to B \times_A G$ is an isomorphism [FrJ, Section 20.2]. We show $\varepsilon\colon Y \to J \times_I X$ is a bijection (hence, a homeomorphism): Let $(j, x) \in J \times_I X$. There is $y \in Y$ such that $\psi(y) = x$. As $\alpha(\beta(y)) = \varphi(x) = \alpha(j)$, there is a unique $b \in \mathrm{Ker}(\alpha) = \beta(K)$ with $\beta(y)^b = j$. Choose $k \in K$ with $\beta(k) = b$. Then $\beta(y^k) = j$ and $\psi(y^k) = \psi(y) = x$. Hence, $\varepsilon(y^k) = (j, x)$. Therefore, $\varepsilon$ is surjective.

Next let $y, y' \in Y$ with $\varepsilon(y) = \varepsilon(y')$. Then $\beta(y) = \beta(y')$ and $\psi(y) = \psi(y')$. As $\psi$ is a cover, there is $k \in K$ with $y^k = y'$. So, $\beta(y)^{\beta(k)} = \beta(y') = \beta(y)$. Hence, $\beta(k) \in \beta(K) \cap S_{\beta(y)} \le \beta(K) \cap B_{\beta(y)} \le \beta(K) \cap B_{\beta(y)} = 1$ (because $\mathrm{Ker}(\alpha) = \beta(K)$ and $\alpha$ is a cover). Thus, $k \in N \cap K = 1$. Therefore, $y = y'$. We conclude that $\varepsilon$ is injective.

Since $\alpha$ is a cover, Lemma 2.12(c) implies that the projection $\psi'\colon \mathbf{B} \times_{\mathbf{A}} \mathbf{G} \to \mathbf{G}$ is a cover. By assumption $\psi$ is a cover. Hence, for each $y \in Y$ and $(j, x) = \varepsilon(y)$ both $\psi'\colon B_j \times_A G_x \to G_x$ and $\psi\colon H_y \to G_x$ are isomorphisms. Since $\psi = \psi \circ \varepsilon'$, so is $\varepsilon\colon H_y \to B_j \times_A G_x$. This concludes the proof that (8) is cartesian. ∎

## 4. Projective Group Structures

The notion "projective group structure" which we introduce here replaces the notion "relatively projective group" of [HaJ, Def. 4.2], also called "strongly relatively projective" in [Pop, p. 4]. A projective group structure is one of the two main objects which we put in duality in this work, the other one being "field-valuation structure with the block approximation condition" (Section 12).

Let $\mathbf{G}$ be a group structure. An **embedding problem** for $\mathbf{G}$ is a pair

$$(1) \qquad\qquad \big(\varphi\colon \mathbf{G} \to \mathbf{A},\ \alpha\colon \mathbf{B} \to \mathbf{A}\big)$$

of morphisms of group structures in which $\alpha$ is a cover. A **solution** of (1) is a morphism $\gamma\colon \mathbf{G} \to \mathbf{B}$ with $\alpha \circ \gamma = \varphi$. The embedding problem is **finite** if $\mathbf{B}$ is finite. We say $\mathbf{G}$ is **projective**, if every finite embedding problem for $\mathbf{G}$ has a solution.

LEMMA 4.1: *Let $\mathbf{G}$ be a group structure. Suppose every finite embedding problem (1) for $\mathbf{G}$ where $\varphi$ is an epimorphism is solvable. Then $\mathbf{G}$ is projective.*

*Proof:* Lemma 3.8 gives a morphism $\bar{\varphi}\colon \hat{\mathbf{A}} \to \mathbf{A}$ of finite group structures and an epimorphism $\hat{\varphi}\colon \mathbf{G} \to \hat{\mathbf{A}}$ satisfying $\varphi = \bar{\varphi} \circ \hat{\varphi}$. Set $\hat{B} = \mathbf{B} \times_{\mathbf{A}} \hat{\mathbf{A}}$. Let $\beta\colon \hat{\mathbf{B}} \to \mathbf{B}$ and $\hat{\alpha}\colon \hat{\mathbf{B}} \to \hat{\mathbf{A}}$ be the projections maps. By Lemma 2.12, $\hat{\alpha}\colon \hat{B} \to \hat{A}$ is a cover. Hence, $(\hat{\varphi}\colon \mathbf{G} \to \hat{\mathbf{A}},\ \hat{\alpha}\colon \hat{\mathbf{B}} \to \hat{\mathbf{A}})$ is a finite embedding problem. By assumption, there is a morphism $\hat{\gamma}\colon \mathbf{G} \to \hat{\mathbf{B}}$ with $\hat{\alpha} \circ \hat{\gamma} = \hat{\varphi}$. Then $\gamma = \beta \circ \hat{\gamma}$ is a solution of (1). Therefore, $\mathbf{G}$ is projective. ∎

Given a profinite group $G$, Gruenberg proved that if every finite embedding problem for $G$ is solvable, then every embedding problem for $G$ is solvable [FrJ, Lemma 20.8]. Gruenberg's proof goes through in the category of group structures almost verbatim.

PROPOSITION 4.2: *Let $\mathbf{G}$ be a projective group structure. Then every embedding problem for $\mathbf{G}$ has a solution.*

*Proof:* Let (1) be an embedding problem for $\mathbf{G}$. Put $K = \mathrm{Ker}(\alpha)$.

PART A: *Suppose $K$ is finite.* Lemma 3.9 gives a cartesian square of group structures

(without the dashed morphisms) in which $\bar{\mathbf{B}}$ and $\bar{\mathbf{A}}$ are finite, $\bar{\alpha}$ is a cover, and $\mathbf{B} = \bar{\mathbf{B}} \times_{\bar{\mathbf{A}}} \mathbf{A}$. Put $\bar{\varphi} = \varphi' \circ \varphi$. Then $(\bar{\varphi} \colon \mathbf{G} \to \bar{\mathbf{A}}, \bar{\alpha} \colon \bar{\mathbf{B}} \to \bar{\mathbf{A}})$ is a finite embedding problem for $\mathbf{G}$. By assumption there is a morphism $\bar{\gamma} \colon \mathbf{G} \to \bar{\mathbf{B}}$ with $\bar{\alpha} \circ \bar{\gamma} = \bar{\varphi}$. Hence, there is a morphism $\gamma \colon \mathbf{G} \to \mathbf{B}$ with $\alpha \circ \gamma = \varphi$ and $\psi' \circ \gamma = \bar{\gamma}$ (Definition 2.10). In particular $\gamma$ solves embedding problem (1).

PART B: *Application of Zorn's lemma.* Suppose (1) is an arbitrary embedding problem for $\mathbf{G}$. By Example 2.5 we may assume $\mathbf{A} = \mathbf{B}/K$ and $\alpha$ is the quotient map. For each normal subgroup $L$ of $B$ contained in $K$ let $\alpha_L \colon \mathbf{B}/L \to \mathbf{A}$ be the quotient map $\mathbf{B}/L \to (\mathbf{B}/L)/(K/L)$. Then, $\alpha_L$ is a cover (Example 2.5) and

$$(2) \qquad\qquad (\varphi \colon \mathbf{G} \to \mathbf{A}, \ \alpha_L \colon \mathbf{B}/L \to \mathbf{A}).$$

is an embedding problem for $\mathbf{G}$. Let $\Lambda$ be the set of pairs $(L, \gamma)$ where $L$ is a closed normal subgroup of $B$ contained in $K$ and $\gamma$ is a solution of (2). The pair $(K, \varphi)$ belongs to $\Lambda$. Partially order $\Lambda$ by $(L', \gamma') \leq (L, \gamma)$ if $L' \leq L$ and $\alpha_{L', L} \circ \gamma' = \gamma$. Here $\alpha_{L', L} \colon \mathbf{B}/L' \to \mathbf{B}/L$ is the cover $\mathbf{B}/L' \to (\mathbf{B}/L')/(L/L')$.

Suppose $\Lambda_0 = \{(L_j, \gamma_j) \mid j \in J\}$ is a descending chain in $\Lambda$. Then $\varprojlim \mathbf{B}/L_j = \mathbf{B}/L$ with $L = \bigcap_{j \in J} L_j$ (Lemma 2.8). The $\gamma_j$'s define a morphism $\gamma \colon \mathbf{G} \to \mathbf{B}/L$ with $\alpha_{L, L_j} \circ \gamma = \gamma_j$ for each $j$. Thus, $(L, \gamma)$ is a lower bound to $\Lambda_0$.

Zorn's lemma gives a minimal element $(L, \gamma)$ for $\Lambda$. It suffices to prove that $L = 1$.

Assume $L \neq 1$. Then $B$ has an open normal subgroup $N$ with $L \nleq N$. Thus, $L' = N \cap L$ is a proper open subgroup of $L$ which is normal in $B$. Then $(\gamma \colon \mathbf{G} \to \mathbf{B}/L, \ \alpha_{L', L} \colon \mathbf{B}/L' \to \mathbf{B}/L)$ is an embedding problem for $\mathbf{G}$. Its kernel $\mathrm{Ker}(\alpha_{L', L}) =$ $L/L'$ is a finite group. Hence, by Part A, it has a solution $\gamma'$. The pair, $(L', \gamma')$ is an element of $\Lambda$ which is strictly smaller than $(L, \gamma)$. This contradiction to the minimality of $(L, \gamma)$ proves that $L = 1$, as desired. ∎

COROLLARY 4.3: *Let* $\psi \colon \mathbf{H} \to \mathbf{G}$ *be a cover of group structures. Suppose* $\mathbf{G}$ *is projective. Then* $\mathbf{H}$ *has a sub-group-structure* $\mathbf{H}'$ *which* $\psi$ *maps isomorphically onto* $\mathbf{G}$.

*Proof:* Suppose $\mathbf{G} = (G, X, G_x)_{x \in X}$ and $\mathbf{H} = (H, Y, H_y)_{y \in Y}$. Proposition 4.2 gives a morphism $\gamma \colon \mathbf{G} \to \mathbf{H}$ with $\psi \circ \gamma = \mathrm{id}_{\mathbf{G}}$. Let $H' = \gamma(G)$ and $Y' = \gamma(X)$. Then $\psi \colon H' \to G$ is an isomorphism and $\psi \colon Y' \to X$ is a homeomorphism. Next, let $x \in X$ and $y' = \gamma(x)$. Then $\psi(y') = x$ and $\gamma(G_x) \leq H_{y'}$. As a cover, $\psi$ maps both $H_{y'}$ and $\gamma(G_x)$ isomorphically onto $G_x$. Hence, $\gamma(G_x) = H_{y'}$. In particular, $H_{y'} \leq H'$. So, $y' \mapsto H_{y'}$ is a continuous map of $Y'$ into $\mathrm{Subgr}(H')$ (Remark 2.6). Thus, $\mathbf{H}' = (H', Y', H_{y'})_{y' \in Y}$ is a sub-group-structure of $\mathbf{H}$ which $\psi$ maps isomorphically onto $\mathbf{G}$. ∎

We shall have several occasions to use the following result of Herfort and Ribes.

PROPOSITION 4.4 ([HeR, Prop. 2 and Thm. B']): *Let* $G = \coprod_{i \in I} G_i$ *be the free profinite product of finitely many profinite groups* $G_i$. *Then* $G_i^g \cap G_j \neq 1$ *implies* $i = j$ *and* $g \in G_i$.

LEMMA 4.5: *Let* $\mathbf{A} = (A, I, A_i)_{i \in I}$ *be a group structure,* $\alpha \colon B \to A$ *an epimorphism of profinite groups, and* $I_0$ *be a finite system of representatives of the* $A$-*orbits of* $I$. *For each* $i \in I_0$ *let* $B_i$ *be a closed subgroup of* $B$ *which* $\alpha$ *maps isomorphically onto* $A_i$. *Then* $\alpha$ *extends to a cover* $\alpha \colon \mathbf{B} \to \mathbf{A}$ *where* $\mathbf{B} = (B, J, B_j)_{j \in J}$ *is a group structure. Moreover, there is a map* $\alpha' \colon I_0 \to J$ *such that* $J = \alpha'(I_0)^B$, $\alpha(\alpha'(i)) = i$, *and* $B_i = B_{\alpha'(i)}$ *for each* $i \in I_0$.

*Proof:* Consider $i \in I_0$. Then $S_i = \{a \in A \mid i^a = i\}$ is a subgroup of $A_i$. Hence, $T_i = \alpha^{-1}(S_i) \cap B_i$ is a subgroup of $B_i$ which $\alpha$ maps bijectively onto $S_i$. Also, the set $\{(i, T_i b) \mid b \in B\}$ bijectively corresponds to the profinite quotient space $B/T_i$. So, $J = \bigcup_{i \in I_0} \{(i, T_i b) \mid b \in B\}$ is a profinite space. The rule $(i, T_i b)^{b'} = (i, T_i b b')$ defines a continuous action of $B$ on $J$. For each $j = (i, T_i b) \in J$ let $B_j = B_i^b$. Then $j \mapsto B_j$ is a strictly continuous, (hence also étale continuous) map from $J$ into $\mathrm{Subgr}(B)$.

Now suppose $(i, T_i b)^{b'} = (i, T_i b)$. Then $T_i bb' = T_i b$. Hence, $b' \in T_i^b \leq B_i^b$. Therefore, $\mathbf{B} = (B, J, B_j)_{j \in J}$ is a group structure.

Next define a map $\alpha\colon J \to I$ by $\alpha(i, T_i b) = i^{\alpha(b)}$. If $\alpha(i', T_{i'} b') = \alpha(i, T_i b)$, then $i^{\alpha(b)} = (i')^{\alpha(b')}$. As $i, i' \in I_0$, this implies $i = i'$ and $\alpha(b) = s_i \alpha(b')$ for some $s_i \in S_i$. Let $t_i$ be the element of $T_i$ with $\alpha(t_i) = s_i$. Then there is $k \in \mathrm{Ker}(\alpha)$ with $b = t_i b' k$. Hence, $(i, T_i b) = (i, T_i t_i b' k) = (i', T_{i'} b')^k$. It follows, $\alpha\colon \mathbf{B} \to \mathbf{A}$ is a cover.

Finally define a map $\alpha'\colon I_0 \to J$ by $\alpha'(i) = (i, T_i)$. Then $(i, T_i b) = \alpha'(i)^b$ for each $i \in I_0$ and $b \in B$, so $J = \alpha'(I_0)^B$. Also, $\alpha(\alpha'(i)) = \alpha(i, T_i) = i$ and $B_{\alpha'(i)} = B_{(i, T_i)} = B_i$ for each $i \in I_0$. ∎

The assumption on a group structure $\mathbf{G}$ to be projective poses some restrictions on $\mathbf{G}$:

PROPOSITION 4.6: Let $\mathbf{G} = (G, X, G_x)_{x \in X}$ be a projective group structure.

(a) Let $x, y \in X$ with $G_x \cap G_y \neq 1$. Then $y = x^g$ for some $g \in G_x$. Hence, $G_x = G_y$.

(b) Let $x \in X$ with $G_x \neq 1$. Then $G_x$ is its own normalizer in $G$.

(c) Suppose $1 \notin \mathrm{StrictClosure}\{G_x \mid x \in X\}$ and $G_x = S_x$ for each $x \in X$. Then $\mathbf{G}$ is a proper structure.

Proof of (a): There is an epimorphism $\bar{\varphi}\colon G \to \bar{A}$ with $\bar{A}$ finite and $\bar{\varphi}(G_x \cap G_y) \neq 1$. Consider an arbitrary epimorphism $\varphi\colon G \to A$ with $A$ finite and $\mathrm{Ker}(\varphi) \leq \mathrm{Ker}(\bar{\varphi})$. Then $\varphi(G_x \cap G_y) \neq 1$.

Use Lemma 3.7 to complete $\varphi$ to an epimorphism $\varphi\colon \mathbf{G} \to \mathbf{A}$ of group structures with $\mathbf{A} = (A, I, A_i)_{i \in I}$ finite. In addition, $\varphi(G_x) = A_{\varphi(x)}$ and $\varphi(x), \varphi(y)$ are not in the same $A$-orbit if $x, y$ are not in the same $G$-orbit.

Assume without loss that $I$ does not contain the symbol 0. Choose a system of representatives $I_0$ for the $A$-orbits of $I$ which does not contain the symbol 0. Put $I'_0 = \{0\} \cup I_0$ and $A_0 = A$. For each $i \in I'_0$ choose an isomorphic copy $B_i$ of $A_i$ and an isomorphism $\alpha_i\colon B_i \to A_i$.

Now consider the free profinite product $B = \coprod_{i \in I'_0} B_i$. Let $\alpha\colon B \to A$ be the unique epimorphism with $\alpha|_{B_i} = \alpha_i$, $i \in I'_0$. Lemma 4.5 extends $B$ to a group structure $\mathbf{B} = (B, J, B_j)_{j \in J}$ with $\{B_j \mid j \in J\} = \{B_i^b \mid i \in I_0\}$ and $\alpha$ to a cover $\alpha\colon \mathbf{B} \to \mathbf{A}$.

Moreover, there is a map $\alpha'\colon I_0 \to J$ such that $J = \alpha'(I_0)^B$, $\alpha(\alpha'(i)) = i$, and $B_i = B_{\alpha'(i)}$ for each $i \in I_0$.

Since $\mathbf{G}$ is projective, Proposition 4.2 gives a morphism $\gamma\colon \mathbf{G} \to \mathbf{B}$ with $\alpha \circ \gamma = \varphi$. In particular $\alpha(\gamma(G_x \cap G_y)) = \varphi(G_x \cap G_y) \neq 1$. Hence, $1 < \gamma(G_x \cap G_y) \leq \gamma(G_x) \cap \gamma(G_y) \leq B_{\gamma(x)} \cap B_{\gamma(y)}$. Write $\gamma(x) = \alpha'(i)^b$ and $\gamma(y) = \alpha'(i')^{b'}$ with $i, i' \in I_0$ and $b, b' \in B$. Then $B_i^b \cap B_{i'}^{b'} = B_{\alpha'(i)}^b \cap B_{\alpha'(i')}^{b'} = B_{\gamma(x)} \cap B_{\gamma(y)} \neq 1$. By Proposition 4.4, $i = i'$. Hence $\gamma(x)$ and $\gamma(y)$ are in the same $B$-orbit. Therefore, $\varphi(x)$ and $\varphi(y)$ are in the same $A$-orbit. The choice of $\varphi$ gives $g \in G$ with $x^g = y$.

Now $B_{\gamma(x)}^{\gamma(g)} = B_{\gamma(x^g)} = B_{\gamma(y)}$. By the preceding paragraph, $B_{\gamma(x)} \cap B_{\gamma(x)}^{\gamma(g)} = B_{\gamma(x)} \cap B_{\gamma(y)} \neq 1$. Since $B_{\gamma(x)}$ is conjugate to $B_i$ for some $i \in I_0$, Proposition 4.4 implies $\gamma(g) \in B_{\gamma(x)}$. Hence, $\varphi(g) \in A_{\varphi(x)}$. Since this relation holds for all $\varphi$ with $\mathrm{Ker}(\varphi) \leq \mathrm{Ker}(\bar{\varphi})$, we have $g \in G_x$, as desired.

*Proof of (b):* Suppose $G_x \neq 1$. Consider $g \in G$ with $G_x^g = G_x$. By (a), there is $a \in G_x$ with $x^{ga} = x$. Then $ga \in G_x$. Hence, $g \in G_x$.

*Proof of (c):* Suppose $G_x \neq G_y$ for some $x, y \in X$, then there is $g \in G_x$ with $y = x^g$. Hence, by assumption, $y = x$. Thus, the forgetful map $\delta_{\mathbf{G}}$ is an étale continuous bijection of $X$ onto $\mathcal{G} = \{G_x \mid x \in X\}$. By Corollary 1.4, $\mathcal{G}$ is étale Hausdorff. Since $X$ is compact, $\delta_{\mathbf{G}}$ is an étale homeomorphism. It follows, $\mathbf{G}$ is proper. ∎

*Example 4.7: Projective structures.*

(a) Projective group. Let $G$ be a profinite group and $X$ the empty space. Then $\mathbf{G} = (G, X,\ )$ is a projective proper group structure if and only if $G$ is a projective group.

(b) Trivial stabilizers. Let $G$ be an arbitrary profinite group. Put $X = G$. Then $X$ is a profinite space and $G$ acts continuously on $X$ by multiplication from the right. In particular, $S_x = 1$ for each $x \in X$. For each $x \in X$ put $G_x = G$. Then $\mathbf{G} = (G, X, G_x)_{x \in X}$ is a projective group structure.

Indeed, let $(\varphi\colon \mathbf{G} \to \mathbf{A},\ \alpha\colon \mathbf{B} \to \mathbf{A})$ be a finite embedding problem for $\mathbf{G}$. Let $i \in I$ and $j \in J$ be elements with $\varphi(1) = i$ and $\alpha(j) = i$. Then $G_1 = G$, so $\varphi(G) \leq A_i$. Also, $\alpha\colon B_j \to A_i$ is an isomorphism. Hence, $\gamma_{\mathrm{g}} = (\alpha|_{B_j})^{-1} \circ \varphi$ is a homomorphism from $G$ to $B$ satisfying $\alpha \circ \gamma_{\mathrm{g}} = \varphi$. Define $\gamma_{\mathrm{s}} \colon X \to J$ by $\gamma_{\mathrm{s}}(x) = j^{\gamma_{\mathrm{g}}(x)}$. Then $\gamma = (\gamma_{\mathrm{g}}, \gamma_{\mathrm{s}})$ is a solution of the embedding problem, as desired.

If $G$ is nontrivial, then $1 \notin \mathrm{StrictClosure}\{G_x \mid x \in X\}$ but $\mathbf{G}$ is not proper. It follows that the assumption $S_x \neq G_x$ in Proposition 4.6(c) (which is violated in our example) is necessary.

(c) Free products of finitely many profinite groups.

Let $K$ be a finite set and $K_0$ a subset. For each $k \in K$ let $G_k$ be a nontrivial profinite group. Write $G = \coprod_{k \in K} G_k$ for the free product of the $G_k$'s. For each $k \in K$ the orbit $\mathcal{G}_k = \{G_k^g \mid g \in G\}$ of $G_i$ under conjugation is a strictly closed subset of $\mathrm{Subgr}(G)$. Hence, $\mathcal{G} = \bigcup_{k \in K_0} \mathcal{G}_k$ is a strictly profinite subspace of $\mathrm{Subgr}(G)$, so strictly closed. In particular, $1 \notin \mathrm{StrictClosure}(\mathcal{G})$. By Proposition 4.4, $H \cap H' = 1$ for all distinct $H, H' \in \mathcal{G}$. It follows from Corollary 1.4 that $\mathcal{G}$ is étale Hausdorff.

Choose a homeomorphic copy $X$ of $\mathcal{G}$ with the strict topology and a strict homeomorphism $\delta \colon X \to \mathcal{G}$. Since the strict topology of $\mathrm{Subgr}(G)$ is finer than its étale topology, $\delta$ is étale continuous. Since $\mathcal{G}$ is étale Hausdorff, $\delta$ is an étale homeomorphism. For each $x \in X$ let $G_x = \delta(x)$. By Proposition 4.4, each $H \in \mathcal{G}$ is its own normalizer in $G$. Thus, in the terminology of Section 2, $S_x = G_x$. Therefore, $\mathbf{G} = (G, X, G_x)_{x \in X}$ is a proper group structure.

We prove $\mathbf{G}$ is projective. To this end consider finite group structures $\mathbf{A} = (A, I, A_i)_{i \in I}$ and $\mathbf{B} = (B, J, B_j)_{j \in J}$, a cover $\alpha \colon \mathbf{B} \to \mathbf{A}$, and an epimorphism $\varphi \colon \mathbf{G} \to \mathbf{A}$. By Lemma 4.1 it suffices to find a morphism $\gamma \colon \mathbf{G} \to \mathbf{B}$ with $\gamma \circ \alpha = \varphi$.

Choose a map $\alpha' \colon I \to J$ with $\alpha(\alpha'(i)) = i$ for each $i \in I$. Now consider $k \in K_0$. Let $i = \varphi(k)$ and $j = \alpha'(i)$. Then $\alpha \colon B_j \to A_i$ is an isomorphism. Hence, $\gamma_k = (\alpha|_{B_j})^{-1} \circ (\varphi|_{G_k})$ is an epimorphism of $G_k$ onto $B_j$ satisfying $\alpha \circ \gamma_k = \varphi|_{G_k}$. The basic property of free products gives a homomorphism $\gamma \colon G \to B$ whose restriction to each $G_k$ is $\gamma_k$. In particular, $\alpha \circ \gamma = \varphi$. Together with the map $\gamma = \alpha' \circ \varphi$ from $X$ to $B$, $\gamma \colon \mathbf{G} \to \mathbf{B}$ is a morphism satisfying $\alpha \circ \gamma = \varphi$, as desired. ∎

## 5. Special Covers

As in Lemma 4.5 we consider a group structure $\mathbf{G} = (G, X, G_x)_{x \in X}$ and an epimorphism of profinite groups $\pi \colon H \to G$. In contrast to Lemma 4.5, we do not assume that $X$ has only finitely many $G$-orbits. Nor do we assume that $X$ has a fundamental domain (beginning of Section 3). Nevertheless, we are able to extend $\pi \colon H \to G$ to a cover $\pi \colon \mathbf{H} \to \mathbf{G}$ in special cases described in Lemma 5.1 below. They occur three times in Galois-theoretic set-ups (Lemma 14.2 and twice in Lemma 15.1).

LEMMA 5.1: *Let $\mathbf{G} = (G, X, G_x)_{x \in X}$ be a group structure and $(G_i, X_i, R_i)_{i \in I_0}$ a special partition of $\mathbf{G}$ (Definition 3.5). Let $\pi \colon H \to G$ be an epimorphism of profinite groups. For each $i \in I_0$ let $H_i$ be a subgroup of $H$ which $\pi$ maps isomorphically onto $G_i$.*

*Then $H$ extends to a profinite group structure $\mathbf{H} = (H, Y, H_y)_{y \in Y}$ and $\pi$ extends to a cover $\mathbf{H} \to \mathbf{G}$. Moreover, for each $i \in I_0$ there is a subspace $Y_i$ of $Y$ such that $\pi \colon Y_i \to X_i$ is a homeomorphism, $H_y \leq H_i$ for each $y \in Y_i$, and $\bigcup_{i \in I_0} Y_i^H = Y$.*

*If, in addition, $\mathbf{G}$ is proper and*

(1) $H_i^\kappa \cap H_i = 1$ *for all $\kappa \in \mathrm{Ker}(\pi)$ with $\kappa \neq 1$ and each $i \in I_0$,*

*then $\mathbf{H}$ is proper.*

*Proof:* The proof has four parts.

PART A: *The space $\hat{Y}$.* Let $X' = \bigcup_{i \in I_0} X_i$. This is a profinite space and hence so is the product $\hat{Y} = X' \times H$. The group $H$ acts continuously on $\hat{Y}$ by $(x, h)^\eta = (x, h\eta)$ and there is a continuous map $\hat{\pi} \colon \hat{Y} \to X$ defined by $\pi(x, h) = x^{\pi(h)}$. Since $X = \bigcup_{i \in I_0} \bigcup_{\rho \in R_i} X_i^\rho$, this map is surjective.

For each $y = (x, h) \in \hat{Y}$ define a subgroup $H_y$ of $H$ in the following way. There is a unique $i \in I_0$ with $x \in X_i$. Then $G_x \leq G_i$. Let $H_x$ be the unique subgroup of $H_i$ satisfying $\pi(H_x) = G_x$. Put $H_y = H_x^h$. Then

(2a) $\hat{\pi}(y^\eta) = \hat{\pi}(y)^{\pi(\eta)}$ for all $y \in \hat{Y}$ and $\eta \in H$,

(2b) $H_y^\eta = H_{y^\eta}$ for all $y \in \hat{Y}$ and $\eta \in H$, and

(2c) $\pi \colon H_y \to G_{\hat{\pi}(y)}$ is an isomorphism, $y \in \hat{Y}$.

CLAIM A1: *The map $\hat{Y} \to \mathrm{Subgr}(H)$ defined by $y \mapsto H_y$ is étale continuous.* It suffices to show that the map $X_i \to \mathrm{Subgr}(H)$ defined by $x \mapsto H_x$ is étale continuous.

By Remark 2.6 we have to show that the corresponding map $X_i \to \mathrm{Subgr}(H_i)$ is étale continuous. Now, by assumption, the map $X \to \mathrm{Subgr}(G)$ given by $x \mapsto G_x$ is étale continuous. Hence, by Remark 2.6, the corresponding map $X_i \to \mathrm{Subgr}(G_i)$ is continuous. As $G_i$ is isomorphic to $H_i$, we get our claim.

EQUIVALENCE RELATION: Define an equivalence relation $\equiv$ on $\hat{Y}$ as follows. Let $(x_1, h_1) \equiv (x_2, h_2)$ if there is a (unique) $i \in I_0$ with $x_1, x_2 \in X_i$, $H_i h_1 = H_i h_2$, and $x_1^{\pi(h_1)} = x_2^{\pi(h_2)}$. This relation satisfies the following rules:

(3a) If $y_1 \equiv y_2$, then $\hat{\pi}(y_1) = \hat{\pi}(y_2)$.

(3b) If $y_1 \equiv y_2$, then $H_{y_1} = H_{y_2}$.

(3c) If $y_1 \equiv y_2$ and $\eta \in H$, then $y_1^\eta \equiv y_2^\eta$.

   Let $K = \mathrm{Ker}(\pi)$.

CLAIM A2: $\hat{\pi}(x_1, h_1) = \hat{\pi}(x_2, h_2)$ *if and only if there is* $k \in K$ *with* $(x_2, h_2) \equiv (x_1, h_1 k)$.

   Indeed, let $i \in I_0$ with $x_1 \in X_i$. If $\hat{\pi}(x_1, h_1) = \hat{\pi}(x_2, h_2)$, then $x_1^{\pi(h_1 h_2^{-1})} = x_2$. Hence, by (2d) and (2f) of Section 3, $\pi(h_1 h_2^{-1}) \in G_i = \pi(H_i)$. Therefore, there is $k_0 \in K$ with $h_1 h_2^{-1} k_0 \in H_i$. Then $k = h_2^{-1} k_0 h_2 \in K$, $H_i h_1 k = H_i h_2$, and $x_1^{\pi(h_1 k)} = x_1^{\pi(h_1)} = x_2^{\pi(h_2)}$. Consequently $(x_2, h_2) \equiv (x_1, h_1 k)$.

   Conversely, if $(x_2, h_2) \equiv (x_1, h_1 k)$, then $x_2^{\pi(h_2)} = x_1^{\pi(h_1 k)} = x_1^{\pi(h_1)}$, so $\hat{\pi}(x_2, h_2) = \hat{\pi}(x_1, h_1)$.

PART B: *The quotient space* $Y$. Let $Y$ be the quotient space of $\hat{Y}$ modulo $\equiv$. For each $y \in Y$ let $H_y = \delta_H(y)$. By (3a), $\hat{\pi} \colon \hat{Y} \to X$ induces a continuous surjection $\pi \colon Y \to X$. By (3b), the map $\hat{Y} \to \mathrm{Subgr}(H)$ induces a étale continuous map $\delta_H \colon Y \to \mathrm{Subgr}(H)$. By (3c), the $H$-action on $\hat{Y}$ induces a continuous action of $H$ on $Y$. By (2),

(4a) $\pi(y^\eta) = \pi(y)^{\pi(\eta)}$ for all $y \in Y$ and $\eta \in H$,

(4b) $H_y^\eta = H_{y^\eta}$ for all $y \in Y$ and $\eta \in H$, and

(4c) $\pi \colon H_y \to G_{\pi(y)}$ is an isomorphism, for each $y \in Y$.

Finally, by Claim A2,

(5) $\pi(y_1) = \pi(y_2)$ if and only if there is $k \in K$ such that $y_2 = y_1^k$.

CLAIM B1: $Y$ *is a profinite space.* Indeed, $\hat{Y}$ is compact, hence so is $Y$. Consider unequivalent $y_1, y_2 \in \hat{Y}$. It suffices to produce an open-closed neighborhood $U$ of $y_1$ which is closed under $\equiv$ and does not contain $y_2$.

If $\hat{\pi}(y_1) \neq \hat{\pi}(y_2)$, we choose an open-closed neighborhood $V$ of $\hat{\pi}(y_1)$ in $X$ which does not contain $\hat{\pi}(y_2)$. Then $U = \hat{\pi}^{-1}(V)$ has the required property.

If $\hat{\pi}(y_1) = \hat{\pi}(y_2)$, we use Claim A2 to replace $y_2$ by an equivalent element of $\hat{Y}$ to assume that $y_1 = (x_1, h_1)$, $y_2 = (x_1, h_1 k)$, where $1 \neq k \in K$. Let $i \in I_0$ such that $x_1 \in X_i$. Then $H_i \cap K = 1$, so $h_1 k h_1^{-1} \notin H_i$. There is an open subgroup $H_i'$ which contains $H_i$ and $h_1 k h_1^{-1} \notin H_i'$. Let $U = X_i \times H_i' h_1$. Then $(x_1, h_1) \in U$ but $h_1 k \notin H_i' h_1$, so $(x_1, h_1 k) \notin U$. Clearly $U$ is an open-closed subset of $\hat{Y}$ closed under $\equiv$.

CLAIM B2: *The stabilizer $S_y$ of each $y \in Y$ is contained in $H_y$.*

Indeed, let $y$ be represented by $(x, h) \in \hat{Y}$. Let $i \in I_0$ with $x \in X_i$. Let $\eta \in H$. Then

$$y^\eta = y \Longrightarrow (x, h)^\eta \equiv (x, h)$$
$$\Longrightarrow (x, h\eta) \equiv (x, h)$$
$$\Longrightarrow H_i h \eta = H_i h \text{ and } x^{\pi(h\eta)} = x^{\pi(h)}$$
$$\Longrightarrow \eta \in H_i^h \text{ and } \pi(\eta) \in S_{x^{\pi(h)}}.$$

Hence, $S_y \leq H_i^\eta$ and $\pi(S_y) \leq S_{x^{\pi(h)}} \leq G_{x^{\pi(h)}} \leq G_i^{\pi(h)}$. In addition, $\pi$ maps $H_i^h$ isomorphically onto $G_i^{\pi(h)}$ and $\pi(H_y) = G_{x^{\pi(h)}}$. Therefore, $S_y \leq H_y$, as claimed.

Claim B2 completes the proof that $\mathbf{H} = (H, Y, H_y)_{y \in Y}$ is a group structure and $\pi \colon \mathbf{H} \to \mathbf{G}$ is a cover.

PART C: *The spaces $Y_i$.* For each $i \in I_0$ let $Y_i$ be the image of $X_i \times 1$ in $Y$. Then, $\pi$ maps $Y_i$ homeomorphically onto $X_i$. By definition, $H_y \leq H_i$ for each $y \in Y_i$. By the assumption on $X$ we have $X = \bigcup_{i \in I_0} X_i^G$. Since $\pi \colon \mathbf{H} \to \mathbf{G}$ is a cover and $\pi(Y_i) = X_i$, we have $Y = \bigcup_{i \in I_0} Y_i^H$.

PART D: $\mathbf{H}$ *is proper under the assumption that* $\mathbf{G}$ *is proper and (1) holds.*

Indeed, let $\mathcal{H} = \{H_y \mid y \in Y\}$ and $\mathcal{G} = \{G_x \mid x \in X\}$. Since $\pi(\mathcal{H}) = \mathcal{G}$, we have $\pi(\text{StrictClosure}(\mathcal{H})) \subseteq \text{StrictClosure}(\mathcal{G})$. Since 1 is not in $\text{StrictClosure}(\mathcal{G})$, it is not in $\text{StrictClosure}(\mathcal{H})$.

Let $y_1, y_2 \in Y$ be distinct. We prove that $H_{y_1}, H_{y_2}$ are distinct and can be separated in the étale topology on $\text{Subgr}(H)$.

First suppose $\pi(y_1) \neq \pi(y_2)$. Then $G_{\pi(y_1)} \neq G_{\pi(y_2)}$. Since $\mathcal{G}$ is étale profinite, there are open subgroups $E_1, E_2$ of $G$ with $\pi(H_{y_i}) = G_{\pi(y_i)} \leq E_i$, $i = 1, 2$, and $\mathcal{G} \cap \mathrm{Subgr}(E_1) \cap \mathrm{Subgr}(E_2) = \emptyset$. Then $F_1 = \pi^{-1}(E_1)$ and $F_2 = \pi^{-1}(E_2)$ are open subgroups of $H$, $H_{y_1} \leq F_1$, $H_{y_2} \leq F_2$, and $\mathcal{H} \cap \mathrm{Subgr}(F_1) \cap \mathrm{Subgr}(F_2) = \emptyset$.

Now suppose $\pi(y_1) = \pi(y_2)$. Since $\pi$ is a cover, there is $\kappa \in K$ with $y_2 = y_1^\kappa$. As $y_1 \neq y_2$, we have $\kappa \neq 1$. Let $y_1$ be represented by $(x, h) \in \hat{Y}$, with $x \in X_i$, where $i \in I_0$, and $h \in H$. Then $H_{y_1} = H_x^h \leq H_i^h$ and $H_{y_2} = H_x^{h\kappa} \leq H_i^{h\kappa}$. By (1), $H_i^{h\kappa h^{-1}} \cap H_i = 1$, that is, $H_i^{h\kappa} \cap H_i^h = 1$. Hence, $H_{y_1} \cap H_{y_2} = 1$. By Corollary 1.4(a), $H_{y_1}$ and $H_{y_2}$ can be separated by the étale topology of $\mathcal{H}$.

It follows that $\mathcal{H}$ is étale Hausdorff and the étale continuous map $\delta_H : Y \to \mathrm{Subgr}(H)$ is bijective. By Claim B1, $Y$ is compact. Hence, $\delta_H$ is a homeomorphism. Consequently, $\mathbf{H}$ is proper. ∎

## 6. Unirationally Closed Fields

Galois duality naturally tranlates each groups structures $\mathbf{G} = (G, X, G_x)_{x \in X}$ with $G = \mathrm{Gal}(K)$ and $K$ a field to "field structures" $\mathbf{K} = (K, X, K_x)_{x \in X}$ with $\mathrm{Gal}(K_x) = G_x$ for each $x \in X$. We give here an arithemtically geometric criterion for $\mathbf{G}$ to be projective. It generalizes Ax's theorem saying that $\mathrm{Gal}(K)$ is projective if $K$ is PAC. The standard proof of Ax's result [FrJ, p. 137] actually uses only the existence of $K$-rational points on varieties over $K$ which become unirational over a finite extension of $K$. Our criterion has the same nature amended with a local-golbal flavor.

Let $K$ be a field. Denote the set of all algebraic (resp. separable algebraic) extensions of $K$ by $\mathrm{AlgExt}(K)$ (resp. $\mathrm{SepAlgExt}(K)$). Galois theory puts $\mathrm{SepAlgExt}(K)$ in a bijective order-reversing correspondence with $\mathrm{Subgr}(\mathrm{Gal}(K))$. It equips $\mathrm{SepAlgExt}(K)$ with two natural topologies, the **strict topology** and the **étale topology**. A basic étale open subset of $\mathrm{SepAlgExt}(K)$ is $\mathrm{SepAlgExt}(L)$, where $L$ is a finite extension of $K$. Thus, $\mathrm{SepAlgExt}(K)$ is not étale Hausdorff unless $K = K_s$. A basic strictly open subset of $\mathrm{SepAlgExt}(K)$ is $\{K' \in \mathrm{SepAlgExt}(K) \mid L \cap K' = L_0\}$ where $L_0$ is a finite separable extension of $K$ and $L$ is a finite Galois extension of $K$ containing $L_0$. $\mathrm{SepAlgExt}(K)$ is a profinite space under the strict topology. Denote the strict closure of a subset $\mathcal{X}$ of $\mathrm{SepAlgExt}(K)$ by $\mathrm{StrictClosure}(\mathcal{X})$.

A **field structure** is a triple $\mathbf{K} = (K, X, \delta)$ consisting of a field $K$, a profinite space $X$, an étale continuous map $\delta \colon X \to \mathrm{SepAlgExt}(K)$, and an étale continuous action (from the right) of $\mathrm{Gal}(K)$ on $X$ satisfying the following condition:

(1a) For each $x \in X$ put $K_x = \delta(x)$. Then $K_{x^\sigma} = K_x^\sigma$ for all $x \in X$ and $\sigma \in \mathrm{Gal}(K)$.

(1b) $x \in X$, $\sigma \in \mathrm{Gal}(K)$, and $x^\sigma = x$ imply $\sigma \in \mathrm{Gal}(K_x)$.

As with group structures, we usually write $\mathbf{K}$ as $(K, X, K_x)_{x \in X}$. The **absolute Galois group structure** associated with $\mathbf{K}$ is $\mathrm{Gal}(\mathbf{K}) = (\mathrm{Gal}(K), X, \mathrm{Gal}(K_x))_{x \in X}$. Conversely, to each absolute Galois group structure $\mathbf{G} = (\mathrm{Gal}(K), X, G_x)_{x \in X}$ we associate a field structure $\mathbf{K} = (K, X, K_x)_{x \in X}$, where $K_x$ is the fixed field of $G_x$ in $K_s$. Then $\mathrm{Gal}(\mathbf{K}) = \mathbf{G}$. We use the correspondence between field structures and absolute Galois group structures to translate the terminology and results obtained so far from group structures to field structures.

*Definition 6.1:* A **unirational arithmetical problem** for a field structure $\mathbf{K} = (K, X, K_x)_{x \in X}$ is a data

$$(2) \qquad \Phi = (V, X_i, L_i, \pi_i \colon U_i \to V \times_K L_i)_{i \in I_0}$$

satisfying these conditions:

(3a) $(\mathrm{Gal}(L_i), X_i)_{i \in I_0}$ is a special partition of $\mathrm{Gal}(\mathbf{K})$ (Definition 3.5).

(3b) $V$ is a smooth affine variety over $K$.

(3c) $U_i$ is a smooth variety over $L_i$ birationally equivalent to $\mathbb{A}_{L_i}^{\dim(V)}$.

(3d) $\pi_i \colon U_i \to V \times_K L_i$ is an étale morphism.

Let $X' = \bigcup_{i \in I_0} X_i$. A **solution** of $\Phi$ is an "extended point" $(\mathbf{a}, \mathbf{b}_x)_{x \in X'}$ with $\mathbf{a} \in V(K)$, $\mathbf{b}_x \in U_i(K_x)$, and $\pi_i(\mathbf{b}_x) = \mathbf{a}$ for each $i \in I_0$ and all $x \in X_i$. Call $\mathbf{K}$ **unirationally closed** if each unirational arithmetical problem for $\mathbf{K}$ has a solution. ∎

LEMMA 6.2 ([HJK, Lemma 3.1]): *Let $L/K$ be a finite Galois extension. Let $\psi \colon B \to \mathrm{Gal}(L/K)$ be an epimorphism of finite groups. Then there exists a finitely generated regular extension $E$ of $K$ and a finite Galois extension $F$ of $E$ containing $L$ such that $B = \mathrm{Gal}(F/E)$ and $\psi$ is the restriction $\mathrm{res}_{F/L} \colon \mathrm{Gal}(F/E) \to \mathrm{Gal}(L/K)$.*

*Moreover, let $K \subseteq L_0 \subseteq L$ and $E \subseteq F_0 \subseteq F$ be fields with $L_0 \subseteq F_0$. Suppose $\psi \colon \mathrm{Gal}(F/F_0) \to \mathrm{Gal}(L/L_0)$ is an isomorphism. Then $F_0$ is a purely transcendental extension of $L_0$ of transcendence degree $|B|$.*

*Proof:* Let $x^\beta$, $\beta \in B$, be algebraically independent elements over $K$. Define an action of $B$ on $F = L(x^\beta \mid \beta \in B)$ by $(x^\beta)^{\beta'} = x^{\beta\beta'}$ and $a^{\beta'} = a^{\psi(\beta')}$ for $a \in L$. Denote the fixed field of $B$ in $F$ by $E$. Then $F/K$ is a finitely generated separable extension. By [Lan, p. 64, Prop. 6], $E/K$ is also a finitely generated separable extension. Also, res: $\mathrm{Gal}(F/E) \to \mathrm{Gal}(L/K)$ coincides with $\psi \colon B \to \mathrm{Gal}(L/K)$. So, $E \cap \tilde{K} = E \cap F \cap \tilde{K} = E \cap L = K$. Therefore, $E/K$ is regular.

Now let $L_0$ and $F_0$ be as in the second paragraph of the lemma. Put $B_0 = \mathrm{Gal}(F/F_0)$. Choose a set of representatives $R$ for the left cosets of $B$ modulo $B_0$. Let $w_1, \ldots, w_m$ be a basis for $L/L_0$. By assumption, $m = |B_0|$. Consider $\rho \in R$. Put

$$t_{\rho j} = \sum_{\beta \in B_0} w_j^\beta x^{\rho\beta}, \qquad j = 1, \ldots, m.$$

Since $\det(w_j^\beta) \neq 0$, each $x^{\rho\beta}$ is a linear combination of $t_{\rho j}$ with coefficients in $L$. Put $\mathbf{t} = (t_{\rho j} \mid \rho \in R, \; j = 1, \ldots, m)$, $\mathbf{x} = (x^\beta \mid \beta \in B)$, and $n = |B|$. Both tuples contain exactly $n$ elements and $L(\mathbf{t}) = L(\mathbf{x}) = F$. So, $L_0(\mathbf{t})$ is a purely transcendental extension of $L_0$.

Each $t_{\rho j}$ is fixed by $H$. So, $L_0(\mathbf{t}) \subseteq F_0$. Moreover, $m = [L : L_0] = [L(\mathbf{t}) : L_0(\mathbf{t})] \geq [F : F_0] = |B_0| = m$. Conclude that $F_0 = L_0(\mathbf{t})$ and $F_0/L_0$ is purely transcendental. ∎

LEMMA 6.3: *Let* $\mathbf{G} = (\mathbf{G}, X, G_x)_{x \in X}$ *and* $\mathbf{A} = (\mathbf{A}, I, A_i)_{i \in I}$ *be groups structures and* $\varphi \colon \mathbf{G} \to \mathbf{A}$ *be an epimorphism. Suppose* $S_x = G_x$ *for each* $x \in X$. *Then* $S_i = A_i$ *for each* $i \in I$.

*Proof:* By assumption, $S_i \leq A_i$. Conversely, let $a \in A_i$. By assumption, there is $x \in X$ with $\varphi(x) = i$ and $\varphi(G_x) = A_i$. Choose $g \in G_x$ with $\varphi(g) = a$. Then $i^a = \varphi(x)^{\varphi(g)} = \varphi(x^g) = \varphi(x) = i$. Thus, $a \in S_i$. ∎

PROPOSITION 6.4: *Let* $\mathbf{K} = (K, X, K_x)_{x \in X}$ *be a unirationally closed field structure. Suppose* $S_x = \mathrm{Gal}(K_x)$ *for each* $x \in X$. *Then* $\mathrm{Gal}(\mathbf{K})$ *is a projective group structure.*

*Proof:* Let $(\varphi \colon \mathrm{Gal}(\mathbf{K}) \to \mathbf{A}, \; \alpha \colon \mathbf{B} \to \mathbf{A})$ be a finite embedding problem for $\mathrm{Gal}(\mathbf{K})$. Thus $\mathbf{B} = (B, J, B_j)_{j \in J}$ and $\mathbf{A} = (A, I, A_i)_{i \in I}$ are finite group structures and $\alpha$ is a cover. By Lemma 4.1, we may assume $\varphi$ is an epimorphism. By assumption, $I$ is discrete and the map $\varphi \colon X \to I$ is continuous. Hence, $X_i = \{x \in X \mid \varphi(x) = i\}$ is an open-closed subset of $X$, $i \in I$ and $X = \bigcup_{i \in I} X_i$. Moreover, $X_i^\sigma = X_{i^{\varphi(\sigma)}}$ for all $i \in I$ and $\sigma \in \mathrm{Gal}(K)$.

Choose a set of representatives $I_0$ for the $A$-orbits of $I$ and for each $i \in I_0$ choose $j(i) \in J$ such that $\alpha(j(i)) = i$.

The rest of the proof has six parts.

PART A: *Replacing* **A** *and* **B** *by Galois structures.* Replace $A$ by $\mathrm{Gal}(L/K)$, where $L$ is a finite Galois extension of $K$, to assume that $\varphi\colon \mathrm{Gal}(K) \to \mathrm{Gal}(L/K)$ is $\mathrm{res}_{K_s/L}$. Denote the fixed field of $A_i$ in $L$ by $L_i$. By assumption, $S_x = \mathrm{Gal}(K_x)$ for each $x \in X$. Hence, by Lemma 6.3, $\mathrm{Gal}(L/L_i) = \{\sigma \in \mathrm{Gal}(L/K) \mid i^\sigma = \sigma\}$. Therefore,

$$(4) \qquad\qquad \mathrm{Gal}(L_i) = \{\sigma \in \mathrm{Gal}(K) \mid X_i^\sigma = X_i\}.$$

and $(\mathrm{Gal}(L_i), X_i, \mathrm{Gal}(K_x))_{x \in X_i}$ is a group structure.

Lemma 6.2 gives a finitely generated regular extension $E$ of $K$ and a finite Galois extension $F$ of $E$ containing $L$ and allows us to replace $B$ by $\mathrm{Gal}(F/E)$ and $\alpha\colon B \to \mathrm{Gal}(L/K)$ by $\mathrm{res}_{F/L}\colon \mathrm{Gal}(F/E) \to \mathrm{Gal}(L/K)$. For each $i \in I_0$ denote the fixed field of $B_{j(i)}$ in $F$ by $F_i$. Since $\alpha$ is a cover, $\mathrm{res}_{F/L}\colon \mathrm{Gal}(F/F_i) \to \mathrm{Gal}(L/L_i)$ is an isomorphism. So, by Lemma 6.2, $F_i$ is a purely transcendental extension of $L_i$ of transcendence degree $r = [F:E]$.

Since $\varphi\colon \mathrm{Gal}(\mathbf{K}) \to \mathbf{A}$ is a morphism, $\mathrm{res}_{K_s/L}\big(\mathrm{Gal}(K_x)\big) \leq \mathrm{Gal}(L/L_i)$ for all $x \in X_i$. Hence, $L_i \leq K_x$.

$$
\begin{array}{ccccc}
E & \!\!\!\text{---}\!\!\! & F_i & \!\!\!\text{---}\!\!\! & F \\
| & & | & & | \\
K & \!\!\!\text{---}\!\!\! & L_i & \!\!\!\text{---}\!\!\! & L
\end{array}
$$

PART B: *Setting up a unirational arithmetical problem.* Choose $y_1, \ldots, y_n \in E$, $z_i \in F_i$, and $\tilde{z} \in F$ satisfying this:

(5a) $E = K(\mathbf{y})$ and $V = \mathrm{Spec}(K[\mathbf{y}])$ is a smooth affine absolutely irreducible subvariety of $\mathbb{A}_K^n$ with generic point $\mathbf{y}$.

(5b) For each $i \in I_0$ the following holds: $F_i = L_i(\mathbf{y}, z_i)$ and $U_i = \mathrm{Spec}(L_i[\mathbf{y}, z_i])$ is a smooth Zariski closed subvariety of $\mathbb{A}_{L_i}^{n+1}$ birationally equivalent to $\mathbb{A}_{L_i}^r$ with generic point $(\mathbf{y}, z_i)$.

(5c) $z_i$ is integral over $L_i[\mathbf{y}]$ and the discriminant of $\mathrm{irr}(z_i, L_i(\mathbf{y}))$ is a unit of $L_i[\mathbf{y}]$. Hence, $L_i[\mathbf{y}, z_i]/L_i[\mathbf{y}]$ is a ring cover in the terminology of [FrJ, Definition 5.4]. Thus, projection on the first $n$ coordinates is an étale morphism $\pi_i\colon U_i \to V \times_K L_i$.

(5d) $F = K(\mathbf{y}, \tilde{z})$ and $L[\mathbf{y}, \tilde{z}]/L[\mathbf{y}]$ is a ring cover.

Then, (2) is a unirational arithmetical problem for $\mathbf{K}$ satisfying Condition (3).

PART C: *A solution of a unirational arithmetical problem.* Since **K** is unirationally closed, Problem (2) has a solution. Thus, there are $\mathbf{a} \in V(K)$ and $\mathbf{b}_x = (\mathbf{a}, c_x) \in U_i(K_x)$ for each $i \in I_0$ and all $x \in X_i$.

Let $i \in I_0$ and $x \in X_i$. Then $\mathrm{Gal}(L_i(c_x))$ is an open subgroup of $\mathrm{Gal}(L_i)$ which contains $\mathrm{Gal}(K_x)$. Also, $W_x = \{x' \in X_i \mid L_i(c_x) \subseteq K_{x'}\}$ is an étale open subset of $X_i$ which contains $x^{\mathrm{Gal}(L_i(c_x))}$. Lemma 3.6 with $\mathrm{Gal}(L_i), X_i, X_i, x, \mathrm{Gal}(L_i(c_x)), W_x$ respectively replacing $G, X, Y, y, G'_y, V_y$, gives

(6a) a finite set $\Lambda_i$, and

(6b) for each $l \in \Lambda_i$ an open-closed subset $X_{il}$ of $X_i$, an element $x_{il} \in X_{il}$, and a finite subset $T_{il}$ of $\mathrm{Gal}(L_i)$,

satisfying the following conditions with $c_{il} = c_{x_{il}}$:

(7a) $\mathrm{Gal}(L_i(c_{il})) = \{\sigma \in \mathrm{Gal}(K) \mid X_{il}^\sigma = X_{il}\}$ for each $l \in \Lambda_i$.

(7b) $\mathrm{Gal}(L_i) = \bigcup_{\tau \in T_{il}} \mathrm{Gal}(L_i(c_{il}))\tau$.

(7c) $X_i = \bigcup_{l \in \Lambda_i} \bigcup_{\tau \in T_{il}} X_{il}^\tau$.

Thus, $(\mathrm{Gal}(L_i(c_{il})), T_{il}, X_{il})_{l \in \Lambda_i}$ is a special partition of $(\mathrm{Gal}(L_i), X_i, \mathrm{Gal}(K_x))_{x \in X_i}$.

PART D: *A homomorphism* $\gamma\colon G \to B$. Since $\mathbf{a}$ is simple on $V$, there is a $K$-place $\rho\colon E \to K \cup \{\infty\}$ with $\rho(\mathbf{y}) = \mathbf{a}$ [JaR, Cor. A2]. Extend $\rho$ to an $L$-place $\rho\colon F \to \tilde{K} \cup \{\infty\}$. Let $\bar{F}$ be the residue field of $\rho$. By (5d), $\bar{F}$ is a finite Galois extension of $K$ containing $L$ [FrJ, Lemma 5.5]. Moreover, there is an embedding $\rho^*\colon \mathrm{Gal}(\bar{F}/K) \to \mathrm{Gal}(F/E)$ with $\rho(\rho^*(\sigma)u) = \sigma(\rho(u))$ for all $\sigma \in \mathrm{Gal}(\bar{F}/K)$ and $u \in F$ with $\varphi(u) \neq \infty$ [FrJ, Lemma 5.5]. Then $\gamma = \rho^* \circ \mathrm{res}_{K_s/\bar{F}}$ is a homomorphism from $\mathrm{Gal}(K)$ to $\mathrm{Gal}(F/E)$ with $\mathrm{res}_{F/L} \circ \gamma = \mathrm{res}_{K_s/L}$.

PART E: *A constant map* $\gamma\colon X_{il} \to J$. Let $i \in I_0$ and $l \in \Lambda_i$. By (5b), $(\mathbf{a}, c_{il})$ is simple on $U_i$. So, there is an $L_i$-place $\rho_{il}\colon F_i \to L_i(c_{il}) \cup \{\infty\}$ with $\rho_{il}(\mathbf{y}, z_i) = (\mathbf{a}, c_{il})$ [JaR, Cor. A2]. Extend it to an $L$-place $\rho_{il}\colon F \to \tilde{K} \cup \{\infty\}$. Since $\rho_{il}|_{EL} = \rho|_{EL}$, there is $\sigma_{il} \in \mathrm{Gal}(F/EL)$ and $\rho_{il} = \rho \circ \sigma_{il}^{-1}$. Define a continuous map $\gamma\colon X_{il} \to J$ by $\gamma(x) = j(i)^{\sigma_{il}}$ for all $x \in X_{il}$. Then $\rho(F_i^{\sigma_{il}}) = \rho \circ \sigma_{il}^{-1}(F_i) = \rho_{il}(F_i) \subseteq L_i(c_{il}) \cup \{\infty\} \subseteq K_x \cup \{\infty\}$. This implies

(8) $\gamma(\mathrm{Gal}(K_x)) \leq \gamma(\mathrm{Gal}(L_i(c_{il}))) \leq \mathrm{Gal}(F/F_i)^{\sigma_{il}} = B_{j(i)^{\sigma_{il}}} = B_{\gamma(x)}, \qquad x \in X_{il}$.

Also, $\alpha(\gamma(x)) = \alpha(j(i)^{\sigma_{il}}) = \alpha(j(i)) = i = \varphi(x)$.

Let $X' = \bigcup_{i \in I_0} \bigcup_{l \in \Lambda_i} X_{il}$. Then $X = \bigcup_{\sigma \in \mathrm{Gal}(K)} (X')^\sigma$. Extend $\gamma \colon X' \to J$ to $X$ by $\gamma(x^\sigma) = \gamma(x)^{\gamma(\sigma)}$ for each $\sigma \in \mathrm{Gal}(K)$. To prove this is a good definition, we have to show that if $x, y \in X'$, $\sigma_1, \sigma_2 \in \mathrm{Gal}(K)$ and $x^{\sigma_1} = y^{\sigma_2}$, then $\gamma(x)^{\gamma(\sigma_1)} = \gamma(y)^{\gamma(\sigma_2)}$. In other words, with $\sigma = \sigma_1 \sigma_2^{-1}$, we have to prove that

(9) $x \in X'$, $\sigma \in \mathrm{Gal}(K)$, and $x^\sigma \in X'$ imply $\gamma(x)^{\gamma(\sigma)} = \gamma(x^\sigma)$.

Indeed, there are $i, i' \in I_0$, $l \in \Lambda_i$, and $l' \in \Lambda_{i'}$ with $x \in X_{il}$ and $x^\sigma \in X_{i'l'}$. By (6b), $x \in X_i$ and $x^\sigma \in X_{i'}$. Hence, $x^\sigma \in X_i^\sigma \cap X_{i'} = X_{i^\sigma} \cap X_{i'}$. Therefore, $i^\sigma = i'$. By the choice of $I_0$, this implies $i = i'$. Hence $x^\sigma \in X_{i^\sigma} \cap X_i$, so $X_i^\sigma = X_i$. By (4), $\sigma \in \mathrm{Gal}(L_i)$. It follows, $x \in X_{il}$ and $x^\sigma \in X_{il'}$, so $x^\sigma \in X_{il}^\sigma \cap X_{il'}$. By (7b), $\sigma = \sigma' \tau$ with $\sigma' \in \mathrm{Gal}(L_i(c_{il}))$ and $\tau \in T_{il}$. Then, by (7c) and (7a), $X_{il'} = X_{il}^{\sigma' \tau} = X_{il}^\tau$. Hence, by (7c), $\tau \in \mathrm{Gal}(L_i(c_{il}))$, so also $\sigma \in \mathrm{Gal}(L_i(c_{il}))$. By (7c), $X_{il'} = X_{il}$. We have therefore proved both $x$ and $x^\sigma$ belong to $X_{il}$. By definition, $\gamma(x) = j(i)^{\sigma_{il}} = \gamma(x^\sigma)$. By (8), $\gamma(\sigma) \in \gamma(\mathrm{Gal}(L_i(c_{il}))) \le B_{\gamma(x)}$. Hence, $\varphi(\sigma) = \alpha(\gamma(\sigma)) \le A_{\varphi(x)}$. By assumption, $S_x = \mathrm{Gal}(K_x)$ for each $x \in X$. By Lemma 6.3, $\alpha(\gamma(\sigma)) \in S_{\varphi(x)}$. Hence, by Lemma 2.2, $\gamma(x)^{\gamma(\sigma)} = \gamma(x)$. Therefore $\gamma(x^\sigma) = \gamma(x) = \gamma(x)^{\gamma(\sigma)}$, as claimed.

PART F: *Conclusion of the proof.* By (9), $\gamma(x)^{\gamma(\sigma)} = \gamma(x^\sigma)$ for all $x \in X$ and $\sigma \in \mathrm{Gal}(K)$. Hence, by (8), $\gamma(\mathrm{Gal}(K_{x^\sigma})) \le B_{\gamma(x^\sigma)}$ for all $x \in X$. Thereofre, $\gamma \colon \mathrm{Gal}(\mathbf{K}) \to \mathbf{B}$ is a morphism. Finally, $\alpha \circ \gamma = \varphi$ on $\mathrm{Gal}(K)$ and on $X'$, hence on $X$. Thus, $\gamma$ solves the embedding problem we posed for $\mathrm{Gal}(\mathbf{K})$. It follows, $\mathrm{Gal}(\mathbf{K})$ is projective. ∎

## 7. Valued Fields

The results of this section are well known, although there is some novelty in the presentation*. We begin with a brief review of inertia and ramification groups.

Denote the residue field of a valued field $(F, v)$ by $\bar{F}$. For each $x \in F$ with $v(x) \geq 0$ let $\bar{x}$ be the residue of $x$ in $\bar{F}$. Finally, let $F_{\text{ins}}$ be the maximal purely inseparable extension of $F$.

Consider a Galois extension $(N, v)/(F, v)$ of Henselian fields. Then $\bar{N}/\bar{F}$ is a normal extension. For each $\sigma \in \text{Gal}(N/F)$ define $\bar{\sigma} \in \text{Aut}(\bar{N}/\bar{F})$ by this rule: $\bar{\sigma}\bar{x} = \overline{\sigma x}$ for $x \in N$ with $v(x) \geq 0$. The map $\sigma \mapsto \bar{\sigma}$ is an epimorphism $\rho \colon \text{Gal}(N/F) \to \text{Aut}(\bar{N}/\bar{F})$ [End, Thm. 19.6]. Its kernel is the **inertia group**:

$$G_0(N/F) = \{\sigma \in \text{Gal}(N/F) \mid v(\sigma x - x) > 0 \text{ for each } x \in N \text{ with } v(x) \geq 0\}.$$

Denote the fixed field in $N$ of $G_0(N/F)$ by $N_0$. Then $\bar{N}_0$ is the maximal separable extension of $\bar{F}$ in $\bar{N}$ [End, Thm. 19.12]. Hence, $\bar{N}_0/\bar{F}$ is Galois and there is a short exact sequence

$$(1) \qquad 1 \longrightarrow \text{Gal}(N/N_0) \longrightarrow \text{Gal}(N/F) \overset{\rho}{\longrightarrow} \text{Gal}(\bar{N}_0/\bar{F}) \longrightarrow 1.$$

Here we have identified each $\bar{\sigma} \in \text{Aut}(\bar{N}/\bar{F})$ with its restriction to $\bar{N}_0$. In addition, $v(N_0^\times) = v(F^\times)$ [End, Cor. 19.14]. Hence, $N_0/F$ is an unramified extension.

The **ramification group** of $\text{Gal}(N/F)$ is

$$G_1(N/F) = \Big\{\sigma \in \text{Gal}(N/F) \mid v\Big(\frac{\sigma x}{x} - 1\Big) > 0 \text{ for each } x \in N\Big\}.$$

It is a normal subgroup of $\text{Gal}(N/F)$ which is contained in $G_0(N/F)$ [End, (20.8)]. Denote the fixed field of $G_1(N/F)$ in $N$ by $N_1$. When $p = \text{char}(\bar{F}) > 0$, $\text{Gal}(N/N_1)$ is the unique $p$-Sylow subgroup of $\text{Gal}(N/N_0)$ [End, Thm. 20.18]. When $\text{char}(\bar{F}) = 0$, $\text{Gal}(N/N_1)$ is trivial. So, in both cases, $\text{char}(\bar{F})$ does not divide $[N_1 : N_0]$.

Suppose now $N = F_s$. Then $N_0 = F_u$ is the **inertia field** and $N_1 = F_r$ is the **ramification field** of $F$. In this case (1) becomes the short exact sequence

$$(2) \qquad 1 \longrightarrow \text{Gal}(F_u) \longrightarrow \text{Gal}(F) \overset{\rho}{\longrightarrow} \text{Gal}(\bar{F}) \longrightarrow 1.$$

---

* This section is a rewrite of [HJK, Sec. 2].

Also, $F \subseteq F_u \subseteq F_r \subseteq F_s$, $F_u/F$ and $F_r/F$ are Galois extensions, $\mathrm{char}(\bar{F}) \nmid [F_r : F_u]$, and $\mathrm{Gal}(F_r)$ is a pro-$p$ group if $p = \mathrm{char}(\bar{F}) \neq 0$.

Consider now a finite extension $(L, v)/(F, v)$ of Henselian fields. Let $e = e(L/F) = (v(L^\times) : v(F^\times))$ be the **ramification index**. There is a positive integer $d$ such that $[L : F] = de[\bar{L} : \bar{F}]$. If $\mathrm{char}(\bar{F}) = p > 0$, then $d$ is a power of $p$ [Art, p. 62, Thm. 10]. If $\mathrm{char}(\bar{F}) = 0$, then $d = 1$. When $d = 1$ we say $L/F$ is **defectless**. An arbitrary algebraic extension $M/F$ is **defectless** if each finite subextension is defectless. This is the case when $\mathrm{char}(\bar{F}) \nmid [M : K]$. For example, $F_r/F_u$ is defectless. In addition, by (2), $[L : F] = [\bar{L} : \bar{F}]$ for each finite subextension $L/F$ of $F_u/F$. Hence, $F_u/F$ is defectless. Consequently, $F_r/F$ is defectless.

**Lemma 7.1:** *Let $(F, v)$ be a Henselian valued field. Use the above notation.*

(a) *There is a field $F'$ with $F_u F' = F_r$ and $F_u \cap F' = F$.*

(b) *The short exact sequence $1 \to \mathrm{Gal}(F_r/F_u) \to \mathrm{Gal}(F_r/F) \to \mathrm{Gal}(F_u/F) \to 1$ splits.*

*Proof:* Statement (b) is a Galois theoretic interpretation of (a). So, we prove (a).

Zorn's lemma gives a maximal extension $F'$ of $F$ in $F_r$ with residue field $\bar{F}$. For each prime number $l \neq \mathrm{char}(\bar{F})$ the value group of $F'$ is $l$-divisible. Otherwise, there is $a \in F'$ with $v(a) \notin lv((F')^\times)$. Put $L = F(\sqrt[l]{a})$. Then $[L : F'] = l$ and $l \leq (v(L^\times) : v((F')^\times))$. Since $e(L/F')[\bar{L} : \overline{F'}] \leq [L : F'] = l$, we have $\bar{L} = \overline{F'} = \bar{F}$. Recall: $\mathrm{Gal}(F_r)$ is a pro-$p$ group if $\mathrm{char}(\bar{F}) = p > 0$ and trivial if $\mathrm{char}(\bar{F}) = 0$. So, $L \subseteq F_r$. This contradicts the maximality of $F'$.

By the discussion preceding Lemma 7.1, $F_u \cap F' = F$. Let $E = F_u F'$. Consider a prime number $l \neq \mathrm{char}(\bar{F})$. Since $E/F'$ is an algebraic extension, $v(E^\times)$ is contained in the divisible hull of $v(F')$. As $v((F')^\times)$ is $l$-divisible, so is $v(E^\times)$. Since $F_u \subseteq E \subseteq F_r \subseteq F_s$ and $\bar{F}_u = \bar{F}_s$, we have $\bar{E} = \bar{F}_r$. Hence, $e(E'/E) = [\overline{E'} : \bar{E}] = 1$ and therefore $[E' : E] = 1$ (because $F_r/E'$ is defectless) for every finite extension $E'$ of $E$ in $F_r$. Consequently, $E = F_r$. ∎

**Lemma 7.2** (Kuhlmann-Pank-Roquette [KPR, Thm. 2.2]): *Let $(F, v)$ be a Henselian field.*

(a) *There is a field $F'$ with $F_r \cap F' = F$ and $F_r F' = F_s$.*

(b) *The short sequence $1 \to \mathrm{Gal}(F_r) \to \mathrm{Gal}(F) \to \mathrm{Gal}(F_r/F) \to 1$ splits.*

*Proof:* Statement (a) is a Galois theoretic interpretation of (b). So, we prove (b). Let $p = \mathrm{char}(\bar{F})$. If $p = 0$, then $F_r = F_s$ and we may take $F' = F$. Suppose $p \neq 0$.

By (2), $\mathrm{Gal}(F_u/F) \cong \mathrm{Gal}(\bar{F})$. By Witt, the $p$-Sylow subgroups of $\mathrm{Gal}(\bar{F})$ are free [Rib, p. 256, Thm. 3.3]. Hence, so are the $p$-Sylow subgroups of $\mathrm{Gal}(F_u/F)$. Since $p \nmid [F_r : F_u]$, restriction maps each $p$-Sylow subgroup of $\mathrm{Gal}(F_r/F)$ isomorphically onto a $p$-Sylow subgroup of $\mathrm{Gal}(F_u/F)$. Hence, each $p$-Sylow subgroup of $\mathrm{Gal}(F_r/F)$ is free. Thus $\mathrm{cd}_p(\mathrm{Gal}(F_r/F)) = 1$ [Rib, p. 207, Cor. 2.2]. As $\mathrm{Gal}(F_r)$ is a pro-$p$ group, the short sequence in (b) splits [Rib, p. 211, Prop. 3.1(iii)']. ∎

**PROPOSITION 7.3:** *Let $(F, v)$ be a valued field.*

(a) *Suppose $(F, v)$ is Henselian. Then the epimorphism $\rho \colon \mathrm{Gal}(F) \to \mathrm{Gal}(\bar{F})$ induced by reduction at $v$ splits.*

(b) *Each subgroup of $\mathrm{Gal}(\bar{F})$ is isomorphic to a subgroup of $\mathrm{Gal}(F)$.*

*Proof of (a):* The map $\rho$ decomposes as $\mathrm{Gal}(F) \xrightarrow{\mathrm{res}} \mathrm{Gal}(F_r/F) \xrightarrow{\mathrm{res}} \mathrm{Gal}(F_u/F) \xrightarrow{\bar{\rho}} \mathrm{Gal}(\bar{F})$. The map $\bar{\rho}$ which is also induced by reduction is an isomorphism (by (2)). By Lemmas 7.1 and 7.2, each of the restriction maps splits. Hence $\rho$ splits.

*Proof of (b):* Let $(F', v)$ be the Henselization of $(F, v)$. Then $\overline{F'} = \bar{F}$. By (a), each subgroup of $\mathrm{Gal}(\bar{F})$ is isomorphic to a subgroup of $\mathrm{Gal}(F')$, hence of $\mathrm{Gal}(F)$. ∎

**PROPOSITION 7.4:** *Let $F/K$ be an extension of fields. Suppose $v$ is a valuation of $F$ which is trivial on $K$ and $\bar{F} = K$. Then*

(a) *res: $\mathrm{Gal}(F) \to \mathrm{Gal}(K)$ is an epimorphism which splits. If, in addition, $(F, v)$ is Henselian, then res is the epimorphism induced by reduction at $v$.*

(b) *$(F, v)$ has a separable algebraic Henselian extension $(F', v)$ such that res: $\mathrm{Gal}(F') \to \mathrm{Gal}(K)$ is an isomorphism and $\overline{F'}$ is a purely inseparable extension of $K$.*

(c) *Suppose $K$ is perfect. Then $(F, v)$ has an algebraic Henselian extension $(F'', v)$ such that $F''$ is perfect, $\overline{F''} = K$, and res: $\mathrm{Gal}(F'') \to \mathrm{Gal}(K)$ is an isomorphism.*

*Proof:* Replace $(F, v)$ by a Henselian closure, if necessary, to assume $(F, v)$ is Henselian. Then assume without loss that $\bar{a} = a$ for each $a \in K_s$. Let $\rho \colon \mathrm{Gal}(F) \to \mathrm{Gal}(K)$ be the epimorphism induced by reduction at $v$. For each $a \in K$ and each $\sigma \in \mathrm{Gal}(K)$ we have, $\sigma a = \overline{\sigma a} = \bar{\sigma} \bar{a} = \bar{\sigma} a = \rho(\sigma) a$. Thus, res: $\mathrm{Gal}(F) \to \mathrm{Gal}(K)$ coincides with $\rho$.

Proposition 7.3(a) gives a section $\rho' \colon \mathrm{Gal}(K) \to \mathrm{Gal}(F)$ of $\rho$. Let $F'$ be the fixed field of $\rho'(\mathrm{Gal}(K))$ in $F_s$. Then $\mathrm{Gal}(F') \to \mathrm{Gal}(K)$ is an isomorphism. Also, for all $u \in F'$ and $\sigma \in \mathrm{Gal}(F')$ we have $\bar{\sigma} \bar{u} = \overline{\sigma u} = \bar{u}$. Hence, $\overline{F'}$ is a purely inseparable extension of $K$. This concludes the proof of (a) and (b).

When $K$ is perfect, $F'' = F'_{\mathrm{ins}}$ satisfies (c). ∎

The following Proposition gives more details to a result of Efrat [Efr. Prop. 4.7].

PROPOSITION 7.5: *Let $K$ be a field, $E_0$ its prime field, and $T$ a set of variables with* card$(T) \geq$ trans.deg$(K/E_0)$. *Let $F_0$ be either $E_0$ or $\mathbb{Q}$. Then there is a field $L$, algebraic over $F_0(T)$, with $G(L) \cong G(K)$.*

*Proof:* There is a unique place $\varphi_0 \colon F_0 \to E_0 \cup \{\infty\}$. Choose a transcendence base $\bar{T}$ for $K/E_0$. By assumption, card$(\bar{T}) \leq$ card$(T)$. Choose a surjective map $\varphi_1 \colon T \to \bar{T}$. Let $F_1 = F_0(T)$ and $E_1 = E_0(\bar{T})$. Extend $\varphi_0$ and $\varphi_1$ to a place $\varphi \colon F_1 \to E_1 \cup \{\infty\}$. Denote the corresponding valuation by $v$. Corollary 7.3(b) gives the desired field $L$. ∎

*Definition 7.6: Rigid Henselian extensions.* Let $K$ be a field and $(L,v)$ a valued field. We say, $(L,v)$ is a **rigid Henselian extension** of $K$ if $(L,v)$ is Henselian, $K \subseteq L$, $v$ is trivial on $K$, $\bar{L}_v = K$, and res: $\mathrm{Gal}(L) \to \mathrm{Gal}(K)$ is an isomorphism. In this case we also call the place $\varphi \colon L \to K \cup \{\infty\}$ associated with $v$ **rigid**.

An arbitrary field extension $L/K$ is a **rigid Henselian extension** if $L$ admits a valuation $v$ such that $(L,v)$ is a rigid Henselian extension $K$. ∎

PROPOSITION 7.7: *Let $F/K$ be a purely transcendental extension. Then:*

(a) *$F$ has a valuation $v$ which is trivial on $K$ and $\bar{F} = K$.*

(b) *$F$ has a separable algebraic extension $F'$ such that res: $\mathrm{Gal}(F') \to \mathrm{Gal}(K)$ is an isomorphism.*

(c) *If $K$ is perfect, then $F$ has a perfect algebraic extension $F'$ which is a rigid Henselian extension of $K$.*

*Proof of (a):* The assertion is evident when $F = K(t)$ and $t$ is transcendental. The general case follows from the special case by transfinite induction and using composition of valuations.

*Proof of (b):*  Apply Proposition 7.4(b).

*Proof of (c):*  Apply Proposition 7.4(c).  ∎

## 8. The Space of Valuations of a Field

Let $K$ be a field. Denote the collection of all valuations of $K$ by $\mathrm{Val}(K)$. We include in $\mathrm{Val}(K)$ also the trivial valuation $v_0$ defined by $v_0(a) = 0$ for each $a \in K^\times$ and $v_0(0) = \infty$. Also, we do not distinguish between equivalent valuations. Thus, we identify valuations with the same valuation rings. Given $a \in K$, we write

$$\mathrm{Val}_a(K) = \{v \in \mathrm{Val}(K) \mid v(a) > 0\}, \qquad \mathrm{Val}'_a(K) = \{v \in \mathrm{Val}(K) \mid v(a) \geq 0\}.$$

Intersections of finitely many sets of these form build a basis for a topology on $\mathrm{Val}(K)$, the so called **patch topology** (see more about the patch topology in [Hoe, Sec. 2]).

The following identities make the use of open subsets of $\mathrm{Val}(K)$ easier:

(1)
$$\mathrm{Val}_{a/b}(K) = \{v \in K \mid v(a) > v(b)\}. \qquad \mathrm{Val}'_{a/b}(K) = \{v \in K \mid v(a) \geq v(b)\}$$
$$\mathrm{Val}'_a(K) = \mathrm{Val}(K) \smallsetminus \mathrm{Val}_{a^{-1}}(K), \quad \mathrm{Val}_0(K) = \mathrm{Val}(K), \quad \mathrm{Val}_1(K) = \emptyset$$

*Example 8.1:* $\mathrm{Val}(\mathbb{Q})$. It consists of $v_p$, with $p$ ranging over all prime numbers, and $v_0$. For each $p$, $\{v \in \mathrm{Val}(\mathbb{Q}) \mid v(p) > 0\} = \{v_p\}$. Thus $v_p$ is a discrete point of $\mathrm{Val}(\mathbb{Q})$. On the other hand, $v_0(a) = 0$ for each $a \in \mathbb{Q}^\times$. So, if $\mathbf{B} = \bigcap_{i=1}^m \mathrm{Val}_{a_i}(\mathbb{Q}) \cap \bigcap_{j=1}^n \mathrm{Val}'_{b_j}(\mathbb{Q})$ contains $v_0$, then we may assume $m = 0$. Hence, $\mathbf{B}$ contains all $v_p$ with $p$ relatively prime to all denominators of $b_j$. This implies, every open neighborhood of $v_0$ consists of almost all elements of $\mathrm{Val}(\mathbb{Q})$. So, $\mathrm{Val}(\mathbb{Q})$ consists of a discrete sequence converging to $v_0$. In particular, $\mathrm{Val}(\mathbb{Q})$ is compact. ∎

The following result generalizes the last conclusion of Example 8.1.

PROPOSITION 8.2: $\mathrm{Val}(K)$ *is profinite.*

*Proof:* The space $\mathrm{Sign}(K) = \prod_{a \in K^\times} \{-1, 0\}$ with the product topology is a profinite space. For each $v \in \mathrm{Val}(K)$ and $a \in K^\times$ let $\mathrm{sign}(v(a))$ be $-1$ if $v(a) < 0$ and $0$ if $v(a) \geq 0$. Define a map $\sigma \colon \mathrm{Val}(K) \to \mathrm{Sign}(K)$ by $\sigma(v)(a) = \mathrm{sign}(v(a))$. It suffices to prove that $\sigma$ is a homeomorphism onto a closed subset of $\mathrm{Sign}(K)$.

Indeed, let $v, v' \in \mathrm{Val}(K)$ with $\sigma(v) = \sigma(v')$. Then $v(a) \geq 0$ if and only if $v'(a) \geq 0$. Hence, $v = v'$. Therefore, $\sigma$ is injective.

A basic open subset of $\sigma(\text{Val}(K))$ has the form

$$\{\sigma(v) \mid v \in \text{Val}(K),\ \text{sign}(v(a_i)) = -1,\ i = 1, \ldots, m,\ \text{sign}(v(b_j)) = 0,\ j = 1, \ldots, n\}$$

with $a_1, \ldots, a_m, b_1, \ldots, b_n \in K^\times$ and $m, n \geq 0$. It is the image of the basic open subset $\bigcap_{i=1}^{m} \text{Val}_{a_i^{-1}}(K) \cap \bigcap_{j=1}^{n} \text{Val}'_{b_j}(K)$. Therefore, $\sigma$ is a homeomorphism.

Next consider an element $f \in \text{Sign}(K)$ which belongs to the closure of $\text{Im}(\sigma)$. We construct $w \in \text{Val}(K)$ with $\sigma(w) = f$. This will conclude the proof of the proposition. Put $O = \{a \in K^\times \mid f(a) = 0\} \cup \{0\}$.

CLAIM: *$O$ is a valuation ring.* Indeed, assume $a, b \in O$ but $a + b \notin O$. Then $a, b \neq 0$ and $\{g \in \text{Sign}(K) \mid g(a) = 0,\ g(b) = 0,\ g(a+b) = -1\}$ is an open neighborhood of $f$. Hence, there is $v \in \text{Val}(K)$ with $\text{sign}(v(a)) = 0$, $\text{sign}(v(b)) = 0$, and $\text{sign}(v(a+b)) = -1$. Thus, $v(a) \geq 0$, $v(b) \geq 0$, and $v(a+b) < 0$. This contradiction proves that $O$ is closed under addition.

Similarly, $O$ is closed under multiplication and contains $0, -1$. Hence, $O$ is a subring of $K$.

Let now $a \in K^\times$. If $f(a) = f(a^{-1}) = -1$, there is $v \in \text{Val}(K)$ with $v(a) < 0$ and $v(a^{-1}) < 0$, a contradiction. Hence, $a \in O$ or $a^{-1} \in O$. Therefore, $O$ is a valuation ring.

Denote the valuation associated with $O$ by $w$. Then $\text{sign}(w(a)) = f(a)$ for each $a \in K$. Therefore, $f = \sigma(w)$. ∎

For a valued field $(K, v)$ and a polynomial $f(X) = \sum_{i=0}^{n} a_i X^i$ with $a_i \in K$ we write $v(f) = \min(v(a_0), \ldots, v(a_n))$. Also, for $\mathbf{x} = (x_1, \ldots, x_n) \in K^n$ we write $v(\mathbf{x}) = \min(v(x_1), \ldots, v(x_n))$.

LEMMA 8.3: *Let $(K, v)$ be a valued field, $(K_v, v_h)$ be a Henselian closure of $(K, v)$, and $L$ a finite separable extension of $K$. Then the following conditions are equivalent:*

(a) *There is a $K$-embedding of $L$ into $K_v$.*

(b) *$L/K$ has a primitive element $x$ such that $\text{irr}(x, K) = X^n + X^{n-1} + a_{n-2}X^{n-2} + \cdots + a_0$ with $v(a_i) > 0$, $i = 0, \ldots, n-2$.*

*Proof of "(a) $\Longrightarrow$ (b)":* The embedding of $L$ into $K_v$ induces a valuation $v_L$ of $L$ which extends $v$. Let $\hat{L}$ be the Galois closure of $L/K$. Choose an extension $\hat{v}$ of $v_L$ to $\hat{L}$. Then, $L$ is contained in the decomposition field $L'$ of $\hat{v}$ over $K$. Hence, $\mathrm{Gal}(\hat{L}/L') \le \mathrm{Gal}(\hat{L}/L)$.

Every extension of $v$ to $\hat{L}$ has the form $\hat{v} \circ \sigma$ with $\sigma \in \mathrm{Gal}(\hat{L}/K)$. We have, $\mathrm{res}_L(\hat{v} \circ \sigma) = \mathrm{res}_L \hat{v}$ if and only if there is $\tau \in \mathrm{Gal}(\hat{L}/L)$ with $\hat{v} \circ \sigma = \hat{v} \circ \tau$, that is, $\sigma\tau^{-1}$ lies in the decomposition group $\mathrm{Gal}(\hat{L}/L')$ of $\hat{v}$. Conclude: $\mathrm{res}_L(\hat{v} \circ \sigma) = \mathrm{res}_L \hat{v}$ if and only if $\sigma \in \mathrm{Gal}(\hat{L}/L)$.

Now use the Chinese remainder theorem [Jar, Lemma 6.7(c)] to find $y \in L$ with $\hat{v}(y) = 0$ and $\hat{v}(\sigma y) > 0$ for each $\sigma \in \mathrm{Gal}(\hat{L}/K) \smallsetminus \mathrm{Gal}(\hat{L}/L)$. Next choose a primitive element $z$ for $L/K$. Multiply $z$ by a suitable element of $K$ to assume

$$(2) \qquad \hat{v}(\sigma z) > \max(0, \hat{v}(y' - y))$$

for each $\sigma \in \mathrm{Gal}(\hat{L}/K)$ and every conjugate $y'$ of $y$ over $K$ with $y' \ne y$. Put $x = y + z$. Then

$$(3) \qquad \hat{v}(x) = 0 \text{ and } \hat{v}(\sigma x) > 0 \text{ for each } \sigma \in \mathrm{Gal}(\hat{L}/K) \smallsetminus \mathrm{Gal}(\hat{L}/L).$$

We prove $L = K(x)$.

To this end consider $\tau \in \mathrm{Gal}(\hat{L}/K(x))$. Then $\tau(y) - y = z - \tau(z)$. Therefore,

$$\hat{v}(\tau(y) - y) \ge \min(\hat{v}(z), \hat{v}(\tau(z))) \ge \min_{\sigma \in \mathrm{Gal}(\hat{L}/K)} \hat{v}(\sigma z).$$

By (2), $\tau(y) = y$. Hence, $\tau(z) = z$. Therefore, $L = K(z) \subseteq K(x) \subseteq L$. It follows that $L = K(x)$, as contended.

Let $x_1, \ldots, x_n$ be the conjugates of $x$ in $L$ with $x_1 = x$. For each $j \ge 2$ there is $\sigma \in \mathrm{Gal}(\hat{L}/K) \smallsetminus \mathrm{Gal}(\hat{L}/L)$ with $\sigma x = x_j$. Hence, by (3),

$$(4) \qquad \hat{v}(x_1) = 0 \text{ and } \hat{v}(x_j) > 0 \text{ if } j > 2.$$

Let $f(X) = X^n + b_{n-1}X^{n-1} + b_{n-2}X^{n-2} + \cdots + b_0 = \mathrm{irr}(x, K)$. By (4), $\hat{v}(b_{n-1}) = \hat{v}(x_1 + x_2 + \cdots + x_n) = 0$, $\hat{v}(b_{n-2}) = \hat{v}(\sum_{j \ne k} x_j x_k) > 0$, $\cdots$, $\hat{v}(b_0) = \hat{v}(x_1 \cdots x_n) > 0$. Obviously, $\frac{x}{b_{n-1}}$ is a primitive element for $L/K$. Its irreducible polynomial over $K$ is

$$X^n + X^{n-1} + \frac{b_{n-2}}{b_{n-1}^2}X^{n-2} + \cdots + \frac{b_0}{b_{n-1}^n} \ .$$

This polynomial has the required form.

*Proof of "(b) $\Longrightarrow$ (a)":* In the notation of (b) let $f = \text{irr}(x, K)$. Then $v(f(-1)) > 0$ and $v(f'(-1)) = v((-1)^{n-1}) = 0$. Hence, by Hensel's Lemma, $f$ has a root $x' \in K_v$. The map $x \mapsto x'$ extends to a $K$-embedding of $L$ into $K_v$. $\blacksquare$

Lemma 8.4 (Open map theorem): *Let $L$ be a field extension of $K$. Then the map $\text{res}_{L/K}\colon \text{Val}(L) \to \text{Val}(K)$ is continuous. If $L/K$ is separable algebraic, then the map is also open.*

*Proof:* By definition, $\text{res}_{L/K}^{-1}(\text{Val}_a(K)) = \text{Val}_a(L)$ and $\text{res}_{L/K}^{-1}(\text{Val}'_a(K)) = \text{Val}'_a(L)$ for each $a \in K^\times$. Hence, restriction of valuations of $L$ to $K$ is a continuous map.

Suppose now $L/K$ is Galois. Put $G = \text{Gal}(L/K)$. Then $G$ acts on $\text{Val}(L)$ continuously and $\text{res}_{L/K}$ induces a continuous bijective map $\rho\colon \text{Val}(L)/G \to \text{Val}(K)$. Since both spaces are profinite, $\rho$ is a homeomorphism. By definition, the canonical map $\pi\colon \text{Val}(L) \to \text{Val}(L)/G$ is open. Thus, $\text{res}_{L/K} = \rho \circ \pi$ is also open.

Finally suppose $L/K$ is separable algebraic. Let $\hat{L}$ be the Galois closure of $L/K$. Then $\text{res}_{\hat{L}/L}$ is continuous and $\text{res}_{\hat{L}/K}$ is open. Let $U$ be an open subset of $\text{Val}(L)$. Then $\text{res}_{L/K}(U) = \text{res}_{\hat{L}/K}(\text{res}_{\hat{L}/L}^{-1}(U))$ is open, as desired. $\blacksquare$

## 9. Locally Uniform $v$-adic Topologies

Every valuation $v$ of a field $K$ gives rise to a topology on $K$ which naturally extends to a topology on $V(K)$ (called the $v$-**topology**) for every variety $V$ defined over $K$. Polynomials $f \in KX]$ and in general morphisms between varieties over $K$ are continuous in the $v$-topology. The proof of continuity depends on only finitely many conditions of the form $v(a) > 0$ and $v(a') \geq 0$. So, it holds for all valuations $v'$ of $K$ satisfying the same conditions. In other words, polynomials are "locally uniform continuous". This observation holds even if we consider the polynomials as functions of valued fields extending $(K, v)$.

The aim of this section is the make this heuristic argument precise. It will be used in Proposition 12.4 to prove that every field-valuation structure satisfying the block approximation condition is unirationally closed.

We start by choosing a **large universal extension** of $K$. This is an algebraically closed field extension $\Omega$ of $K$ with $\mathrm{trans.deg}(\Omega/K) > \mathrm{card}(K)$. Denote the set of all field extensions $L$ of $K$ with $L \subseteq \Omega$ and $\mathrm{trans.deg}(L/K) \leq \mathrm{card}(K)$ by $\mathrm{Extend}(K)$. For each $v \in \mathrm{Val}(K)$ denote the class of all valued fields $(L, w)$ extending $(K, v)$ such that $L \in \mathrm{Extend}(K)$ by $\mathrm{Extend}(K, v)$. For each subset $\mathbf{B}$ of $\mathrm{Val}(K)$ let $\mathrm{Extend}(K, \mathbf{B}) = \bigcup_{v \in \mathbf{B}} \mathrm{Extend}(K, v)$. In addition, let $\mathrm{Hensel}(K, \mathbf{B})$ be the set of all Henselian fields $(L, w)$ in $\mathrm{Extend}(K, \mathbf{B})$.

The reason for working in $\Omega$ is to avoid using classes, especially to avoid operations with classes which may led to set theoretic paradoxes.

Denote the collection of all subsets of a set $A$ by $\mathrm{Subset}(A)$. Consider a reduced scheme of finite type $V$ over $K$ and a subset $\mathbf{B}$ of $\mathrm{Val}(K)$. Let

$$\mathrm{Set}(K, V, \mathbf{B}) = \prod_{(L,v) \in \mathrm{Extend}(K, \mathbf{B})} \mathrm{Subset}(V(L)).$$

Thus, each element of $\mathrm{Set}(K, V, \mathbf{B})$ is a set valued function $\mathcal{V}$ from $\mathrm{Extend}(K, \mathbf{B})$ satisfying $\mathcal{V}(L, v) \subseteq V(L)$ for each $(L, v) \in \mathrm{Extend}(K, \mathbf{B})$. Regard $V$ itself as an element of $\mathrm{Set}(K, V, \mathbf{B})$.

Let $\mathcal{V}, \mathcal{V}' \in \mathrm{Set}(K, V, \mathbf{B})$. We write $\mathcal{V} \subseteq \mathcal{V}'$ if $\mathcal{V}(L, v) \subseteq \mathcal{V}'(L, v)$ for all $(L, v) \in \mathrm{Extend}(K, \mathbf{B})$,

The **restriction** of $\mathcal{V}$ to a subset $\mathbf{B}_0$ of $\mathbf{B}$ is the function $\mathcal{V}|_{\mathbf{B}_0} \in \mathrm{Set}(K, V, \mathbf{B}_0)$ defined by $\mathcal{V}|_{\mathbf{B}_0}(L, v) = \mathcal{V}(L, v)$ for each $L \in \mathrm{Extend}(K, \mathbf{B}_0)$.

Define **unions** and **intersections** in $\mathrm{Set}(K, V, \mathbf{B})$ via unions and intersections of sets:

$$\big(\bigcup_{i \in I} \mathcal{V}_i\big)(L, v) = \bigcup_{i \in I} \mathcal{V}_i(L, v) \qquad \big(\bigcap_{i \in I} \mathcal{U}_i\big)(L, v) = \bigcap_{i \in I} \mathcal{U}_i(L, v)$$

These operations satisfy the usual de-Morgan laws. Similarly define the direct product of $\mathcal{U} \in \mathrm{Set}(K, U, \mathbf{B})$ with $\mathcal{V} \in \mathrm{Set}(K, V, \mathbf{B})$ by the rule $(\mathcal{U} \times \mathcal{V})(L, v) = \mathcal{U}(L, v) \times \mathcal{V}(L, v)$.

Let $\mathbf{a} \in K^n$, $\mathbf{c} \in (K^\times)^m$, and $f_1, \ldots, f_m \in K[X_1, \ldots, X_n]$. Define an element $\mathcal{O}_{\mathbf{a}, \mathbf{c}, \mathbf{f}, \mathbf{B}}$ in $\mathrm{Set}(K, \mathbb{A}^n, \mathbf{B})$ in the following way. For all $(L, v) \in \mathrm{Extend}(K, \mathbf{B})$

$$\mathcal{O}_{\mathbf{a}, \mathbf{c}, \mathbf{f}, \mathbf{B}}(L, v) = \{\mathbf{x} \in L^n \mid v(f_i(\mathbf{x}) - f_i(\mathbf{a})) > v(c_i), \ i = 1, \ldots, m\}$$

Note that $\mathcal{O}_{\mathbf{a}, \mathbf{c}, \mathbf{f}, \mathbf{B}}(L, v)$ is a $v$-open neighborhood of $\mathbf{a}$ in $L^n$. If we embed $K^n$ diagonally in $\prod_{(L,v) \in \mathrm{Extend}(K, \mathbf{B})} L^n$, then $\mathbf{a}$ belongs to $\prod_{(L,v) \in \mathrm{Extend}(K, \mathbf{B})} \mathcal{O}_{\mathbf{a}, \mathbf{c}, \mathbf{f}, \mathbf{B}}(L, v)$. So, we call $\mathbf{O}_{\mathbf{a}, \mathbf{c}, \mathbf{f}, \mathbf{B}}$ a **basic open neighborhood** of $\mathbf{a}$ in $\mathrm{Set}(K, \mathbb{A}^n, \mathbf{B})$. The intersection of finitely many basic open neighborhoods of $\mathbf{a}$ is again a basic open neighborhood of $\mathbf{a}$ in $\mathrm{Set}(K, \mathbb{A}^n, \mathbf{B})$. Define an **open neighborhood** of $\mathbf{a}$ in $\mathrm{Set}(K, \mathbb{A}^n, \mathbf{B})$ to be a union of basic open neighborhoods of $\mathbf{a}$ in $\mathrm{Set}(K, \mathbb{A}^n, \mathbf{B})$.

An example of an open neighborhood of $\mathbf{a}$ in $\mathrm{Set}(K, \mathbb{A}^n, \mathbf{B})$ is an **open ball**:

$$\mathcal{B}_{\mathbf{a}, c, \mathbf{B}}(L, v) = \{\mathbf{x} \in L^n \mid v(\mathbf{x} - \mathbf{a}) > v(c)\}.$$

Let $V$ be a Zariski closed subset of $\mathbb{A}^n_K$, $\mathbf{a} \in V(K)$, and $\mathcal{V}$ an open neighborhood of $\mathbf{a}$ in $\mathrm{Set}(K, \mathbb{A}^n, \mathbf{B})$. Refer to $V \cap \mathcal{V}$ as an **open neighborhood** of $\mathbf{a}$ in $\mathrm{Set}(K, V, \mathbf{B})$. ∎

*Remark 9.1:* Let $v \in \mathrm{Val}(K)$, $\mathbf{a} \in K^n$, and $c, d \in K$. Suppose $v(c) \le v(d)$. Then $\mathbf{B} = \{w \in \mathrm{Val}(K) \mid w(c) \le w(d)\}$ is an open neighborhood of $v$ in $\mathrm{Val}(K)$. Moreover, $\mathcal{B}_{\mathbf{a}, d, \mathbf{B}}(L, w) \subseteq \mathcal{B}_{\mathbf{a}, c, \mathbf{B}}(L, w)$ for all $(L, w) \in \mathrm{Extend}(K, \mathbf{B})$. ∎

*Definition 9.2: Uniform local topology on schemes.* Let $V$ be a Zariski closed subset in $\mathbb{A}^m_K$, $W$ a Zariski closed subset in $\mathbb{A}^n_K$, and $\varphi \colon V \to W$ be a $K$-morphism. Then there are polynomials $f_1, \ldots, f_n \in K[X_1, \ldots, X_m]$ with $\varphi(\mathbf{x}) = (f_1(\mathbf{x}), \ldots, f_n(\mathbf{x}))$ for all $L \in \mathrm{Extend}(K)$ and $\mathbf{x} \in V(L)$.

Let $\mathbf{B}$ be a subset of $\mathrm{Val}(K)$. For each $\mathcal{V} \in \mathrm{Set}(K, V, \mathbf{B})$ define $\varphi(\mathcal{V})$ to be the element of $\mathrm{Set}(K, W, \mathbf{B})$ given by $(\varphi(\mathcal{V}))(L, v) = \varphi(\mathcal{V}(L, v))$. Similarly, for each $\mathcal{W} \in \mathrm{Set}(K, W, \mathbf{B})$ define $\varphi^{-1}(\mathcal{W})$ to be the element of $\mathrm{Set}(K, V, \mathbf{B})$ defined by $\varphi^{-1}(\mathcal{W})(L, v) = \varphi^{-1}(\mathcal{W}(L, v))$.

As an example, let $\mathbf{a} \in V(K)$, $\mathbf{b} = \varphi(\mathbf{a})$, and $g_1, \ldots, g_k \in K[Y_1, \ldots, Y_n]$. Then $\mathbf{g} \circ \varphi = (h_1, \ldots, h_k)$ with $h_i(\mathbf{X}) = g_i(f_1(\mathbf{X}), \ldots, f_n(\mathbf{X}))$ and

$$\varphi^{-1}(W \cap \mathcal{O}_{\mathbf{b}, \mathbf{c}, \mathbf{g}, \mathbf{B}}) = V \cap \mathcal{O}_{\mathbf{a}, \mathbf{c}, \mathbf{g} \circ \varphi, \mathbf{B}}.$$

So, the inverse image under $\varphi$ of any open neighborhood of $\mathbf{b}$ in $\mathrm{Set}(K, W, \mathbf{B})$ is an open neighborhood of $\mathbf{a}$ in $\mathrm{Set}(K, V, \mathbf{B})$. In particular, if $\varphi$ is an isomorphism, $\mathcal{V}$ is an open neighborhood of $\mathbf{a}$ in $\mathrm{Set}(K, V, \mathbf{B})$, and $\mathcal{W} = \varphi(\mathcal{V})$, then $\mathcal{W}$ is an open neighborhood of $\mathbf{b}$ in $\mathrm{Set}(K, W, \mathbf{B})$ and $\varphi^{-1}(\mathcal{W}) = \mathcal{V}$.

Let now $V$ be a reduced scheme of finite type over $K$ and $\mathbf{a} \in V(K)$. Choose a Zariski $K$-open affine neighborhood $V_0$ of $\mathbf{a}$ in $V$. Each open neighborhood of $\mathbf{a}$ in $\mathrm{Set}(K, V_0, \mathbf{B})$ is an **open neighborhood** of $\mathbf{a}$ in $\mathrm{Set}(K, V, \mathbf{B})$. The observation of the preceding paragraph shows this definition is independent of $V_0$.  ∎

LEMMA 9.3 (Local uniform continuity of polynomials): *Let $(K, v)$ be a valued field $g \in K[X_1, \ldots, X_n]$, $\mathbf{a}, \mathbf{x} \in K^n$, and $e \in K^\times$. Suppose $v(g) \geq 0$, $v(\mathbf{a}, \mathbf{x}) \geq 0$, and $v(\mathbf{x} - \mathbf{a}) > v(e)$. Then $v(g(\mathbf{x}) - g(\mathbf{a}))) > v(e)$.*

*Proof:* We prove the Lemma by induction on $n$.

Suppose first $n = 1$. Write $g(X) = \sum_{i=0}^{r} c_i X^i$ with $c_i \in K$ satisfying $v(c_i) \geq 0$, $i = 0, \ldots, r$. Then

$$\begin{aligned}
v(g(x) - g(a)) &= v\Big(\sum_{i=0}^{r} c_i(x^i - a^i)\Big) \\
&\geq \min_{1 \leq i \leq r}(v(c_i) + v(x - a) + v(x^{i-1} + x^{i-2}a + \cdots + a^{i-1})) \\
&> v(e).
\end{aligned}$$

Assume now $n > 2$ and the statement holds up to $n - 1$. Then

$$v(g(\mathbf{x}) - g(\mathbf{a})) \geq \min\left(v(g(x_1, \ldots, x_{n-1}, x_n) - g(x_1, \ldots, x_{n-1}, a_n)),\right.$$

$$\left. v(g(x_1, \ldots, x_{n-1}, a_n) - g(a_1, \ldots, a_{n-1}, a_n))\right) > v(e).$$

This concludes the induction. ∎

As a consequence we show that open balls are locally basic open neighborhoods of $K$-rational points on varieties over $K$.

LEMMA 9.4: *Let $K$ be a field, $V$ a Zariski closed subset of $\mathbb{A}^n_K$, $\mathbf{a} \in V(K)$, $\mathbf{B}$ a closed subset of $\mathrm{Val}(K)$, and $\mathcal{V}$ an open neighborhood of $\mathbf{a}$ in $\mathrm{Set}(K, V, \mathbf{B})$. Then there is a partition $\mathbf{B} = \bigcup_{i=1}^m \mathbf{B}_i$ with $\mathbf{B}_i$ closed and for each $i$ there is an open ball $\mathcal{B}_{\mathbf{a}, c_i, \mathbf{B}_i}$ in $\mathrm{Set}(K, \mathbb{A}^n, \mathbf{B}_i)$ such that $V(L) \cap \mathcal{B}_{\mathbf{a}, c_i, \mathbf{B}_i}(L, v) \subseteq \mathcal{V}(L, v)$ for each $(L, v) \in \mathrm{Extend}(K, \mathbf{B}_i)$.*

*Proof:* Assume without loss $V = \mathbb{A}^n_K$. Choose $c'_1, \ldots, c'_l \in K^\times$ and $f_1, \ldots, f_l \in K[X_1, \ldots, X_n]$ such that $\mathcal{O}_{\mathbf{a}, \mathbf{c}', \mathbf{f}, \mathbf{B}}$ is an open neighborhood of $\mathbf{a}$ in $\mathcal{V}$. For each $v \in \mathbf{B}$ choose $e_v \in K^\times$ with $v(e_v \mathbf{a}) \geq 0$. Put $g_{vk}(\mathbf{X}) = f_k(\frac{1}{e_v} \mathbf{X})$, $k = 1, \ldots, l$. Next choose $d_v \in K^\times$ with $v(d_v g_{vk}) \geq 0$ for $k = 1, \ldots, l$. Finally choose $c_v \in K^\times$ with $v(c_v e_v) \geq 0$ and $v\left(\frac{c_v e_v}{d_v c'_k}\right) \geq 0$ for $k = 1, \ldots, l$. Then

$$\mathbf{B}'_v = \{v' \in \mathbf{B} \mid v'(e_v \mathbf{a}) \geq 0, \ v'(d_v g_{vk}) \geq 0, \ v'(c_v e_v) \geq 0, \ v'\left(\frac{c_v e_v}{d_v c'_k}\right) \geq 0, \ k = 1, \ldots, l\}$$

is an open neighborhood of $v$ in $\mathbf{B}$.

By Lemma 8.2, $\mathbf{B}$ is profinite. Hence, $\mathbf{B}'_v$ has a subset $\mathbf{B}_v$ which is open-closed in $\mathbf{B}$ and contains $v$. Compactness of $\mathbf{B}$ gives $v_1, \ldots, v_m \in \mathbf{B}$ with $\mathbf{B} = \bigcup_{i=1}^m \mathbf{B}_{v_i}$. Let $\mathbf{B}_1 = \mathbf{B}_{v_1}$ and $\mathbf{B}_i = \mathbf{B}_{v_i} \smallsetminus (\mathbf{B}_{v_1} \cup \cdots \cup \mathbf{B}_{v_{i-1}})$, $i = 2, \ldots, m$. Then $\mathbf{B}_i$ is closed in $\mathrm{Val}(K)$, $\mathbf{B}_i \subseteq \mathbf{B}'_{v_i}$, $i = 1, \ldots, m$, and $\mathbf{B} = \bigcup_{i=1}^m \mathbf{B}_i$.

Consider now an $i$ between $1$ and $m$. Put $c_i = c_{v_i}$, $d_i = d_{v_i}$, $e_i = e_{v_i}$, and $g_{ik} = g_{v_i k}$ for $k = 1, \ldots, l$. It suffices to prove that

$$\mathcal{B}_{\mathbf{a}, c_i, \mathbf{B}_i}(L, w) \subseteq \mathcal{O}_{\mathbf{a}, \mathbf{c}', \mathbf{f}, \mathbf{B}_i}(L, w)$$

for each $(L, w) \in \mathrm{Extend}(K, \mathbf{B}_i)$.

Indeed, our choices imply

(5)     $w(e_i \mathbf{a}) \geq 0, \ w(d_i g_{ik}) \geq 0, \ w(c_i e_i) \geq 0, \ w(c_i e_i) \geq w(d_i c_k'), \ k = 1, \ldots, l.$

Let $\mathbf{x} \in \mathcal{B}_{\mathbf{a}, \mathbf{c}, \mathbf{B}_i}(L, w)$. Then $w(\mathbf{x} - \mathbf{a}) > w(c_i)$. So, by (5), $w(e_i \mathbf{x} - e_i \mathbf{a}) > w(c_i e_i) \geq 0$. Hence, by (5), $w(e_i \mathbf{x}) \geq 0$. It follows from (5) and Lemma 9.3 that

$$w(d_i g_{ik}(e_i \mathbf{x}) - d_i g_{ik}(e_i \mathbf{a})) > w(c_i e_i) \geq w(d_i c_k'), \ k = 1, \ldots, l.$$

Thus, $w(f_k(\mathbf{x}) - f_k(\mathbf{a})) > w(c_k'), \ k = 1, \ldots, l.$ This means $\mathbf{x} \in \mathcal{B}_{\mathbf{a}, \mathbf{c}', \mathbf{f}, \mathbf{B}_i}(L, w)$, as claimed. ∎

## 10. Locally Uniform Hensel's Lemma

Let $(K, v)$ be a valued field, $\varphi\colon V \to W$ a morphism of absolutely irreducible varieties over $K$, $\mathbf{a} \in V_{\mathrm{simp}}(K)$, $\mathbf{b} \in W_{\mathrm{simp}}(K)$, and $\varphi(\mathbf{a}) = \mathbf{b}$. Suppose $\varphi$ is étale at $\mathbf{a}$. Let $(L, v)$ be a Henselian extension of $(K, v)$. Then $\mathbf{a}$ has a $v$-open neighborhood $\mathcal{V}$ in $V(L)$ and $\mathbf{b}$ has a $v$-open neighborhood $\mathcal{W}$ in $W(L)$ such that $\varphi\colon \mathcal{V}(L) \to \mathcal{W}(L)$ is a $v$-homeomorphism [GPR, Thm. 9.4]. The proof of this result relies on a higher dimensional Hensel's Lemma.

We strengthen this result by making $\mathcal{V}$ and $\mathcal{W}$ uniform on an open neighborhood of $v$ in $\mathrm{Val}(K)$. The proof reduces the general case to the case where $V$ is a hypersurface in $\mathbb{A}^{r+1}$, $W = \mathbb{A}^r$, and $\varphi$ is the projection on the first $r$ coordinates. Then we use a sharper form of Hensel's lemma.

LEMMA 10.1: *Let $(L, w)$ be a Henselian field and $f \in L[T_1, \ldots, T_r, X]$ monic in $X$. Put $f' = \frac{\partial f}{\partial X}$. Assume $w(f) \geq 0$ (hence $w(f') \geq 0$). Let $\mathbf{b}_0, \mathbf{b} \in L^r$, $c_0 \in L$, and $\varepsilon \geq \delta \geq 0$ be in $w(L^\times)$. Suppose*

(1a) $$w(\mathbf{b}_0, c_0) \geq 0,$$

(1b) $$w\big(f'(\mathbf{b}_0, c_0)\big) = \delta,$$

(1c) $$w\big(f(\mathbf{b}_0, c_0)\big) > \delta + \varepsilon, \text{ and}$$

(1d) $$w(\mathbf{b} - \mathbf{b}_0) > \delta + \varepsilon.$$

*Then $w(\mathbf{b}) \geq 0$ and there is a unique $c \in L$ with $f(\mathbf{b}, c) = 0$ and $w(c - c_0) > \varepsilon$. In particular, $w(c) \geq 0$ and $w\big(f'(\mathbf{b}, c)\big) = \delta$.*

*Proof:* By (1a) and (1d), $w(\mathbf{b}) \geq 0$. By (1d) and Lemma 9.3

(2a) $$w\big(f'(\mathbf{b}, c_0) - f'(\mathbf{b}_0, c_0)\big) > \delta + \varepsilon \geq \delta,$$

(2b) $$w\big(f(\mathbf{b}, c_0) - f(\mathbf{b}_0, c_0)\big) > \delta + \varepsilon,$$

Hence by (1a) and (1c)

(3) $$w\big(f'(\mathbf{b}, c_0)\big) = \delta, \qquad w\big(f(\mathbf{b}, c_0)\big) > \delta + \varepsilon = 2\delta + (\varepsilon - \delta).$$

A sharp form of Hensel's lemma [Jar, Prop. 11.1(e)] gives a unique $c \in L$ such that $f(\mathbf{b}, c) = 0$ and $w(c - c_0) > \delta + (\varepsilon - \delta) = \varepsilon \geq \delta$. By (1a), $w(c) \geq 0$. By (1d) and Lemma 9.3, $w(f'(\mathbf{b}, c) - f'(\mathbf{b}_0, c_0)) > \delta$. Hence, by (1b), $w\big(f'(\mathbf{b}, c)\big) = \delta$. ∎

For each $f \in K[X_1, \ldots, X_n]$ let $V(f)$ be the hypersurface in $\mathbb{A}^n$ defined by $f = 0$.

LEMMA 10.2: *Let* $f \in K[T_1, \ldots, T_r, X]$, $v \in \mathrm{Val}(K)$, *and* $(\mathbf{b}_0, c_0) \in K^{r+1}$. *Put* $V = V(f)$ *and* $f' = \frac{\partial f}{\partial X}$. *Suppose* $f$ *is monic in* $X$,

$$(4) \qquad v(f) \geq 0, \ v(\mathbf{b}_0, c_0) \geq 0, \quad and \quad v(f(\mathbf{b}_0, c_0)) > 2v(f'(\mathbf{b}_0, c_0)).$$

*Then* $v$ *has an open neighborhood* $\mathbf{B}$ *in* $\mathrm{Val}(K)$, $\mathbf{b}_0$ *has an open neighborhood* $\mathcal{B}$ *in* $\mathrm{Set}(K, \mathbb{A}^r, \mathbf{B})$, *and* $c_0$ *has an open neighborhood* $\mathcal{C}$ *in* $\mathrm{Set}(K, \mathbb{A}^1, \mathbf{B})$ *satisfying this: For each* $(L, w) \in \mathrm{Hensel}(K, \mathbf{B})$ *the projection*

$$(5) \qquad\qquad \mathrm{pr}\colon (\mathcal{B}(L) \times \mathcal{C}(L)) \cap V(L) \to \mathcal{B}(L)$$

*is a* $w$-*homeomorphism.*

Proof: The sharp inequality in (4) implies $f'(\mathbf{b}_0, c_0) \neq 0$. Hence,

$$\mathbf{B} = \{w \in \mathrm{Val}(K) \mid\ w(f) \geq 0, \ w(\mathbf{b}_0, c_0) \geq 0, \ \text{and} \ w(f(\mathbf{b}_0, c_0)) > 2w(f'(\mathbf{b}_0, c_0))\}$$

is an open neighborhood of $v$ in $\mathrm{Val}(K)$.

Consider $(L, w) \in \mathrm{Hensel}(K, \mathbf{B})$. Let $\delta = w(f'(\mathbf{b}_0, c_0))$. Then

$$w(f) \geq 0, \quad w(\mathbf{b}_0, c_0) \geq 0, \quad \text{and} \quad w(f(\mathbf{b}_0, c_0)) > 2\delta.$$

Let

$$\mathcal{B}(L) = \{\mathbf{b} \in L^r \mid\ w(\mathbf{b} - \mathbf{b}_0) > 2\delta\} \quad \text{and} \quad \mathcal{C}(L) = \{c \in L \mid\ 2(c - c_0) > \delta\}.$$

By Lemma 5.1 (with $\varepsilon = \delta$) the map $\mathrm{pr}$ in (5) is bijective. As a projection map, $\mathrm{pr}$ is continuous. We prove that $\mathrm{pr}^{-1}$ is continuous.

Consider $\mathbf{b}_1 \in \mathcal{B}(L)$. Let $c_1$ be the unique element of $L$ with $(\mathbf{b}_1, c_1) \in \big(\mathcal{B}(L) \times \mathcal{C}(L)\big) \cap V(L)$. Let $\varepsilon \in w(L^\times)$ with $\delta \leq \varepsilon$. By Lemma 10.1, $w(\mathbf{b}_1, c_1) \geq 0$ and $w(f'(\mathbf{b}_1, c_1)) = \delta$. Let $\mathbf{b} \in \mathcal{B}(L)$ with $w(\mathbf{b} - \mathbf{b}_1) > \delta + \varepsilon$. Then the unique element $c \in L$ which Lemma 10.1 (with $(\mathbf{b}_1, c_1)$ replacing $(\mathbf{b}_0, c_0)$) gives satisfies $f(\mathbf{b}, c) = 0$ and $w(c - c_1) > \varepsilon$. In particular, $c \in \mathcal{C}(L)$, $\mathrm{pr}(\mathbf{b}, c) = \mathbf{b}$, and $w(\mathbf{b}, c) - (\mathbf{b}_1, c_1)) > \varepsilon$, as desired. ∎

PROPOSITION 10.3: *Let* $\varphi \colon V \to W$ *be a morphism of absolutely irreducible varieties over* $K$, $v \in \mathrm{Val}(K)$, $\mathbf{a} \in V_{\mathrm{simp}}(K)$, *and* $\mathbf{b} \in W_{\mathrm{simp}}(K)$. *Suppose* $\varphi$ *is étale at* $\mathbf{a}$ *and* $\varphi(\mathbf{a}) = \mathbf{b}$. *Then* $v$ *has an open neighborhood* $\mathbf{B}_v$ *in* $\mathrm{Val}(K)$, $\mathbf{a}$ *has an open neighborhood* $\mathcal{V}_v$ *in* $\mathrm{Set}(K, V, \mathbf{B}_v)$, *and* $\mathbf{b}$ *has an open neighborhood* $\mathcal{W}_v$ *in* $\mathrm{Set}(K, W, \mathbf{B}_v)$ *satisfying this: For each* $(L, w) \in \mathrm{Hensel}(K, \mathbf{B}_v)$ *the map* $\varphi \colon \mathcal{V}_v(L) \to \mathcal{W}_v(L)$ *is a* $w-$*homeomorphism.*

*Proof:* Let $r = \dim(W) = \dim(V)$.

PART A: *Suppose* $W = \mathbb{A}^r$. By [Ray, p. 60], $\varphi$ is locally standard étale. That is, there are a Zariski $K$-open neighborhood $A$ of $\mathbf{b}$ in $\mathbb{A}^r$, a Zariski $K$-open affine neighborhood $V_0$ of $\mathbf{a}$ in $V$, a polynomial $f \in K[T_1, \ldots, T_r, X]$ which is monic in $X$ (and absolutely irreducible), an element $c \in K$, and an isomorphism $\theta \colon V_0 \to (A \times \mathbb{A}^1) \cap V(f)$ over $K$ with $f(\mathbf{b}, c) = 0$, $\frac{\partial f}{\partial X}(\mathbf{b}, c) \neq 0$, $\theta(\mathbf{a}) = (\mathbf{b}, c)$, $\varphi(V_0) = A$, and $\mathrm{pr} \circ \theta = \varphi$. Multiply the point $(\mathbf{b}, c)$ with an appropriate element $u$ of $K^\times$ and the coefficients of $f$ with powers of $u$ and replace $\theta$ by $\theta \circ \mu_u$, where $\mu_u(\mathbf{x}, y) = (u\mathbf{x}, uy)$, to assume $v(f) \geq 0$ and $v(\mathbf{b}, c) \geq 0$.

Lemma 10.2 gives an open neighborhood $\mathbf{B}$ of $v$ in $\mathrm{Val}(K)$, an open neighborhood $\mathcal{B}$ of $\mathbf{b}$ in $\mathrm{Set}(K, \mathbb{A}^r, \mathbf{B})$, an open neighborhood $\mathcal{C}$ of $c$ in $\mathrm{Set}(K, \mathbb{A}^1, \mathbf{B})$ satisfying this:

(6) For each $(L, w) \in \mathrm{Hensel}(K, \mathbf{B})$ the projection $\mathrm{pr} \colon (\mathcal{B}(L) \times \mathcal{C}(L)) \cap V(f)(L) \to \mathcal{B}(L)$ is a $w$-homeomorphism.

Replace $\mathcal{B}$ by $\mathcal{B} \cap A$, if necessary, to assume $\mathcal{B} \subseteq A$.

Put $\mathcal{V} = \theta^{-1}\big((\mathcal{B} \times \mathcal{C}) \cap V(f)\big)$. By Definition 9.2, $\mathcal{V}$ is an open neighborhood of $\mathbf{a}$ in $\mathrm{Set}(K, V, \mathbf{B})$. Also, for each $(L, w) \in \mathrm{Hensel}(K, \mathrm{Val}(K))$ the map $\theta \colon \mathcal{V}(L) \to (\mathcal{B}(L) \times \mathcal{C}(L)) \cap V(f)(L)$ is a $w$-homeomorphism. If, in addition, $w|_K \in \mathbf{B}$, (6) implies the map $\varphi \colon \mathcal{V}(L) \to \mathcal{B}(L)$ is a $w$-homeomorphism.

PART B: *The general case.* As **b** is simple on $W$, the maximal ideal $\mathfrak{m}_{W,\mathbf{b}}$ of the local ring of $W$ has $r$ generators $t_1, \ldots, t_r$, $\tau = (t_1, \ldots, t_r)$ is an étale map of $W$ into $\mathbb{A}^r$ at **b** and $\tau(\mathbf{b}) = \mathbf{o} = (0, \ldots, 0)$ [Mum, p. 255, Thm. 1]. Part A gives an open neighborhood $\mathbf{B}_1$ of $v$ in $\mathrm{Val}(K)$, an open neighborhood $\mathcal{W}_1$ of **b** in $\mathrm{Set}(K, W, \mathbf{B}_1)$, and an open neighborhood $\mathcal{A}_1$ of **o** in $\mathrm{Set}(K, \mathbb{A}^r_K, \mathbf{B}_1)$ satisfying this: For each $(L, w) \in \mathrm{Hensel}(K, \mathbf{B}_1)$ the map $\tau \colon \mathcal{W}_1(L) \to \mathcal{A}_1(L)$ is a $w$-homeomorphism.

By [Hrt, p. 268, Prop. 10.1(b)], $\tau \circ \varphi$ is an étale morphism of $V$ into $\mathbb{A}^r_K$ at **a**. Part A gives an open neighborhood $\mathbf{B}_2$ of $v$ in $\mathrm{Val}(K)$, an open neighborhood $\mathcal{V}_2$ of **a** in $\mathrm{Set}(K, V)$, and an open neighborhood $\mathcal{A}_2$ of **o** in $\mathrm{Set}(K, \mathbb{A}^r_K)$ satisfying this: For all $(L, w) \in \mathrm{Hensel}(K, \mathbf{B}_2)$ the map $\tau \circ \varphi \colon \mathcal{V}_2(L) \to \mathcal{A}_2(L)$ is a $w$-homeomorphism.

Let $\mathbf{B} = \mathbf{B}_1 \cap \mathbf{B}_2$, $\mathcal{A} = \mathcal{A}_1 \cap \mathcal{A}_2$, $\mathcal{W} = \mathcal{W}_1 \cap \tau^{-1}(\mathcal{A}_2)$, and $\mathcal{V} = \varphi^{-1}(\mathcal{W})$. Then $\mathbf{B}, \mathcal{V}, \mathcal{W}$ satisfy the requirements of the lemma. ∎

COROLLARY 10.4: *Let* $\varphi \colon V \to W$ *be a morphism of absolutely irreducible varieties over* $K$, $\mathbf{a} \in V_{\mathrm{simp}}(K)$, $\mathbf{b} \in W_{\mathrm{simp}}(K)$, *and* $\mathbf{B}$ *a closed subset of* $\mathrm{Val}(K)$. *Suppose* $\varphi$ *is étale at* **a** *and* $\varphi(\mathbf{a}) = \mathbf{b}$. *Then there is a partition* $\mathbf{B} = \bigcup_{i=1}^n \mathbf{B}_i$ *with* $\mathbf{B}_i$ *closed,* **a** *has an open neighborhood* $\mathcal{V}_i$ *in* $\mathrm{Set}(K, V, \mathbf{B}_i)$, *and* **b** *has an open neighborhood* $\mathcal{W}_i$ *in* $\mathrm{Set}(K, W, \mathbf{B}_i)$, $i = 1, \ldots, n$, *satisfying this: For all* $i$, $w \in \mathbf{B}_i$, *and* $(L, w) \in \mathrm{Hensel}(K, w)$ *the map* $\varphi \colon \mathcal{V}_i(L) \to \mathcal{W}_i(L)$ *is a* $w$–*homeomorphism.*

*Proof:* For each $v \in \mathbf{B}$ let $\mathbf{B}_v$, $\mathcal{V}_v$, and $\mathcal{W}_v$ as in Proposition 10.3. Choose an openclosed subset $\mathbf{B}'_v$ of $\mathrm{Val}(K)$ with $v \in \mathbf{B}'_v \subseteq \mathbf{B}_v$. Then, the collection of all $\mathbf{B}'_v$ is an open covering of $\mathbf{B}$. As $\mathbf{B}$ is closed in $\mathrm{Val}(K)$ and $\mathrm{Val}(K)$ is compact (Proposition 8.2), $\mathbf{B}$ is compact. Thus there are $v_1, \ldots, v_n \in \mathbf{B}$ with $\mathbf{B} = \bigcup_{i=1}^n \mathbf{B}'_{v_i}$. Let $\mathbf{B}_i = B'_{v_i} \smallsetminus B'_{v_1} \cup \cdots \cup B'_{v_{i-1}}$, $\mathcal{V}_i = \mathcal{V}_{v_i}$, and $\mathcal{W}_i = \mathcal{W}_{v_i}$. They satisfy the conclusion of the corollary. ∎

## 11. Field-Valuation Structures

We extend field structures to "field-valuation structures" by quipping each local field with a valuation.

A **field-valuation structure** is a structure $\mathbf{K} = (K, X, K_x, v_x)_{x \in X}$ satisfying the following conditions:

(1a) $(K, X, K_x)_{x \in X}$ is a field structure. Thus, for each finite separable extension $L$ of $K$ the set $X_L = \{x \in X \mid L \subseteq K_x\}$ is open in $X$.

(1b) $v_x$ is a valuation of $K_x$ satisfying $v_{x^\sigma} = v_x^\sigma$ for all $x \in X$ and $\sigma \in \mathrm{Gal}(K)$. Here $v_x^\sigma(u^\sigma) = v_x(u)$ for each $u \in K_x$.

(1c) For each finite separable extension $L$ of $K$ define a map $\nu_L \colon X_L \to \mathrm{Val}(L)$ by $\nu_L(x) = v_x|_L$. Then $\nu_L$ is continuous.

The **absolute Galois structure** associated with $\mathbf{K}$ is the same associated with the underlying field structure, namely $\mathrm{Gal}(\mathbf{K}) = (\mathrm{Gal}(K), X, \mathrm{Gal}(K_x))_{x \in X}$. We call $\mathbf{K}$ **proper** if $\mathrm{Gal}(\mathbf{K})$ is proper. Call $\mathbf{K}$ **Henselian** if $(K_x, v_x)$ is Henselian for each $x \in X$.

Denote the maximal purely inseparable extension of a field $K$ by $K_{\mathrm{ins}}$.

LEMMA 11.1: *Let $\mathbf{K} = (K, X, K_x, v_x)_{x \in X}$ be a field-valuation structure.*

(a) *Let $K'$ be a separable algebraic extension of $K$ and $X'$ a closed subset of $X$. Suppose $X'$ is closed under the action of $\mathrm{Gal}(K')$ and $K' \subseteq K_x$ for each $x \in X'$. Then $\mathbf{K}' = (K', X', K_x, v_x)_{x \in X'}$ is a field-valuation structure.*

(b) *For each $x \in X$ let $v_{x,\mathrm{ins}}$ be the unique extension of $v_x$ to $K_{x,\mathrm{ins}}$. Then $\mathbf{K}_{\mathrm{ins}} = (K_{\mathrm{ins}}, X, K_{x,\mathrm{ins}}, v_{x,\mathrm{ins}})$ is a field-valuation structure. Moreover, there is an isomorphism* res: $\mathrm{Gal}(\mathbf{K}_{\mathrm{ins}}) \to \mathrm{Gal}(\mathbf{K})$ *of group structures.*

*Proof of (a):* By Remark 2.6, $(K', X, K_x)_{x \in X'}$ is a field structure. It remains to prove that $\nu_{L'} \colon X'_{L'} \to \mathrm{Val}(L')$ is continuous for each finite separable extension $L'$ of $K'$. It suffices to consider $u \in L'$ and to prove that each of the sets $Y = \{x \in X'_{L'} \mid v_x(u) > 0\}$ and $Y' = \{x \in X'_{L'} \mid v_x(u) \geq 0\}$ is open in $X'$. To this end choose a finite separable extension $L$ of $K$ containing $u$ with $L' = K'L$. Then $Y = X' \cap \{x \in X_L \mid v_x(u) > 0\}$, so $Y$ is open by (1c). Similarly, $Y'$ is open.

*Proof of (b):* It suffices to consider the case when $p = \mathrm{char}(K) > 0$. Let $L'$ be a finite extension of $K_{\mathrm{ins}}$ and $u \in L'$. Put $L = K_s \cap L'$. Then $L_{\mathrm{ins}} = L'$ and there is a power $q$ of $p$ with $u^q \in L$. Thus, $\{x \in X \mid L' \subseteq K_{x,\mathrm{ins}}\} = \{x \in X \mid L \subseteq K_x\}$ is open. Also, $v_{x,\mathrm{ins}}(u) = \frac{1}{q} v_x(u^q)$. This implies, $\mathbf{K}_{\mathrm{ins}}$ is a field-valuation structure. ∎

When all $(K_x, v_x)$ are Henselian, we may replace Condition (1c) by a more convenient condition:

Lemma 11.2: *Let $(K, X, K_x)_{x \in X}$ be a field structure. For each $x \in X$ let $v_x$ be a Henselian valuation on $K_x$ such that $v_{x^\sigma} = v_x^\sigma$ for all $x \in X$ and $\sigma \in \mathrm{Gal}(K)$. Extend each $v_x$ to $K_s$ in the unique possible way. Then $(K, X, K_x, v_x)_{x \in X}$ is a field-valuation structure if and only if*

(2) *the map $\nu\colon X \to \mathrm{Val}(K_s)$ defined by $x \mapsto v_x$ is continuous.*

*Proof:* By the uniqueness of the extension of $v_x$ from $K_x$ to $K_s$, the equality $v_{x^\sigma} = v_x^\sigma$ holds in $\mathrm{Val}(K_s)$ for all $x \in X$ and $\sigma \in \mathrm{Gal}(K)$.

Field-valuation structure implies (2): Let $u \in K_s^\times$ and let $x \in X$. We have to show that if $v_x(u) > 0$ (resp. $v_x(u) \geq 0$), and $x' \in X$ is sufficiently close to $x$, then $v_{x'}(u) > 0$ (resp. $v_{x'}(u) \geq 0$).

Let $f(X) = X^n + a_{n-1}X^{n-1} + \cdots + a_0$ be the irreducible polynomial of $u$ over $K_x$. Then

(3) $$u = -a_{n-1} - a_{n-2}u^{-1} - \cdots - a_0(u^{-1})^{n-1}.$$

Let $u_1, \ldots, u_n$ be the roots of $f$ in $K_s$. Since $v_x$ uniquely extends to $K_s$, we have

(4) $$v_x(u_1) = \cdots = v_x(u_n) = v_x(u).$$

Let $N/K$ be a finite separable extension of $K$ containing $u_1, \ldots, u_n$. Put $L = N \cap K_x$. Then $a_0, \ldots, a_{n-1} \in L$. We distinguish between two cases.

(a) Suppose $v_x(u) > 0$. We have $L \subseteq K_x$. Since $a_0, \ldots, a_{n-1}$ are the elementary symmetric functions in $u_1, \ldots, u_n$, (4) implies that $v_x(a_0), \ldots, v_x(a_{n-1}) > 0$. Hence, by (1), if $x' \in X$ is sufficiently close to $x$, then $L \subseteq K_{x'}$ and $v_{x'}(a_0), \ldots, v_{x'}(a_{n-1}) > 0$. It follows $v_{x'}(u) > 0$. Indeed, if $v_{x'}(u) \leq 0$, then $v_{x'}(u^{-1}) \geq 0$, hence, by (3), $v_{x'}(u) > 0$, a contradiction.

(b) Suppose $v_x(u) \geq 0$. We have $L \subseteq K_x$. By (4), $v_x(a_0), \ldots, v_x(a_{n-1}) \geq 0$. Hence, by (1), if $x' \in X$ is sufficiently close to $x$, then $L \subseteq K_{x'}$ and $v_{x'}(a_0), \ldots, v_{x'}(a_{n-1}) \geq \blacksquare$ 0. It follows that $v_{x'}(u) \geq 0$. Indeed, if $v_{x'}(u) < 0$, then $v_{x'}(u^{-1}) > 0$, and by (3), $v_{x'}(u) \geq 0$, a contradiction.

(2) IMPLIES FIELD-VALUATION STRUCTURE: Let $L$ be a finite separable extension of $K$. Then, $\nu\colon X \to \mathrm{Val}(K_s)$ and res: $\mathrm{Val}(K_s) \to \mathrm{Val}(L)$ are continuous (Lemma 8.4). Hence, $\nu_L = \mathrm{res} \circ \nu|_{X_L}$ is continuous. $\quad\blacksquare$

## 12. Block Approximation

Let $\mathbf{K} = (K, X, K_x, v_x)_{x \in X}$ be a field valuation structure. Put $\mathcal{K} = \{K_x \mid x \in X\}$. Suppose $K$ is P$\mathcal{K}$C, $X$ has only finitely many $\mathrm{Gal}(K)$-orbits, and the restriction of the corresponding valuations to $K$ are independent. An application of the local homeomorphism theorem [GPR, Thm. 9.4] for varieties over Henselian fields and the weak approximation theorem prove that $K$ is unirationally closed [HaJ3, Prop. 3.2]. In the general case, when $X$ has possibly infinitely many $\mathrm{Gal}(K)$-orbits, the block approximation condition substitutes all three conditions. It says roughly that finitely many algebraic points of a variety $V$ over $K$, each associated with an open-closed subset of $X$ (a "block") can be simultanuously approximated within the block by a single $K$-rational point of $V$. Here is the precise definition:

*Definition 12.1: Block approximation condition.* A **block approximation problem** for a field-valuation structure $\mathbf{K} = (K, X, K_x, v_x)_{x \in X}$ is a data $(V, X_i, L_i, \mathbf{a}_i, c_i)_{i \in I_0}$ satisfying this:

(1a) $(\mathrm{Gal}(L_i), X_i)_{i \in I_0}$ is a special partition of $\mathrm{Gal}(\mathbf{K})$.

(1b) $V$ is a smooth affine variety over $K$.

(1c) $\mathbf{a}_i \in V(L_i)$.

(1d) $c_i \in K^\times$.

An analogous condition where valuations are replaced by orderings appears in [Pre, p. 354] and [FHV, Prop. 1.2].

A **solution** of the problem is a point $\mathbf{a} \in V(K)$ with $v_x(\mathbf{a} - \mathbf{a}_i) > v_x(c_i)$ for all $i \in I_0$ and $x \in X_i$. We say $\mathbf{K}$ satisfies the **block approximation condition** if each block approximation problem for $\mathbf{K}$ is solvable. ∎

The block approximation condition has several interesting consequences.

*Definition 12.2: Pseudo-$\mathcal{K}$-closed fields.* Let $K$ be a field and $\mathcal{K}$ a set of field extensions of $K$. We say $K$ is **P$\mathcal{K}$C** if this holds: Every smooth absolutely irreducible variety $V$ over $K$ with a $K'$-rational point for each $K' \in \mathcal{K}$ has a $K$-rational point. ∎

PROPOSITION 12.3: *Let $\mathbf{K} = (K, X, K_x, v_x)_{x \in X}$ be a Henselian field-valuation structure satisfying the block approximation condition.*

(a) *Put $\mathcal{K} = \{K_x \mid x \in X\}$. Then $K$ is P$\mathcal{K}$C.*

(b) *Let $x_1, \ldots, x_n \in X$ lie in distinct $\mathrm{Gal}(K)$-orbits. Then $v_{x_1}|_K, \ldots v_{x_n}|_K$ satisfies the weak approximation theorem.*

(c) *Suppose $x, y \in X$ lie in distinct $\mathrm{Gal}(K)$-orbits. Then $v_x|_K$ and $v_y|_K$ are independent.*

(d) *Suppose $X$ has more than one $\mathrm{Gal}(K)$-orbit. Then the trivial valuation is not in $\nu_K(X)$.*

(e) *For each $x \in X$, $K$ is $v_x$-dense in $K_x$; and*

(f) *$(K_x, v_x)$ is a Henselian closure of $(K, v_x|_K)$.*

(g) *Suppose $K_x \neq K_s$. Then $\mathrm{Aut}(K_x/K) = 1$.*

*Proof of (a):* Let $V$ be a smooth absolutely irreducible variety over $K$ with a point $\mathbf{a}_x \in V(K_x)$ for each $x \in X$. Then $\mathrm{Gal}(K(\mathbf{a}_x))$ is an open subgroup of $\mathrm{Gal}(K)$ containing $\mathrm{Gal}(K_x)$. Lemma 3.6 gives a special partition $(\mathrm{Gal}(K(\mathbf{a}_{x_i})), X_i)_{i \in I_0}$ with $x_i \in X_i$ for each $i \in I_0$. Thus, $(V, X_i, K(\mathbf{a}_{x_i}), \mathbf{a}_{x_i}, 1)_{i \in I_0}$ is a block approximation problem for $\mathbf{K}$. Our assumption gives a point $\mathbf{a} \in V(K)$. It follows, $K$ is P$\mathcal{K}$C.

*Proof of (b):* Put $v_i = v_{x_i}|_K$, $i = 1, \ldots, n$. Let $a_i, c_i$ be elements of $K$ with $c_i \neq 0$. Since $X/\mathrm{Gal}(K)$ is profinite, there are open-closed distinct $\mathrm{Gal}(K)$-invariant subsets $X_1, \ldots, X_n$ of $X$ with $x_i \in X_i$, $i = 1, \ldots, n$. Let $I = \{0, 1, \ldots, n\}$ and $X_0 = X \smallsetminus X_1 \cup \cdots \cup X_n$. Then let $a_0 = 0$ and $c_0 = 1$. Then $(\mathbb{A}^1_K, X_i, K, a_i, c_i)_{i \in I_0}$ is a block approximation problem for $\mathbf{K}$.

By assumption, there is $a \in K$ with $v_i(a - a_i) > v_i(c_i)$, $i = 1, \ldots, n$. It follows, $v_1, \ldots, v_n$ satisfy the weak approximation theorem.

*Proof of (c):* Use (b).

*Proof of (d):* Assume $v_0 = v_x|_K$ is trivial for some $x \in X$. Choose $y \in X$ outside the $\mathrm{Gal}(K)$-orbit of $x$. By (c), $v_1 = v_y|_K$ is nontrivial. So, there is $a_1 \in K$ with $v_1(a_1) < 0$. Statement (b) gives $a \in K$ with $v_0(a - a_1) > 0$ and $v_1(a) > 0$. By the first inequality, $a = a_1$. So, by the second inequality, $v_1(a_1) > 0$, in contradiction to the choice of $a_1$.

*Proof of (e):* Let $x \in X$, $a_1 \in K_x$, and $c_1 \in K^\times$. We have to find $a \in K$ satisfying $v_x(a - a_1) > v_x(c_1)$.

By assumption $K(a_1) \leq K_x$ and the stabilizer of $x$ is contained in $\mathrm{Gal}(K_x)$. Therefore, Lemma 3.4 gives an open-closed neighborhood $X_1$ of $x$ in $X$ which is invariant under $\mathrm{Gal}(K(a_1))$ such that $X_1^{\mathrm{Gal}(K)} = \bigcup_{\rho \in R_1} X_1^\rho$ for each $R_1 \subseteq \mathrm{Gal}(K)$ with $\mathrm{Gal}(K) = \bigcup_{\rho \in R_1} \mathrm{Gal}(K(a_1))\rho$. Thus, $\mathrm{Gal}(K(a_1)) = \{\sigma \in \mathrm{Gal}(K) \mid X_1^\sigma = X_1\}$. Put $L_1 = K(a_1)$.

Let $I_0 = \{0, 1\}$, $X_0 = X \smallsetminus X_1^{\mathrm{Gal}(K)}$, $a_0 = 0$, $c_0 = 1$, and $L_0 = K$. Then $(\mathbb{A}^1_K, X_i, L_i, a_i, c_i)_{i=0,1}$ is a block approximation problem for $\mathbf{K}$.

By assumption, there is $a \in K$ with $v_x(a - a_1) > v_x(c_1)$, as desired.

*Proof of (f):* By assumption, $(K_x, v_x)$ is Henselian. Choose a Henselian closure $(K', v_x)$ of $(K, v_x|_K)$ in $(K_x, v_x)$. Consider $a \in K_x$. Let $a_1, \ldots, a_n$ be the conjugates of $a$ over $K'$. By (e) there is $b \in K$ with $v_x(b - a) > \max_{i \neq j} v_x(a_i - a_j)$. Hence, by Krasner's Lemma [Jar, Lemma 12.1], $K'(a) \subseteq K'(b) = K'$. Therefore, $K_x = K'$.

*Proof of (g):* Let $\sigma \in \mathrm{Aut}(K_x/K)$. Then both $v_x$ and $v_x^\sigma$ are Henselian valuations of $K_x$. Therefore, $K$ has a nontrivial valuation $w$ which is coarser than both $v_x$ and $v_x^\sigma$ [Jar, Lemma 13.2]. In particular, the $v$-topology of $K$ coincides with the $w$-topology of $K$ [Jar, Lemma 3.2]. Hence, by (e), $K$ is $w$-dense in $K_x$.

Assume there exists $b \in K_x$ with $b \neq b^\sigma$. Then there exists $c \in K^\times$ with $v(c) > v(b - b^{\sigma^{-1}})$ and there exists $a \in K$ with $w(a - b) > w(c)$. Since $w$ is coarser than both $v$ and $v^\sigma$, we have $v(a - b) > v(c)$ and $v^\sigma(a - b) > v(c)$. Hence, $v(a - b^{\sigma^{-1}}) > v(c)$. Therefore, $v(b - b^{\sigma^{-1}}) > v(c)$, in contradiction to the choice of $c$. ∎

PROPOSITION 12.4: *Let $\mathbf{K}$ be a Henselian field-valuation structure that satisfies the block approximation condition. Then $\mathbf{K}$ is unirationally closed.*

*Proof:* Consider a unirational arithmetical problem

$$\Phi = (V, X_i, L_i, \pi_i \colon U_i \to V \times_K L_i)_{i \in I_0}$$

for $\mathbf{K}$ as in Definition 6.1. Let $X' = \bigcup_{i \in I_0} X_i$. We find a solution $(\mathbf{a}, \mathbf{b}_x)_{x \in X'}$ of $\Phi$.

To this end consider $i \in I_0$. Since $\nu_{L_i}$ is continuous, $\mathbf{B}_i = \nu_{L_i}(X_i)$ is a closed subset of $\mathrm{Val}(L_i)$. For each $x \in X_i$ put $v_{x,i} = \nu_{L_i}(x)$. Then $(K, v_x|_K) \subseteq (L_i, v_{x,i}) \subseteq (K_x, v_x)$.

Since $U_i$ is birationally equivalent to $\mathbb{A}^r_{L_i}$, there exists $\mathbf{a}_i \in U_i(L_i)$. Then $\mathbf{b}_i = \pi_i(\mathbf{a}_i) \in V(L_i)$. By definition, $\pi_i$ is étale at $\mathbf{a}_i$ (see (3d) of Section 6). Thus, Corollary 10.4 (with $L_i$, $\pi_i\colon U_i \to V \times_K L_i$ replacing $K$, $\varphi\colon V \to W$) gives a partition $\mathbf{B}_i = \bigcup_{j \in J_i} \mathbf{B}_{ij}$ with $\mathbf{B}_{ij}$ closed in $\mathrm{Val}(L_i)$, an open neighborhood $\mathcal{U}_{ij}$ of $\mathbf{a}_i$ in $\mathrm{Set}(L_i, U_i, \mathbf{B}_{ij})$, and an open neighborhood $\mathcal{V}_{ij}$ of $\mathbf{b}_i$ in $\mathrm{Set}(L_i, V \times_K L_i, \mathbf{B}_{ij})$, $j \in J_i$ satisfying this:

(2) For all $j \in J_i$ and $x \in X_i$ with $v_{x,i} \in \mathbf{B}_{ij}$ the map $\pi_i\colon \mathcal{U}_{ij}(K_x) \to \mathcal{V}_{ij}(K_x)$ is a $v_x$-homeomorphism.

For all $i \in I_0$ and $j \in J_i$ Lemma 9.4 gives a partition $\mathbf{B}_{ij} = \bigcup_{l \in \Lambda_{ij}} \mathbf{B}_{ijl}$ with $\Lambda_{ijl}$ finite, $\mathbf{B}_{ijl}$ closed, and $c_{ijl} \in L_i^\times$, $l \in \Lambda_{ij}$, such that $\mathcal{B}_{\mathbf{b}_i, c_{ijl}, \mathbf{B}_{ijl}}(M, w) \subseteq \mathcal{V}_{ij}(M, w)$ for each $(M, w) \in \mathrm{Hensel}(L_i, \mathbf{B}_{ijl})$, $l \in \Lambda_{ij}$. For all $l \in \Lambda_{ij}$ put $L_{ijl} = L_i$, $X_{ijl} = \nu_{L_i}^{-1}(\mathbf{B}_{ijl})$, and $\mathbf{b}_{ijl} = \mathbf{b}_i$. Then $X_{ijl}$ is a closed subset of $X$, $X_i = \bigcup_{j \in J_i} \bigcup_{l \in \Lambda_{ij}} X_{ijl}$, and

(3) $$\{\mathbf{b} \in V(K_x) \mid v_x(\mathbf{b} - \mathbf{b}_i) > v_x(c_{ijl})\} \subseteq \mathcal{V}_{ij}(K_x)$$

for all $x \in X_{ijl}$ and $l \in \Lambda_{ij}$.

Since $X_i$ is open-closed in $X$, so are $X_{ijl}$. If $\sigma \in \mathrm{Gal}(L_i)$, then $X_{ijl}^\sigma = X_{ijl}$. Indeed, let $x \in X_{ijl}$. Then $\nu_{L_i}(x^\sigma) = v_{x^\sigma}|_{L_i} = v_x^\sigma|_{L_i} = v_x|_{L_i} \in \mathbf{B}_{ijl}$, so $x^\sigma \in X_{ijl}$. If $\sigma \in \mathrm{Gal}(K)$, $i, i' \in I_0$, $j \in J_i$, $j' \in J_{i'}$, $l \in \Lambda_{ij}$, $l' \in \Lambda_{i'j'}$, and $X_{ijl}^\sigma \cap X_{i'j'l'} \neq \emptyset$, then $X_i^\sigma \cap X_{i'} \neq \emptyset$, so $i' = i$ and $\sigma \in \mathrm{Gal}(L_i)$. Thus our assumption becomes $X_{ijl} \cap X_{ij'l'} \neq \emptyset$. Therefore, $j = j'$ and $l = l'$. It follows that

$$(V, X_{ijl}, L_{ijl}, \mathbf{a}_{ijl}, c_{ijl})_{i \in I_0, \, j \in J_i, \, l \in \Lambda_{ij}}$$

is a block approximation problem for $\mathbf{K}$.

The block approximation condition gives $\mathbf{b} \in V(K)$ with $v_x(\mathbf{b} - \mathbf{a}_i) > v_x(c_{ijl})$ for all $i \in I_0$, $j \in J_i$, $l \in \Lambda_{ij}$, and $x \in X_{ijl}$. By (3), $\mathbf{b} \in \mathcal{V}_{ij}(K_x)$. By (2), there is $\mathbf{a}_x \in \mathcal{U}_{ij}(K_x)$ with $\pi_i(\mathbf{a}_x) = \mathbf{b}$. In particular, $\mathbf{a}_x \in U_i(K_x)$. Thus, $(\mathbf{b}, \mathbf{a}_x)_{x \in X'}$ is a solution of $\Phi$. ∎

THEOREM 12.5: *Let* $\mathbf{K} = (K, X, K_x, v_x)_{x \in X}$ *be a proper Henselian field-valuation structure. Suppose* $\mathbf{K}$ *satisfies the block approximation condition. Then* $\mathrm{Gal}(\mathbf{K})$ *is a projective group structure.*

Proof: By Proposition 12.4, $\mathbf{K}$ is unirationally closed. Since $\mathrm{Gal}(\mathbf{K})$ is a proper group structure, $S_x = \mathrm{Gal}(K_x)$ for each $x \in X$ (Remark 2.1). Hence, by Proposition 6.4, $\mathrm{Gal}(\mathbf{K})$ is projective. ∎

This completes the proof of Part (a) of the Main Theorem. The rest of the work is devoted to the proof of Part (b) of the Main Theorem.

## 13. Rigid Henselian Extensions

This section is a continuation of Section 7. It contains various results about valued fields which are needed in the proof of Part (b) of the Main Theorem.

For a field extension $F/K$ let $\mathrm{Val}(F/K)$ be the space of all valuations of $F$ (including the trivial one) which are trivial on $K$. Denote the valuation ring of a valuation $w$ of $K$ by $O_w$, its maximal ideal by $M_w$, and its residue field by $\bar{K}_w$. Another valuation $v$ of $K$ is said to be **finer than** $w$ if $O_v \subseteq O_w$, equivalently if $M_w \subseteq M_v$. Thus, $w(x) < w(y)$ implies $v(x) < v(y)$, for all $x, y \in K$. Then $E = \bar{K}_w$ has a unique valuation $\bar{v}$ satisfying $\bar{v}(x + M_w) = v(x)$ for $x \in O_w$. In particular, $\bar{K}_v = \bar{E}_{\bar{v}}$. Denote $\bar{v}$ by $v/w$.

Conversely, given a valuation $\bar{v}$ of $E$, there is a unique valuation $v$ of $K$ which is finer than $w$ for which $v/w = \bar{v}$ [Jar, §3]. Then the place $\varphi_v \colon K \to \bar{K}_v \cup \{\infty\}$ corresponding to $v$ is the compositum of the place $\varphi_w \colon K \to \bar{K}_w \cup \{\infty\}$ and the place $\varphi_{\bar{v}} \colon \bar{K}_w \to \bar{K}_v \cup \{\infty\}$. We write $v = \bar{v} \cdot w$.

LEMMA 13.1: *Let $K$ be a field, $\tilde{K}$ its algebraic closure, and $T$ a set of indeterminates with $\mathrm{card}(T) \geq \mathrm{card}(\tilde{K})$. Put $F = K(T)$. Then, for each algebraic extension $L$ of $K$ there exists $v \in \mathrm{Val}(F/K)$ with $\bar{F}_v = L$.*

*Proof:* Put $m = \mathrm{card}(T)$. Choose a well ordered transfinite sequence $(a_\alpha)_{\alpha < m}$ which generates $L$ over $K$. Well-order $T$ as $(t_\alpha)_{\alpha < m}$. For each $\beta \leq m$ let $F_\beta = K(t_\alpha \mid \alpha \leq \beta)$ and $L_\beta = L(a_\alpha \mid \alpha \leq \beta)$.

Consider $\gamma \leq m$. Inductively suppose for each $\beta < \gamma$ there is a $v_\beta \in \mathrm{Val}(F_\beta/K)$ with $\bar{F}_\beta = L_\beta$ such that $v_{\beta'}$ extends $v_\beta$ whenever $\beta \leq \beta'$.

If $\gamma$ is a limit cardinal, then the union of all $v_\beta$ is a valuation $v_\gamma$ of $F_\gamma$ with residue field $L_\gamma$. Otherwise, $\gamma = \beta + 1$, $F_\gamma = F_\beta(t_\gamma)$, and $t_\gamma$ is transcendental over $F_\beta$. Extend $v_\beta$ to a place $v'$ of $F_\gamma$ with residue field $L_\beta(t_\gamma)$ with $t_\gamma$ being its own residue [Bou, Chap. VI, §10.1, Lemma 1, p. 434]. Let $w$ be the $L_\beta$-valuation of $L_\beta(t_\gamma)$ with $\bar{t}_\gamma = a_\gamma$ and $\overline{L_\beta(t_\gamma)} = L_\beta(a_\gamma) = L_\gamma$. Then $\varphi_w \circ \varphi_{v'}$. Hence, $v_\gamma = w \cdot v'$ extends $v_\beta$ and has $L_\gamma$ as residue field. This completes the induction.

The valuation $v = v_m$ of $F$ is trivial on $K$ and satisfies $\bar{F}_v = L$. ∎

LEMMA 13.2: *Consider a perfect field $K$.*

(a) *Let $L$ be an extension of $K$ and $v \in \mathrm{Val}(L/K)$. Suppose $(L,v)$ is Henselian, $\bar{L}_v$ is an algebraic extension of $K$, and res: $\mathrm{Gal}(L) \to \mathrm{Gal}(K)$ is an isomorphism. Then, $\bar{L}_v = K$.*

(b) *Let $L$ be a rigid Henselian extension of a $K$ and $L'$ a separable algebraic extension of $L$. Then $L'$ is a rigid Henselian extension of $L' \cap \tilde{K}$.*

(c) *Suppose $L/K$ and $M/L$ are rigid Henselian extensions. Then so is $M/K$.*

(d) *Let $K$ be a field and $I$ a totally ordered set. For each $i \in I$ let $(L_i, v_i)$ be a rigid Henselian extension of $K$. Suppose $(L_i, v_i) \subseteq (L_j, v_j)$ if $i \le j$. Put $(L,v) = \bigcup_{i \in I}(L_i, v_i)$. Then $(L,v)$ is a rigid Henselian extension of $K$.*

*Proof of (a):* By Lemma 7.4(a), reduction modulo $v$ defines an epimorphism $\rho$: $\mathrm{Gal}(L) \to \mathrm{Gal}(\bar{L}_v)$ and $\rho = \mathrm{res}_{L_s/\tilde{K}}$. Hence,

$$\mathrm{Gal}(\bar{L}_v) = \rho(\mathrm{Gal}(L)) = \mathrm{res}_{L_s/\tilde{K}}(\mathrm{Gal}(L)) = \mathrm{Gal}(K).$$

Therefore, $K = \bar{L}_v$.

*Proof of (b):* By definition, $L$ has a valuation $v$ such that $(L,v)$ is Henselian, $\bar{L}_v = K$, and res: $\mathrm{Gal}(L) \to \mathrm{Gal}(K)$ is an isomorphism. Denote the unique extension of $v$ to $L'$ by $v$. Then, $(L',v)$ is Henselian, $\overline{L'}_v/K$ is algebraic, and res: $\mathrm{Gal}(L') \to \mathrm{Gal}(L' \cap \tilde{K})$ is an isomorphism. By (a), $\overline{L'}_v = L' \cap \tilde{K}$. Therefore, $(L',v)$ is a rigid Henselian extension of $L' \cap \tilde{K}$.

*Proof of (c):* By assumption, $L$ admits a valuation $v$ and $M$ admits a valuation $w$ such that $(L,v)$ is a rigid Henselian extension of $K$ and $(M,w)$ is a rigid Henselian extension of $L$. Let $w' = v \cdot w$. Then $(M,w')$ is Henselian and $\bar{M}_{w'} = K$ [Jar, Prop. 13.1]. Also, $\varphi_{w'}(a) = \varphi_v(\varphi_w(a)) = a$ for each $a \in K$. Hence, $w'$ is trivial on $K$. Finally, res: $\mathrm{Gal}(M) \to \mathrm{Gal}(L)$ and res: $\mathrm{Gal}(L) \to \mathrm{Gal}(K)$ are isomorphisms. Therefore, res: $\mathrm{Gal}(M) \to \mathrm{Gal}(K)$ is an isomorphism. Consequently, $(M,w')$ is a rigid Henselian extension of $K$.

*Proof of (d):* Routine check.  ∎

An earlier version of the following result appears on page 24 of [Pop] without a proof.

LEMMA 13.3: *Let $K$ be a field and $T$ a set of indeterminates with $\mathrm{card}(T) \geq \mathrm{card}(\tilde{K})$. Put $F = K(T)_{\mathrm{ins}}$. Then, for each perfect algebraic extension $L$ of $K$ there are $v \in \mathrm{Val}(F/K)$ and a Henselian closure $(F_v, v)$ of $(F, v)$ which is a rigid Henselian extension of $L$.*

*Proof:* We may replace $K$ by $K_{\mathrm{ins}}$, if necessary, to assume $K$ is perfect. Write $T = \bigcup_{i=1}^{\infty} T_i$ with $\mathrm{card}(T_i) = \mathrm{card}(T)$ for each $i$. Inductively define $K_0 = K$ and $K_i = K_{i-1}(T_i)_{\mathrm{ins}}$ for $i = 1, 2, 3 \ldots$. Then $K_i$ is perfect and $\mathrm{card}(\tilde{K}_i) = \mathrm{card}(T_i)$ $i = 1, 2, 3, \ldots$. Also, $F = \bigcup_{i=1}^{\infty} K_i$.

Let $v_0$ be the trivial valuation of $K$. Put $K'_0 = K_0$ and $L_0 = L$. Suppose by induction we have constructed algebraic extensions $K'_i \subseteq L_i$ of $K_i$ and a valuation $v_i$ of $L_i$ satisfying this:

(1a) $(K'_i, v_i)$ is a Henselian closure of $(K_i, v_i|_{K_i})$.

(1b) $(L_i, v_i)$ is a rigid Henselian extension of $L$.

(1c) $L_{i-1} \subseteq K'_i$.

(1d) $v_i$ extends $v_{i-1}$.

Lemma 13.1 gives a valuation $w \in \mathrm{Val}(K_{i+1}/K_i)$ with residue field $L_i$. Let $v_{i+1} = v_i \cdot w$. Since $L_i/K_i$ is separable, $(K_{i+1}, w)$ has a Henselian closure $E$ which contains $L_i$. Since $v_{i+1}$ is finer than $w$, there is a Henselian closure $(K'_{i+1}, v_{i+1})$ of $(K_{i+1}, v_{i+1})$ which contains $E$ [Jar, Cor. 14.4], hence $L_i$. By Proposition 7.4(c), $K'_{i+1}$ has an algebraic extension $L_{i+1}$ such that res: $\mathrm{Gal}(L_{i+1}) \to \mathrm{Gal}(L)$ is an isomorphism. Denote the unique extension of $v_{i+1}$ to $L_{i+1}$ again by $v_{i+1}$. Then $(L_{i+1}, v_{i+1})$ is a rigid Henselian extension of $L$ (Lemma 13.2(b)).

Let $F_v = \bigcup_{i=1}^{\infty} K'_i$, $L_\infty = \bigcup_{i=1}^{\infty} L_i$, and $v = \bigcup_{i=1}^{\infty} v_i$. Then $v$ is a valuation of $F_v$ over $K$, $(F_v, v)$ is a Henselian closure of $(F, v)$, $F_v = L_\infty$, $L \subseteq L_\infty$, and res: $\mathrm{Gal}(L_\infty) \to \mathrm{Gal}(L)$ is an isomorphism. Thus, $(F_v, v)$ is a rigid Henselian extension of $L$. ∎

LEMMA 13.4: *Let $(K, v)$ be a valued field and $(E, w)$ a Henselian closure. Suppose $E \neq K_s$ and for each separable algebraic extension $F \neq K_s$ of $E$ the residue field $\bar{F}$ of $F$ under the unique extension of $w$ to $F$ is not separably closed. Then $\mathrm{Aut}(E/K) = 1$ and $EE^\sigma = K_s$ for each $\sigma \in \mathrm{Gal}(K) \smallsetminus \mathrm{Gal}(E)$.*

*Proof:* By assumption, $\bar{E}$ is not separably closed. Hence, by F. K. Schmidt - Engler, $\mathrm{Aut}(E/K) = 1$ [Jar, Prop. 14.5].

Consider now $\sigma \in \mathrm{Gal}(K)$. Put $E' = E^\sigma$ and $w' = w^\sigma$. Then $(E', w')$ is also a Henselian closure of $(K, v)$. Let $F = EE'$. Denote the unique extension of $w$ (resp. $w'$) to $F$ by $w_F$ (resp. $w'_F$). Then both $w_F$ and $w'_F$ extend $v$. We prove: Either $\sigma \in \mathrm{Gal}(E)$ or $F = K_s$.

CASE A: $w_F = w'_F$. Denote the unique extension of $w_F$ to $K_s$ by $w_s$. It coincide with the unique extension $w'_s$ of $w'_F$ to $K_s$. In addition, $w_s^\sigma$ is the unique extension of $w'$ to $K_s$, so also the unique extension of $w'_F$ to $K_s$. Thus, $w_s = w'_s = w_s^\sigma$. Therefore, $\sigma$ belongs to the decomposition group of $w_s$ over $E$, which is $\mathrm{Gal}(E)$.

CASE B: $w_F \neq w'_F$. By Engler, $w_F$ and $w'_F$ are incomparable [Jar, Prop. 6.6]. Since $F$ is Henselian with respect to both $w_F$ and $w'_F$, the field $\bar{F}_{w_F}$ is separably closed [Jar, Prop. 13.4]. Hence, by assumption, $F = K_s$. ∎

## 14. Projective Group Structures as Absolute Galois Structures

Part (b) of the Main Theorem gives for each proper projective group structure $\mathbf{G} = (G, X, G_x)_{x \in X}$ a proper field-valuation structure $\mathbf{L}$ and isomorphism $\lambda\colon \mathbf{G} \to \mathrm{Gal}(\mathbf{L})$. We call $\lambda$ a **Galois isomorphism of G**. An obvious necessary condition for the existence of a Galois isomorphism of $\mathbf{G}$ is the existence of a **Galois approximation** of $\mathbf{G}$. This is a rigid epimorphism $\kappa\colon \mathbf{G} \to \mathrm{Gal}(\mathbf{K})$ where $\mathbf{K}$ is a field structure. In this section we generalize [Pop, Thm. 3.4] and "lift" each Galois approximation of $\mathbf{G}$ to an isomorphism $\kappa'\colon \mathbf{G} \to \mathrm{Gal}(\mathbf{K}')$ where $\mathbf{K}'$ is a field structure. Then, in Section 14, we lift $\kappa'$ further to a Galois isomorphism $\lambda$ as above.

Here $\kappa'\colon \mathbf{G} \to \mathrm{Gal}(\mathbf{K}')$ is said to **lift** $\kappa$ if $K \subseteq K'$ and res: $\mathrm{Gal}(K') \to \mathrm{Gal}(K)$ extends to a rigid epimorphism $\rho\colon \mathrm{Gal}(\mathbf{K}') \to \mathrm{Gal}(\mathbf{K})$ with res $\circ\, \kappa' = \kappa$.

LEMMA 14.1: *Let $G$ be a profinite group, $H$ an open subgroup, $K$ a closed normal subgroup, and $\mathcal{G}$ a étale compact subset of $\mathrm{Subgr}(G)$. Suppose $\Gamma \cap K = 1$ for each $\Gamma \in \mathcal{G}$. Then $G$ has an open normal subgroup $N$ with $N \leq H$ and $\Gamma N \cap K N = N$ for each $\Gamma \in \mathcal{G}$.*

*Proof:* Let $\mathcal{N}$ be the set of open normal subgroups of $G$ containing $K$. Assume without loss $H \triangleleft G$. Now consider $\Delta \in \mathcal{G}$. Assume, for each $M \in \mathcal{N}$, the closed subset $\Delta \cap M \smallsetminus H$ of $G$ is nonempty. Then, by compactness of $G$, $\bigcap_{M \in \mathcal{N}} \Delta \cap M \smallsetminus H \neq \emptyset$. On the other hand, $\bigcap_{M \in \mathcal{N}} \Delta \cap M = \Delta \cap \bigcap_{M \in \mathcal{N}} M = \Delta \cap K = 1$. This contradiction gives $M_\Delta \in \mathcal{N}$ with $\Delta \cap M_\Delta \smallsetminus H = \emptyset$. In other words, $\Delta \cap M_\Delta \leq H$. It follows that $\Delta(H \cap M_\Delta) \cap M_\Delta \leq H$.

Now consider the étale open neighborhood $\mathcal{U}_\Delta = \mathrm{Subgr}\big(\Delta(H \cap M_\Delta)\big) \cap \mathcal{G}$ of $\Delta$ in $\mathcal{G}$. For each $\Gamma \in \mathcal{U}_\Delta$ we have $\Gamma \cap M_\Delta \leq \Delta(H \cap M_\Delta) \cap M_\Delta \leq H$.

Since $\mathcal{G}$ is étale compact, there are $\Delta_1, \ldots, \Delta_r \in \mathcal{G}$ with $\mathcal{G} = \bigcup_{i=1}^{r} \mathcal{U}_{\Delta_i}$. Then $N = H \cap \bigcap_{i=1}^{r} M_{\Delta_i}$ is the desired open normal subgroup of $G$. Indeed, let $\Gamma \in \mathcal{G}$. Then $\Gamma \in \mathcal{U}_{\Delta_j}$ for some $j$. So, $\Gamma \cap KN \leq \Gamma \cap M_{\Delta_j} \leq H$. Thus, $\Gamma \cap KN \leq H \cap \bigcap_{i=1}^{r} M_{\Delta_i} = N$. Therefore, $\Gamma N \cap K N = N$. ∎

In the following Lemma and its applications we use the relation $A \subset B$ between sets to mean "$A$ is a proper subset of $B$".

LEMMA 14.2: *Let $\mathbf{G} = (G, X, G_x)_{x \in X}$ be a proper projective group structure, $\kappa \colon \mathbf{G} \to$ Gal($\mathbf{K}$) a Galois approximation, and $G_0$ an open subgroup of $G$. Then $\kappa$ can be lifted to a Galois approximation $\varepsilon \colon \mathbf{G} \to$ Gal($\mathbf{E}'$) with $K \subset E'$, Ker($\varepsilon$) $\leq G_0$, and trans.deg($E'/K$) $< \infty$.*

*Proof:* Replacing $\mathbf{K}$ by $\mathbf{K}_{\mathrm{ins}}$ (Lemma 11.1), if necessary, we may assume $K$ is perfect. The rest of the proof has three parts.

PART A: *Replace* Gal($K$) *by a relative Galois group.* By definition, $G_x \cap \mathrm{Ker}(\kappa) = 1$ for each $x \in X$. Hence, Lemma 14.1 gives an open normal subgroup $N$ of $G$ contained in $G_0$ with

(1) $$G_x N \cap \mathrm{Ker}(\kappa)N = N \qquad \text{for each } x \in X.$$

Put $B = G/N$, $A = \mathrm{Gal}(K)/\kappa(N)$, let $\beta \colon G \to B$ and $\iota \colon \mathrm{Gal}(K) \to A$ be the quotient maps, and $\alpha \colon B \to A$ the epimorphism induced by $\kappa$. Then $\alpha \circ \beta = \iota \circ \kappa$. Let $\bar{G} = B \times_A \mathrm{Gal}(K)$. Then let $\bar{\kappa} \colon \bar{G} \to \mathrm{Gal}(K)$ and $\bar{\beta} \colon \bar{G} \to B$ be the coordinate projections. There is a unique morphism $\rho \colon G \to \bar{G}$ with $\bar{\kappa} \circ \rho = \kappa$ and $\bar{\beta} \circ \rho = \beta$.

Since $\mathrm{Ker}(\iota \circ \kappa) = N\mathrm{Ker}(\kappa) = \mathrm{Ker}(\beta)\mathrm{Ker}(\kappa)$, we may assume that $\bar{G} = G/N \cap \mathrm{Ker}(\kappa)$ and $\rho$ is the quotient map [FrJ, Section 20.2].

Let $L$ be the fixed field of $\kappa(N)$ in $\tilde{K}$. Identify $A$ with $\mathrm{Gal}(L/K)$ and $\iota$ with $\mathrm{res}_{\tilde{K}/L}$. Lemma 6.2 gives a regular extension $E$ of $K$ of transcendence degree equal to $|B|$ (in particular, $E \neq K$) and a finite Galois extension $F$ of $E$ containing $L$ with $B = \mathrm{Gal}(F/E)$ and $\alpha = \mathrm{res}_{F/L}$. Since $E/K$ is regular, $\bar{G} = \mathrm{Gal}(F/E) \times_{\mathrm{Gal}(L/K)} \mathrm{Gal}(K) = \mathrm{Gal}(F\tilde{K}/E)$, $\bar{\beta} = \mathrm{res}_{F\tilde{K}/F}$, and $\bar{\kappa} = \mathrm{res}_{F\tilde{K}/\tilde{K}}$.

Extend $\bar{G}$ to a group structure $\bar{\mathbf{G}} = \mathbf{G}/\mathrm{Ker}(\rho)$ and $\rho$ to the quotient map $\rho \colon \mathbf{G} \to \bar{\mathbf{G}}$. Then $\bar{\kappa}$ extends to a rigid epimorphism $\bar{\kappa} \colon \bar{\mathbf{G}} \to \mathrm{Gal}(\mathbf{K})$ such that $\kappa = \bar{\kappa} \circ \rho$.

PART B: *The cover* $\pi\colon \operatorname{Gal}(\mathbf{E}) \to \bar{\mathbf{G}}$. Write $\bar{\mathbf{G}}$ as $(\bar{G}, Y, \bar{G}_y)_{y\in Y}$. Put $\bar{N} = \rho(N)$. For each $y \in Y$ choose $x \in X$ such that $\rho(x) = y$. Then $\bar{G}_y = \rho(G_x)$ and $\bar{G}_y \bar{N} = \rho(G_x N)$ is an open subgroup of $\bar{G}$ which contains $\bar{G}_y$. Let $L_y$ be the fixed field of $\bar{\kappa}(\bar{G}_y\bar{N}) = \kappa(G_x N)$ in $\tilde{K}$ and $F_y$ the fixed field of $\bar{\beta}(\bar{G}_y\bar{N}) = \beta(G_x N) = G_x N/N$ in $F$. Then $\kappa(G_x N) = \operatorname{Gal}(L_y)$, $\beta(G_x N) = \operatorname{Gal}(F/F_y)$, and $\bar{G}_y\bar{N} = \operatorname{Gal}(F\tilde{K}/F_y)$. Since $\alpha = \operatorname{res}_{F/L}$ maps $\operatorname{Gal}(F/F_y)$ onto $\operatorname{res}_{\tilde{K}/L}(\operatorname{Gal}(L_y)) = \operatorname{Gal}(L/L_y)$, we have $L_y \subseteq F_y$. Also, $\operatorname{Ker}(\alpha) = \operatorname{Ker}(\kappa)N/N$. Hence, by (1), $\alpha$ is injective on $G_x N/N$. Thus $\alpha$ maps $\operatorname{Gal}(F/F_y)$ isomorphically onto $\operatorname{Gal}(L/L_y)$. By Lemma 6.2, $F_y/L_y$ is a purely transcendental extension.

Proposition 7.7(c) gives a perfect algebraic extension $E_y$ of $F_y$ which is a rigid Henselian extension of $L_y$. In particular, $\operatorname{res}_{\tilde{E}/\tilde{K}}\colon \operatorname{Gal}(E_y) \to \operatorname{Gal}(L_y)$ is an isomorphism. Therefore, $\tilde{E} = E_y\tilde{K}$ and $F\tilde{K} = F_yL\tilde{K} = F_y\tilde{K}$. Consequently, res: $\operatorname{Gal}(E_y) \to \operatorname{Gal}(F\tilde{K}/F_y)$ is an isomorphism.

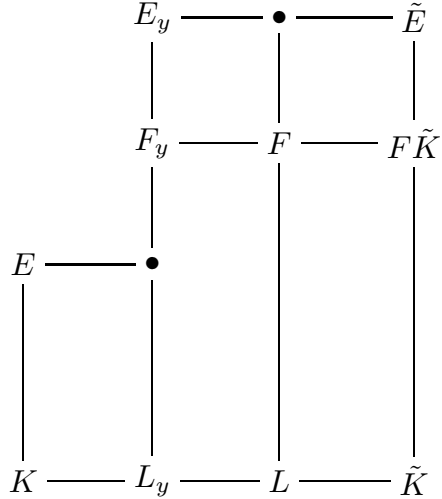

Lemma 3.6 gives a finite subset $\{y_i \mid i \in I_0\}$ of $Y$ and a special partition $(\bar{G}_i, Y_i, R_i)_{i\in I_0}$ of $\bar{\mathbf{G}}$ (Definition 3.5) such that $\bar{G}_i = \bar{G}_{y_i}\bar{N}$ and $y_i \in Y_i$ for each $i \in I_0$. Thus, $\operatorname{res}_{\tilde{E}/F\tilde{K}}\colon \operatorname{Gal}(E_{y_i} \to \bar{G}_i$ is an isomorphism, $i \in I_0$. Therefore, Lemma 5.1 extends $\operatorname{res}_{\tilde{E}/F\tilde{K}}\colon \operatorname{Gal}(E) \to \bar{G}$ to a cover of group structures. This means there is a field structure $\mathbf{E}$ on $E$ and a cover $\pi\colon \operatorname{Gal}(\mathbf{E}) \to \bar{\mathbf{G}}$.

PART C: *Applying projectivity.* Since $\mathbf{G}$ is projective, there is a morphism $\varepsilon\colon \mathbf{G} \to \operatorname{Gal}(\mathbf{E})$ with $\pi \circ \varepsilon = \rho$. Replace $\operatorname{Gal}(\mathbf{E})$ by $\varepsilon(\mathbf{G})$ and replace $\mathbf{E}$ by an appropriate sub-field-structure $\mathbf{E}'$ to assume that the underlying maps of $\varepsilon$ are surjective. Since both $\rho$ and $\pi$ are covers, $\rho\colon G_x \to \bar{G}_{\rho(x)}$ and $\pi\colon \mathrm{Gal}(E'_{\varepsilon(x)}) \to \bar{G}_{\rho(x)}$ are isomorphisms, so $\varepsilon\colon G_x \to \mathrm{Gal}(E'_{\varepsilon(x)})$ is an isomorphism for each $x \in X$. Thus $\varepsilon$ is a rigid epimorphism, hence $\varepsilon$ is a Galois approximation of $\mathbf{G}$ which lifts $\kappa$.

Since $\mathrm{Ker}(\rho) \leq N \leq G_0$, also $\mathrm{Ker}(\varepsilon) \leq G_0$. Finally, since $E$ is a proper extension of $K$, so is $E'$. ∎

PROPOSITION 14.4: *Let $\mathbf{G}$ be a proper projective group structure and $\kappa\colon \mathbf{G} \to \mathrm{Gal}(\mathbf{K})$ a Galois approximation. Then $\kappa$ can be lifted to a Galois isomorphism $\lambda\colon \mathbf{G} \to \mathrm{Gal}(\mathbf{L})$ with an underlying perfect field.*

*Proof:* Let $\{G_\alpha \mid \alpha < m\}$ be a well ordering of all open subgroups of $G$. By transfinite induction we construct for each $\alpha \leq m$ a Galois approximation $\kappa_\alpha\colon \mathbf{G} \to \mathrm{Gal}(\mathbf{K}_\alpha)$ such that $\kappa_0 = \kappa$, $\kappa_\beta$ lifts $\kappa_\alpha$ if $\alpha \leq \beta \leq m$, the underlying field $K_\alpha$ of $\mathbf{K}_\alpha$ is perfect, and $\mathrm{Ker}(\kappa_{\alpha+1}) \leq G_\alpha$.

Indeed, suppose $\beta$ is an ordinal number at most $m$ and $\kappa_\alpha$ have already been constructed for each $\alpha < \beta$. If $\beta = \alpha + 1$ is a successor ordinal, use Lemma 14.2 to construct a Galois approximation $\kappa_\beta\colon \mathbf{G} \to \mathrm{Gal}(\mathbf{K}_\beta)$ and a rigid epimorphism $\rho_{\beta,\alpha}\colon \mathrm{Gal}(\mathbf{K}_\beta) \to \mathrm{Gal}(\mathbf{K}_\alpha)$ with $\rho_{\beta,\alpha} \circ \kappa_\beta = \kappa_\alpha$ such that $K_\beta$ is perfect, $K_\alpha \subseteq K_\beta$, $\rho_{\beta,\alpha}\colon \mathrm{Gal}(K_\beta) \to \mathrm{Gal}(K_\alpha)$ is the restriction map, and $\mathrm{Ker}(\kappa_\beta) \leq G_\alpha$. If $\beta$ is a limit ordinal, then $\{\mathrm{Gal}(\mathbf{K}_\alpha), \rho_{\alpha',\alpha} \mid \alpha \leq \alpha' < \beta\}$ is an inverse system of Galois group structures with $K_\alpha \subseteq K_{\alpha'}$, $\rho_{\alpha,\alpha'}\colon \mathrm{Gal}(K'_\alpha) \to \mathrm{Gal}(K_\alpha)$ are the restriction maps, and $\rho_{\alpha,\alpha'}\colon \mathrm{Gal}(\mathbf{K}'_\alpha) \to \mathrm{Gal}(\mathbf{K}_\alpha)$ are rigid epimorphisms. Then $\mathrm{Gal}(\mathbf{K}_\beta) = \varprojlim \mathrm{Gal}(\mathbf{K}_\alpha)$ is a group structure with $K_\beta = \bigcup_{\alpha < \beta} K_\alpha$ (Remark 2.7) and with rigid projections $\rho_{\beta,\alpha}\colon \mathrm{Gal}(\mathbf{K}_\beta) \to \mathrm{Gal}(\mathbf{K}_\alpha)$. Moreover, the inverse limit of the $\kappa_\alpha$'s gives a Galois approximation $\kappa_\beta\colon \mathbf{G} \to \mathrm{Gal}(\mathbf{K}_\beta)$ with $\rho_{\beta,\alpha} \circ \kappa_\beta = \kappa_\alpha$ for each $\alpha < \beta$.

Having completed the transfinite induction, we put $\mathbf{L} = \mathbf{K}_m$ and $\lambda = \kappa_m$. Then the underlying field of $\mathbf{L}$ is perfect and $\lambda\colon \mathbf{G} \to \mathrm{Gal}(\mathbf{L})$ is a Galois approximation lifting $\kappa$ (Remark 2.7). Moreover, $\mathrm{Ker}(\lambda) \le \bigcap_{\alpha < m} G_\alpha = 1$. Since $\mathbf{G}$ is proper, $\lambda$ is an isomorphism (Remark 2.1). ∎

*Remark 14.3: Cardinality of $L$.* We may assume that the cardinality of $L$ in Proposition 14.4 is not smaller than any given cardinality $m$. Indeed, without loss $\kappa$ is an isomorphism. Hence, if $\lambda$ lifts $\kappa$, then $\lambda$ is an isomorphism. Put $\lambda_0 = \kappa$. By transfinite induction construct a family of Galois approximations $\lambda_\alpha \colon \mathbf{G} \to \mathrm{Gal}(\mathbf{L}_\alpha)$ with underlying fields $L_\alpha$ such that $\lambda_\beta$ lifts $\lambda_\alpha$ and $L_\alpha \subset L_\beta$ for all $\alpha \le \beta \le m$. Namely, if $\beta$ is a limit ordinal, put $L_\beta = \bigcup_{\alpha < \beta} L_\beta$ and $L_{\beta, x} = \bigcup_{\alpha < \beta} L_{\alpha, x}$; otherwise use Lemma 14.2 to construct a lifting $\lambda_\beta$ of $\lambda_{\beta - 1}$. Then $\lambda = \lambda_m$ has the required property. ∎

## 15. From Field Structures to Field-Valuation Structures

Having lifted a given Galois approximation $\kappa\colon \mathbf{G} \to \mathrm{Gal}(\mathbf{K})$ of a proper projective group structure to a Galois isomorphism $\varepsilon\colon \mathbf{G} \to \mathrm{Gal}(\mathbf{E})$, we wish to extend $\mathbf{E}$ to a proper field-valuation structure $\mathbf{L}$ which satisfies the block approximation condition and res: $\mathrm{Gal}(\mathbf{L}) \to \mathrm{Gal}(\mathbf{K})$ is an isomorphism.

The crucial step in the construction is, starting from a field-valuation structure $\mathbf{K}$ and a data $(V, X_i, L_i, \mathbf{b}_i)_{i \in I_0}$ satisfying (2) below, to extend $\mathbf{K}$ to a field-valuation structure with a point $\mathbf{z} \in V(K')$ blockwise approximating each $\mathbf{b}_i$ infinitely well over $K$; that is, $\mathbf{z}$ satisfies Condition (3c) below.

Let $\mathbf{K} = (K, X, K_x, v_x)_{x \in X}$ and $\mathbf{K}' = (K', X', K'_x, v'_x)_{x \in X}$ be field-valuation structures. We say $\mathbf{K}'$ **extends** $\mathbf{K}$ and write $\mathbf{K} \subseteq \mathbf{K}'$ if $K \subseteq K'$, $K_x \subseteq K'_x$, and $v_x = v'_x|_{K_x}$ for each $x \in X$.

LEMMA 15.1: *Let* $\mathbf{K} = (K, X, K_x, v_x)_{x \in X}$ *and* $\bar{\mathbf{K}} = (\bar{K}, X, \bar{K}_x, \bar{v}_x)_{x \in X}$ *be proper Henselian field-valuation structures satisfying this:*

(1a) *$\bar{K}$ and $K$ are perfect.*

(1b) *$\bar{\mathbf{K}} \subseteq \mathbf{K}$ and the map $\mathrm{res}_{K_s/\bar{K}_s}\colon \mathrm{Gal}(\mathbf{K}) \to \mathrm{Gal}(\bar{\mathbf{K}})$ (with the identity map $X \to X$) is an isomorphism.*

(1c) *$\mathrm{Gal}(\bar{\mathbf{K}})$ is projective.*

(1d) *$\bar{v}_x$ is the trivial valuation of $\bar{K}_x$, $x \in X$.*

(1e) *$\bar{K}_x$ is the residue field of $(K_x, v_x)$, $x \in X$.*

*Consider a data $(V, X_i, L_i, \mathbf{b}_i)_{i \in I_0}$ satisfying this:*

(2a) *$(\mathrm{Gal}(L_i), X_i)_{i \in I_0}$ is a special partition of $\mathrm{Gal}(\mathbf{K})$.*

(2b) *$V$ is a smooth affine variety over $K$.*

(2c) *$\mathbf{b}_i \in V(L_i)$.*

*Then $\mathbf{K}$ has a proper field-valuation extension $\mathbf{K}' = (K', X, K'_x, v'_x)_{x \in X}$ with $K'$ perfect satisfying this:*

(3a) *$(K'_x, v'_x)$ is a Henselian field with residue field $\bar{K}_x$, $x \in X$.*

(3b) *$\mathrm{res}_{\bar{K}'/\bar{K}}\colon \mathrm{Gal}(K') \to \mathrm{Gal}(K)$ together with the identity map $X \to X$ form an isomorphism $\mathrm{Gal}(\mathbf{K}') \to \mathrm{Gal}(\mathbf{K})$.*

(3c) *There is* $\mathbf{z} \in V(K')$ *with* $v_x'(\mathbf{z} - \mathbf{b}_i) > v_x'(c)$ *for all* $i \in I_0$, $x \in X_i$, *and* $c \in K^\times$.

*Proof:* Suppose first $X = \{x\}$. By Remark 2.1, $K_x = K$. Let $(V, X_i, L_i, \mathbf{b}_i)_{i \in I_0}$ is a data satisfying (2). Then $I_0 = \{i\}$ and $L_i = K$. Hence, $\mathbf{K}' = \mathbf{K}$ and $\mathbf{z} = \mathbf{b}_i$ satisfy (3). We may therefore suppose $X$ has at least two elements.

We construct an extension $F$ of $K$ of large transcendence degree such that $V(F)$ contains a generic point $\mathbf{z}$ of $V$ over $K$. Then we extend $F$ to a proper field-valuation structure $\mathbf{F} = (F, Y, F_y, w_y)$ with a cover $\pi\colon \mathrm{Gal}(\mathbf{F}) \to \mathrm{Gal}(\mathbf{K})$ such that $(F_y, w_y)$ is a rigid Henselian extension of $(K_{\pi(y)}, v_{\pi(y)})$, $y \in Y$. Since $\mathrm{Gal}(\mathbf{K})$ is projective, $\mathbf{F}$ has a extension $\mathbf{F}' = (K', X', F_y, w_y)_{y \in X'}$ such that $\pi\colon \mathrm{Gal}(\mathbf{K}') \to \mathrm{Gal}(\mathbf{K})$ is an isomorphism. Renaming $X'$ as $X$ gives the desired extension $\mathbf{K}'$ of $\mathbf{K}$. In this construction, the valuations $w_y$ are defined in such a manner that $\mathbf{z}$ blockwise approximates the $\mathbf{b}_i$'s infinitely well over $K$. The construction has six parts.

PART A: *The field* $F$. Let $\mathbf{z}$ be a generic point of $V$ over $K$. Put $E = K(\mathbf{z})_{\mathrm{ins}}$. Since $V$ is absolutely irreducible and $K$ is perfect, $E/K$ is a regular extension. Hence, res: $\mathrm{Gal}(E) \to \mathrm{Gal}(K)$ is an epimorphism. Let $i \in I_0$. By [JaR, p. 456, Cor. A2], there is an $L_i$-place $\bar{\rho}_i\colon L_i(\mathbf{z}) \to L_i \cup \{\infty\}$ with $\bar{\rho}_i(\mathbf{z}) = \mathbf{b}_i$. By Proposition 7.4(c), there is a perfect algebraic extension $E_i$ of $L_i(\mathbf{z})$ and an extension of $\bar{\rho}_i$ to a rigid Henselian place $\rho_i\colon E_i \to L_i \cup \{\infty\}$. In particular, $E \subseteq E_i$ and $\rho_i(\mathbf{z}) = \mathbf{b}_i$.

Choose a set $T$ of indeterminates with $\mathrm{card}(T) \geq \mathrm{card}(E)$. Put $F = E(T)_{\mathrm{ins}}$. Then $F$ is a regular extension of $E$, hence of $K$. Therefore, res: $\mathrm{Gal}(F) \to \mathrm{Gal}(K)$ is an epimorphism. In addition, $\mathbf{z} \in V(F)$.

PART B: *The field structure* $(F, Y, F_y)_{y \in Y}$. Lemma 13.3 gives for each $i \in I_0$ a valuation $w_i'$ of $F$ with residue field $E_i$ and a Henselian closure $(F_i, w_i')$ of $(F, w_i')$ such that the corresponding place $\varphi_i\colon F_i \to E_i \cup \{\infty\}$ is rigid.

Put $\psi_i = \rho_i \circ \varphi_i$. Then $\psi_i\colon F_i \to L_i \cup \{\infty\}$ is a rigid $L_i$-place (Lemma 13.2(c)). In particular, res: $\mathrm{Gal}(F_i) \to \mathrm{Gal}(L_i)$ is an isomorphism. Moreover, $\psi_i$ extends to a $\tilde{K}$-place $\psi_i\colon \tilde{F} \to \tilde{K} \cup \{\infty\}$ with $\psi_i(F') = (F' \cap \tilde{K}) \cup \{\infty\}$ for each algebraic extension $F'$ of $F_i$. Denote the corresponding valuation by $w_i'$. Thus, if $F'$ is not algebraically closed, then the residue field of $F'$ with respect to $w_i'$ is not algebraically closed. By

Lemma 13.4,

$$(4) \qquad\qquad F_i F_i^\kappa = \tilde{F}$$

for each $\kappa \in \mathrm{Gal}(F) \smallsetminus \mathrm{Gal}(F_i)$.

By Lemma 5.1, $\mathrm{Gal}(F)$ extends to a proper group structure

$$(5) \qquad\qquad \mathrm{Gal}(\mathbf{F}) = (\mathrm{Gal}(F), Y, \mathrm{Gal}(F_y))_{y \in Y}$$

and res: $\mathrm{Gal}(F) \to \mathrm{Gal}(K)$ extends to a cover $\pi\colon \mathrm{Gal}(\mathbf{F}) \to \mathrm{Gal}(\mathbf{K})$ of group structures.

PART C: *The field-valuation structure* $\mathbf{F} = (F, Y, F_y, w_y)_{y \in Y}$. In addition to the cover $\pi$ mentioned in Part B, Lemma 5.1 gives for each $i \in I_0$ a subspace $Y_i$ of $Y$ such that $\pi(Y_i) = X_i$, $F_i \le F_y$ for each $y \in Y_i$, and $Y = \bigcup_{i \in I_0} Y_i^{\mathrm{Gal}(E)}$.

Consider $i \in I_0$ and $y \in Y_i$. Let $x = \pi(y)$. Then $F_i \le F_y$, $L_i \le K_x$, and res: $\mathrm{Gal}(F_y) \to \mathrm{Gal}(K_x)$ is an isomorphism (because $\pi$ is a cover). Thus, $F_y = F_i K_x$ and $K_x = F_y \cap \tilde{K}$. Since $\psi_i\colon F_i \to L_i \cup \{\infty\}$ is a rigid $L_i$-place (Part B), $\psi_i(F_y) = K_x \cup \{\infty\}$ (Lemma 13.2(b)). Since $\mathbf{K}$ is proper and $X$ has at least two elements, $K_x \ne \tilde{K}$ (Remark 2.1), so $F_y \ne \tilde{F}$.

By assumption, $(K_x, v_x)$ is Henselian. Hence, $v_x$ uniquely extends to a valuation $v_x$ of $\tilde{K}$. Let $w_y = v_x \cdot w_i'$ be the unique valuation of $\tilde{F}$ such that $w_y(u) = v_x(\psi_i(u))$ for each $u \in \tilde{F}$ with $\psi_i(u) \in \tilde{K}$. Then $O_{w_y} = \{u \in \tilde{F} \mid v_x(\psi_i(u)) \ge 0\}$, Thus, if $u \in \tilde{F}$ satisfies $\psi_i(u) = \infty$, then $w_y(u) < 0$. If $u \in \tilde{K}$, then $\psi_i(u) = u$, so $w_y(u) = v_x(u)$. Hence, $w_y$ extends $v_x$ (See also the beginning of Section 13.) Since $(K_x, v_x)$ and $(F_y, w_i')$ are Henselian, $(F_y, w_y)$ is Henselian [Jar, Prop. 13.1]. In addition, $\bar{K}_x$ is the residue field of $F_y$ at $w_y$.

We would like to define $(F_y, w_y)$ for all $y \in Y$. So, we consider $\sigma \in \mathrm{Gal}(F)$ and suppose, in addition to the assumption made above, that $y^\sigma \in Y_j$ for some $j \in I_0$. We prove that $w_{y^\sigma} = w_y^\sigma$.

Indeed, $\pi(y) \in X_i$ and $\pi(y)^{\pi(\sigma)} \in X_j$. Hence, $X_i^{\pi(\sigma)} \cap X_j \neq \emptyset$. By (2g) of Section 3, $i = j$ and $\pi(\sigma) \in \mathrm{Gal}(L_i)$. Hence, there are $\zeta \in \mathrm{Gal}(F_i)$ and $\kappa \in \mathrm{Gal}(F\tilde{K})$ with $\sigma = \kappa\zeta$. Since $y^\sigma \in Y_i$, we have $F_i \subseteq F_{y^\sigma} = F_y^\sigma$. Therefore, $F_i F_i^{\kappa^{-1}} = F_i F_i^{\zeta^{-1}\sigma^{-1}} = F_i F_i^{\sigma^{-1}} \subseteq F_y \subset \tilde{F}$. By (4), $\kappa = 1$, so $\sigma \in \mathrm{Gal}(F_i)$. Now consider $u \in F_y^\sigma$ with $\psi_i(u) \in K_x$. Since $\psi_i$ is rigid, $\psi_i(u^{\sigma^{-1}}) = \psi_i(u)^{\pi(\sigma)^{-1}}$ (Proposition 7.4(a)). Therefore, $w_y^\sigma(u) = w_y(u^{\sigma^{-1}}) = v_x(\psi_i(u^{\sigma^{-1}})) = v_x(\psi_i(u)^{\pi(\sigma)^{-1}}) = v_x^{\pi(\sigma)}(\psi_i(u)) = v_{x^{\pi(\sigma)}}(\psi_i(u)) = v_{\pi(y^\sigma)}(\psi_i(u)) = w_{y^\sigma}(u)$. It follows, $w_y^\sigma = w_{y^\sigma}$ on $F_{y^\sigma}$, and therefore also on $\tilde{F}$, as claimed.

For an arbitrary $y' \in Y$ there are $\tau \in \mathrm{Gal}(F)$, $i \in I_0$, and $y \in Y_i$ with $y' = y^\tau$. Since $(F, Y, F_y)_{y \in Y}$ is a field structure, $F_{y'} = F_y^\tau$. Define $w_{y'}$ to be $w_y^\tau$. By the preceding paragraph, this is a good definition. Thus, with $x' = \pi(y')$, the valued field $(F_{y'}, w_{y'})$ is a rigid Henselian extension of $(K_{x'}, v_{x'})$. Moreover, $w_{(y')^\sigma} = w_{y'}^\sigma$ for all $\sigma \in \mathrm{Gal}(F)$.

PART D: *Continuity of the maps $\nu_F \colon Y_F \to \mathrm{Val}(\tilde{F})$.* For each $x \in X$ let $\nu_K(x) = v_x$. Since $\mathbf{K}$ is a Henselian field-valuation structure, the map $\nu_K \colon X \to \mathrm{Val}(\tilde{K})$ is continuous (Lemma 11.2). Similarly, for each $y \in Y$ let $\nu_F(y) = w_y$. By Lemma 11.2, it suffices to prove that the map $\nu_F \colon Y \to \mathrm{Val}(\tilde{F})$ is continuous.

We start by proving that for each $i \in I_0$, the restriction of $\nu_F$ to $Y_i$ is continuous. Let $y \in Y_i$ and let $u \in \tilde{F}$ such that $w_y(u) > 0$. By Part C, $\psi_i(u) \neq \infty$ and $v_{\pi(y)}(\psi_i(u)) = w_y(y) > 0$. If $y' \in Y$ is sufficiently close to $y$, then $\pi(y')$ is sufficiently close to $\pi(y)$, and hence $v_{\pi(y')}(\psi_i(u)) > 0$ (because $\nu_K$ is continuous). Thus $w_{y'}(u) > 0$. Similarly one shows that if $w_y(u) \geq 0$ and $y'$ is sufficiently close to $y$, then $w_{y'}(u) \geq 0$.

It follows that the map $\nu_i \colon Y_i \times \mathrm{Gal}(F) \to \mathrm{Val}(\tilde{F})$ given by $\nu_i(y, \tau) = w_y^\tau$ is continuous. Indeed, let $a \in \tilde{F}$ and suppose $w_y^\tau(a) > 0$. If $y' \in Y_i$ is sufficiently close to $y$ and $\tau' \in \mathrm{Gal}(F)$ is sufficiently closed to $\tau$, then, by the preceding paragraph, $w_{y'}^{\tau'}(a) = w_{y'}(a^{(\tau')^{-1}}) = w_{y'}(a^{\tau^{-1}}) = w_y^\tau(a) > 0$. Similar statement holds for $\geq$ replacing $>$.

Let $\tilde{\nu}_i$ be the restriction of $\nu_{\tilde{\mathbb{F}}}$ to $Y_i^{\mathrm{Gal}(F)}$. Let $\mu\colon Y_i \times \mathrm{Gal}(F) \to Y_i^{\mathrm{Gal}(F)}$ be the map defined by $\mu(y,\tau) = y^\tau$. By Part C, $\nu_i = \tilde{\nu}_i \circ \mu$. Also, $\mu$ a continuous map between profinite spaces, hence closed. By the preceding paragraph, for each closed subset $C$ of $\mathrm{Val}(F)$, the set $\nu_i^{-1}(C)$ is closed in $Y_i \times \mathrm{Gal}(F)$. Therefore, $\tilde{\nu}_i^{-1}(C) = \mu(\nu_i^{-1}(C))$ is a closed subset of $Y_i^{\mathrm{Gal}(F)}$. Consequently, $\tilde{\nu}_i$ is continuous.

Since $Y = \bigcup_{i \in I_0} Y_i^{\mathrm{Gal}(F)}$, the preceding paragraph implies $\nu_F\colon Y \to \mathrm{Val}(\tilde{F})$ is continuous, as claimed.

PART E: *The proper group structure* $\mathbf{G}'$. By (1b) and (1c), $\mathrm{Gal}(\mathbf{K})$ is projective. By Part B, $\pi\colon \mathrm{Gal}(\mathbf{F}) \to \mathrm{Gal}(\mathbf{K})$ is a cover of group structures. Hence, by Corollary 4.3, $\mathrm{Gal}(\mathbf{F})$ has a proper sub-group-structure

$$\mathbf{G}' = (\mathrm{Gal}(K'), X', \mathrm{Gal}(F_{x'}))_{x' \in X'} \ ,$$

where $K'$ is an algebraic extension of $F$ and $X' \subseteq Y$ such that $\pi\colon \mathbf{G}' \to \mathrm{Gal}(\mathbf{K})$ is an isomorphism. In particular, $\mathrm{res}\colon \mathrm{Gal}(K') \to \mathrm{Gal}(K)$ is an isomorphism and $\pi\colon X' \to X$ is a homeomorphism. Then $\mathbf{F}' = (K', X', F_{x'}, w_{x'})_{x' \in X'}$ is a field-valuation structure.

PART F: *The proper field-valuation structure* $\mathbf{K}'$. For each $x \in X$ let $x'$ be the unique element of $X'$ with $\pi(x') = x$. Put $K_x' = F_{x'}$ and $v_x' = w_{x'}$. Then $\mathbf{K}' = (K', X, K_x', v_x')_{x \in X}$ is a proper field structure isomorphic to $\mathbf{F}'$. In addition, $\mathbf{K}'$ extends $\mathbf{K}$ and satisfies Conditions (3a) and (3b).

We still have to prove Condition (3c) (block approximation). Let $\mathbf{z} = (z_1, \ldots, z_n)$ and $\mathbf{b}_i = (b_{i1}, \ldots, b_{in})$, $i \in I_0$. Then $\psi_i(\mathbf{z}) = \rho_i(\mathbf{z}) = \mathbf{b}_i$. Let $y \in Y_i$ and put $x = \pi(y)$. Then, for all $c \in K^\times$ and $1 \leq j \leq n$ we have

$$w_y\Big(\frac{z_j - b_{ij}}{c}\Big) = v_x\Big(\frac{\psi_i(z_j) - b_{ij}}{c}\Big) = v_x\Big(\frac{0}{c}\Big) > 0.$$

Therefore, $w_y(\mathbf{z} - \mathbf{b}_i) > w_y(c)$.

Finally, consider $x \in X_i$. Choose $x' \in X'$ and $y \in Y_i$ with $\pi(x') = x = \pi(y)$. Then there is $\kappa \in \mathrm{Gal}(F\tilde{K})$ with $x' = y^\kappa$. By the preceding paragraph, $v_x'(\mathbf{z} - \mathbf{b}_i) =$ $w_{x'}(\mathbf{z} - \mathbf{b}_i) = w_y(\mathbf{z}^{\kappa^{-1}} - \mathbf{b}^{\kappa^{-1}}) = w_y(\mathbf{z} - \mathbf{b}_i) > w_y(c) = v_x(c) = v'_x(c)$ for all $c \in K^\times$. This concludes the proof of the Lemma. ∎

We apply Lemma 15.1 in each step of a transfinite induction. In the rest of this section we write res: $\mathrm{Gal}(\mathbf{L}) \to \mathrm{Gal}(\mathbf{K})$ for proper field structures $\mathbf{K} \subseteq \mathbf{L}$ to denote the unique morphism that extend the homomorphism res: $\mathrm{Gal}(L) \to \mathrm{Gal}(K)$ (Remark 2.1).

LEMMA 15.2: *Let* $\mathbf{K} = (K, X, K_x, v_x)_{x \in X}$ *and* $\bar{\mathbf{K}} = (\bar{K}, X, \bar{K}_x, \bar{v}_x)$ *be proper Henselian field-valuation structures satisfying (1). Then* $\mathbf{K}$ *has a proper field-valuation extension* $\mathbf{L} = (L, X, L_x, w_x)_{x \in X}$ *with* $L$ *perfect satisfying this:*

(7a) $(L_x, w_x)$ *is Henselian with residue field* $\bar{K}_x$, $x \in X$.

(7b) *res:* $\mathrm{Gal}(\mathbf{L}) \to \mathrm{Gal}(\mathbf{K})$ *is an isomorphism.*

(7c) $\mathbf{L}$ *satisfies the block approximation condition.*

*Proof:* Well-order all data satisfying (2) in a transfinite sequence

$$(V_\alpha, X_{\alpha,i}, K_{\alpha,i}, \mathbf{b}_{\alpha,i})_{i \in I_\alpha}, \quad \alpha < m.$$

Use transfinite induction and Lemma 15.1 to construct for each ordinal number $\alpha \leq m$ a proper field-valuation structure $\mathbf{K}_\alpha = (K_\alpha, X, K_{\alpha,x}, v_{\alpha,x})_{x \in X}$ with $K_\alpha$ perfect satisfying these conditions:

(8a) $(K_{\alpha,x}, v_{\alpha,x})$ is a Henselian field with residue field $\bar{K}_x$, $x \in X$.

(8b) $\mathbf{K}_\alpha \subseteq \mathbf{K}_\beta$ and res: $\mathrm{Gal}(\mathbf{K}_\beta) \to \mathrm{Gal}(\mathbf{K}_\alpha)$ is an isomorphism for all $\alpha < \beta \leq m$.

(8c) $\mathbf{K}_\beta = \bigcup_{\alpha < \beta} \mathbf{K}_\alpha$ for each limit ordinal $\beta \leq m$.

(8d) For each ordinal number $\alpha < m$ there is a point $\mathbf{z} \in V_\alpha(K_{\alpha+1})$ with $v_{\alpha+1,x}(\mathbf{z} - \mathbf{b}_{\alpha,i}) > v_{\alpha+1,x}(c)$ for all $i \in I_\alpha$, $x \in X_{\alpha,i}$, and $c \in K_\alpha^\times$.

Rewrite $\mathbf{K}_m$ as $\mathbf{L}_1 = (L_1, X, L_{1,x}, w_{1,x})_{x \in X}$. Then:

(9a) $(L_1, v_{1,x})$ is a Henselian field with residue field $\bar{K}_x$, $x \in X$.

(9b) $\mathbf{K} \subseteq \mathbf{L}_1$ and res: $\mathrm{Gal}(\mathbf{L}_1) \to \mathrm{Gal}(\mathbf{K})$ is an isomorphism.

(9c) Each approximation problem $(V, X_i, K_i, \mathbf{b}_i)_{i \in I_0}$ for $\mathbf{K}$ has a solution $\mathbf{z} \in V(L_1)$.

Finally use usual induction to construct an ascending sequence of proper field-valuation structures $\mathbf{L}_j$, $j = 1, 2, 3, \ldots$, such that $\mathbf{L}_{j+1}$ relates to $\mathbf{L}_j$ in the same way that $\mathbf{L}_1$ relates to $\mathbf{K}$, $j = 1, 2, 3, \ldots$. The structure $\mathbf{L} = \bigcup_{j=1}^\infty \mathbf{L}_j$ satisfies (7). ∎

PROPOSITION 15.3: *Let* $\mathbf{K} = (K, X, K_x)_{x \in X}$ *be a proper field structure with* $\mathrm{Gal}(\mathbf{K})$ *projective. Then there is a proper field-valuation structure* $\mathbf{L} = (L, X, L_x, w_x)_{x \in X}$ *with* $L$ *perfect having these properties:*

(10a) $(L, X, L_x)_{x \in X}$ *extends* $\mathbf{K}$.

(10b) $(L_x, w_x)$ *is Henselian with residue field* $(K_x)_{\mathrm{ins}}$, $x \in X$.

(10c) res: $\mathrm{Gal}(\mathbf{L}) \to \mathrm{Gal}(\mathbf{K})$ *is an isomorphism.*

(10d) $\mathbf{L}$ *satisfies the block approximation condition.*

*Proof:*  Replace $\mathbf{K}$ by $\mathbf{K}_{\mathrm{ins}}$, if necessary, to assume $K$ is perfect. Identify $\mathbf{K}$ with $(K, X, K_x, v_x)_{x \in X}$, where $v_x$ the trivial valuation on $K_x$ for each $x \in X$. Put $\bar{K}_x = K_x$, $\bar{v}_x = v_x$, $\bar{K} = K$, and $\bar{\mathbf{K}} = \mathbf{K}$. Then $(\bar{\mathbf{K}}, \mathbf{K})$ satisfies (1). Lemma 15.2 gives $\mathbf{L}$ satisfying (10).    ■

We are finally ready to prove Part (b) of the Main Theorem:

THEOREM 15.4: *Let* $\mathbf{K}$ *be a field structure,* $\mathbf{G}$ *a projective group structure, and* $\kappa$: $\mathbf{G} \to \mathrm{Gal}(\mathbf{K})$ *a rigid epimorphism. Then there exists a proper Henselian field-valuation structure* $\mathbf{L} = (L, X, L_x, w_x)_{x \in X}$ *and an isomorphism* $\psi$: $\mathbf{G} \to \mathrm{Gal}(\mathbf{L})$ *with* $L$ *perfect having these properties:*

(11a) $\mathbf{K} \subseteq \mathbf{L}$ *and* res $\circ \psi = \kappa$.

(11b) $w_x$ *is trivial on* $K$, $x \in X$.

(11c) $\mathbf{L}$ *satisfies the block approximation condition.*

*Proof:*  Replace $K$ by $K_x$ and $K_{\mathrm{ins}}$ by $(K_x)_{\mathrm{ins}}$, if necessary, to assume $K$ is perfect. Proposition 14.4 gives a proper field structure $\mathbf{K}'$ which extends $\mathbf{K}$ and an isomorphism $\kappa'$: $\mathbf{G} \to \mathrm{Gal}(\mathbf{K}')$ with $\mathrm{res}_{\bar{K}'/\bar{K}} \circ \kappa' = \kappa$. Proposition 15.3 extends $\mathbf{K}'$ to a proper Henselian field-valuation structure $\mathbf{L} = (L, X, L_x, w_x)_{x \in X}$ that satisfies the block approximation theorem such that $L$ is perfect and $\mathrm{res}_{\bar{L}/\bar{K}'}$: $\mathrm{Gal}(\mathbf{L}) \to \mathrm{Gal}(\mathbf{K}')$ is an isomorphism. Thus, there is an isomorphism $\psi$: $\mathbf{G} \to \mathrm{Gal}(\mathbf{L})$ with $\mathrm{res}_{L_s/K'_s} \circ \psi = \varphi'$. This establishes (11a), (11b), and (11c).    ■

An easy consequence of Theorem 15.4 is the realization of profinite products of finitely many absolute Galois groups as an absolute Galois group. Of course, one may get away with a much reduced machinery than the one we have developed here. See [Ers], [Koe], or [HJK].

THEOREM 15.5: *For each $i$ in a set $I_0$ let $K_i$ be a field which is not separably closed. Then there is a proper Henselian field-valuation structure* $\mathbf{L} = (L, X, L_x, w_x)_{x \in X}$ *with* $\mathrm{char}(L) = 0$ *satisfying the block approximation condition and* $G(L) \cong \coprod_{i \in I_0} \mathrm{Gal}(K_i)$.

*Proof:* Choose a set $T$ of cardinality at least the transcendence degree of $K_i$ over its prime field for all $i \in I_0$. By Proposition 7.5, $\mathbb{Q}(T)$ has an algebraic extension $K_i'$ with $\mathrm{Gal}(K_i) \cong \mathrm{Gal}(K_i')$. Let $K = \bigcap_{i \in I_0} K_i'$. Then replace $K_i$ by $K_i'$, if necessary, to assume all $K_i$ are algebraic extension of $K$ and $\mathrm{Gal}(K) = \langle \mathrm{Gal}(K_i) \mid i \in I_0 \rangle$.

For each $i \in I_0$ let $G_i$ be an isomorphic copy of $\mathrm{Gal}(K_i)$ and $\kappa_i \colon G_i \to \mathrm{Gal}(K_i)$ an isomorphism. Example 4.7 constructs a proper projective group structure $\mathbf{G} = (G, X, G_x)_{x \in X}$ with $G = \coprod_{i \in I_0} \mathrm{Gal}(K_i)$. Let $\kappa \colon G \to \mathrm{Gal}(K)$ be the epimorphism whose restriction to $G_i$ is $\kappa_i$. By Example 2.5, $\mathbf{G}/\mathrm{Ker}(\kappa)$ is a group structure and the quotient map $\mathbf{G} \to \mathbf{G}/\mathrm{Ker}(\kappa)$ is a cover. Thus, there is a field structure $\mathbf{K} = (K, Y, K_y)$ and $\kappa$ extend to a cover $\kappa \colon \mathbf{G} \to \mathrm{Gal}(\mathbf{K})$. Theorem 15.4 gives the desired field-valuation structure $\mathbf{L}$. ∎

24 September,  2002